\documentclass{amsart}

\usepackage{amssymb,amsmath,latexsym,amsthm,amscd}
\usepackage[english]{babel}

\newtheorem{theorem}{Theorem}
\newtheorem{lemma}[theorem]{Lemma}
\newtheorem{proposition}[theorem]{Proposition}
\newtheorem{corollary}[theorem]{Corollary}

\newtheorem{definition}[theorem]{Definition}
\newtheorem{remark}[theorem]{Remark}

\numberwithin{theorem}{section}
\numberwithin{equation}{section}

\newcommand{\F}{\mathbb{F}}
\newcommand{\Q}{\mathbb{Q}}

\newcommand{\Z}{\mathbb{Z}}

\newcommand{\Gal}{\mathrm{Gal}}
\newcommand{\G}{\mathrm{G}}
\newcommand{\M}{\mathrm{M}}
\newcommand{\GL}{\mathrm{GL}}
\newcommand{\PGL}{\mathrm{PGL}}
\newcommand{\SL}{\mathrm{SL}}

\newcommand{\Aut}{\mathrm{Aut}}

\newcommand{\Ker}{\mathrm{Ker}\, }
\newcommand{\Norm}{\mathrm{Norm}}
\newcommand{\ord}{\mathrm{ord}\, }

\newcommand{\tr}{\mathrm{Tr}\, }
\newcommand{\trace}{\mathrm{Tr}}
\newcommand{\Spec}{\mathrm{Spec}}
\newcommand{\z}{\mathrm{Z}}
\newcommand{\W}{\mathrm{W}}

\newcommand{\resp}{resp.\ }

\newcommand{\fix}{\mathrm{Fix}}
\newcommand{\Conj}{\mathrm{Conj}}

\begin{document}

\title[On uniform lower bound of the Galois images]
{On uniform lower bound of 
the Galois images associated to elliptic curves}


\author[Keisuke {\sc Arai}]{{\sc Keisuke} ARAI}
\address{Keisuke {\sc Arai}\\
Graduate School of Mathematical Sciences\\
The University of Tokyo\\
Tokyo 153-8914, Japan}
\email{araik@ms.u-tokyo.ac.jp}

\maketitle

\begin{abstract}
Let $p$ be a prime and $K$ be a number field. 
Let $\rho_{E,p}:\G_K \longrightarrow \Aut(T_p E)\cong \GL_2(\Z_p)$ 
be the Galois representation 
given by the Galois action on the $p$-adic Tate module of an 
elliptic curve $E$ over $K$.
Serre showed that the image of $\rho_{E,p}$ is open 
if $E$ has no complex multiplication. 
For an elliptic curve $E$ over $K$ whose $j$-invariant does not appear 
in an exceptional finite set, 
we give an explicit uniform lower bound 
of the size of the image of $\rho_{E,p}$. 
\end{abstract}

\bigskip


\section{Introduction}
\label{sec:intro}

Let $k$ be a field of characteristic $0$, 
and let $\G_k=\Gal(\overline{k}/k)$ be the absolute Galois 
group of $k$ where $\overline{k}$ is an algebraic closure 
of $k$. 
Let $p$ be a prime number. 
For an elliptic curve $E$ over $k$, let $T_p E$ be 
the $p$-adic Tate module of $E$, and let 
$$\rho_{E,p}:\G_k \longrightarrow \Aut(T_p E)\cong \GL_2(\Z_p)$$ 
be the $p$-adic representation determined by the action of $\G_k$ on 
$T_p E$. 
By a number field we mean a finite extension of $\Q$.

We recall a famous theorem proved by Serre:

\begin{theorem}
\label{thSe}

([17, IV-11])
Let $K$ be a number field, $E$ be an elliptic curve over $K$ with no 
complex multiplication and $p$ be a prime. 
Then the representation 
$\rho_{E,p}$ 
has an open image i.e. there exists a positive integer $n$ depending 
on $K$, $E$ and $p$ such that 
$$\rho_{E,p}(\G_K)\supseteq 1+p^n \M_2(\Z_p).$$

\end{theorem}

In this paper, we show that 
there exists a uniform 
bound of such $n$ if we let $E$ vary for fixed $K$ and $p$.

\begin{theorem}
\label{bound1}

Let $K$ be a number field and $p$ be a prime. 
Then there exists a positive integer $n$ 
depending on $K$ and $p$ such that for any elliptic curve $E$ 
over $K$ with no complex multiplication, we have 
$$\rho_{E,p}(\G_{K})\supseteq 1+p^n \M_2(\Z_p).$$

\end{theorem}

We will deduce Theorem \ref{bound1} from the following more 
precise result of this 
paper in Proposition \ref{deduction}.
We use the following conventions:
\begin{align*}
1+p^0\Z_p&:=\Z_p^{\times},\\
1+p^0\M_2(\Z_p)&:=\GL_2(\Z_p),\\
1+p^0\M_2(\Z/p\Z)&:=\GL_2(\Z/p\Z).
\end{align*}

\begin{theorem}
\label{bound2}

For a prime $p$, 
there exists an integer $n\geq 0$ satisfying the following 
condition $(\textbf{C})_p$. 

\noindent
$(\textbf{C})_p$:
Let $K$ be a number field. 
Then there exists a finite subset $\Sigma\subseteq K$ depending on 
$p$ such that 
for an elliptic curve $E$ over $K$, 
the condition $j(E)\not\in \Sigma$ implies 
$$\rho_{E,p}(\G_{K})\supseteq 
(1+p^n \M_2(\Z_p))^{\det=1}.$$

Let $n(p)\geq 0$ be the minimum integer $n$ satisfying $(\textbf{C})_p$. 
Then we have 
\begin{align*}
n(p)&\leq \begin{cases}
  0 & \text{if $p\geq 23$},\\
  1 & \text{if $p=19$, $17$, $13$, $11$},\\
  2  & \text{if $p=7$},\\
  3  & \text{if $p=5$},\\
  5  & \text{if $p=3$},\\
  11  & \text{if $p=2$}.
  \end{cases}
\end{align*}

\end{theorem}

The estimate in Theorem 1.3 is best possible for $p\geq 7$. 

\begin{corollary}
\label{n(p)}

We have 
\begin{align*}
n(p)&= \begin{cases}
  0 & \text{if $p\geq 23$},\\
  1 & \text{if $p=19$, $17$, $13$, $11$},\\
  2  & \text{if $p=7$}.
  \end{cases}
\end{align*}

\end{corollary}

\begin{proof}

Since $g(X_0(19))=1$, the integer $n(19)$ cannot be $0$. 
We also have 
$g(X_0(17))=1$, $g(X_0(13))=0$, $g(X_0(11))=1$ and 
$g(X_0(7^2))=1$. 
Thus we get the result in a similar way. 

\end{proof}

\begin{remark}
\label{known-result}

If $p\geq 23$ and 
\begin{equation*}
  \begin{cases}
    K\nsupseteq \text{(the quadratic subfield of $\Q(\zeta_p)$)}
    & \text{when $p\not\equiv\pm 3\mod{8}$},\\
    \text{there exists an inclusion $K\hookrightarrow \Q_p$}
    & \text{when $p\equiv\pm 3\mod{8}$},
  \end{cases}
\end{equation*}
then the result of Theorem 1.3 follows from 
[3, Theorem 7; 21, Proposition 1.40, 1.43; 
9, p.116-118].

If $p=17$, $19$ and $K=\Q$, the result of Theorem 1.3 follows from 
[3, Theorem 7; 8, Theorem (4.1); 
9, p.116-118].

If $p=2$ and $K=\Q$, the Galois image is studied in [14], 
assuming that an elliptic curve has no complex multiplication, 
and all its $2$-torsion points are rational. 

\end{remark}

\begin{remark}

There are many other studies of the Galois images associated to 
elliptic curves over number fields or of rational points on 
modular curves in [2,4,6,7,10-13,15,16,18,20]. 
Several questions related to the subject 
of this paper are raised in [19, p.187].

\end{remark}

The contents of this paper are as follows:

In Section 2, we deduce Theorem 1.2 from Theorem 1.3 by studying 
the determinants. 

In Section 3, we regard elliptic curves as rational points 
on modular curves, and reduce Theorem \ref{bound2} to 
a genus estimate. 
Replacing $K$ by its finite extension, we may assume that 
$K$ contains a primitive $p^{n(p)+1}$-st root $\zeta_{p^{n(p)+1}}$ 
of unity. 
Suppose an elliptic curve $E/K$ does not satisfy 
$\rho_{E,p}(\G_{K})\supseteq 
(1+p^{n(p)} \M_2(\Z_p))^{\det=1}$. 
Then we have 
$\rho_{E,p}(\G_{K})\mod{p^{n(p)+1}}\subseteq H$ 
for some subgroup 
$H\subseteq \SL_2(\Z/p^{n(p)+1}\Z)$ 
satisfying 
$H\nsupseteq (1+p^{n(p)}\M_2(\Z/p^{n(p)+1}\Z))^{\det=1}$. 
Thus $E/K$ determines a rational point on the modular curve $X_H$ 
corresponding to $H$. 
If the genus $g(X_H)$ of $X_H$ is greater than or equal to $2$, 
we conclude that $X_H$ has only finitely 
many rational points by Mordell's conjecture ([3, Theorem 7]). 
Since there are only finitely many subgroups $H$ as above, 
the number of the $j$-invariants of $E/K$ 
not satisfying 
$\rho_{E,p}(\G_{K})\supseteq (1+p^{n(p)} \M_2(\Z_p))^{\det=1}$
is finite. 
Thus Theorem 1.3 will follow. 

In Section 4 - 7, we prove $g(X_H)\geq 2$. 
In section 4, we prepare for estimating $g(X_H)$. 
Put 
$G=\SL_2(\Z/p^{n(p)+1}\Z),\;
\sigma =\left(\begin{matrix} 0&1\\ -1&0 \end{matrix}\right),\;
\tau =\left(\begin{matrix} 1&1\\ -1&0 \end{matrix}\right),\;
u=\left(\begin{matrix} 1&1\\ 0&1 \end{matrix}\right)$. 
For $\alpha\in G$, 
let $\Conj(\alpha)$ be the conjugacy class containing $\alpha$. 
If $H\ni -1$, we have 
$$g(X_H)=1+\frac{1}{12}[G:H]\left(1
  -3\frac{\sharp H\cap\Conj(\sigma)}{\sharp \Conj(\sigma)}
  -4\frac{\sharp H\cap\Conj(\tau)}{\sharp \Conj(\tau)}
  -6\frac{\sharp \langle u \rangle\backslash G/H}{[G:H]}\right).$$ 
We calculate the number of elements conjugate to 
$\sigma,\tau,u$ contained in maximal subgroups of $\SL_2(\Z/p\Z)$.

In Section 5, we calculate the number of elements conjugate to 
$\sigma,\tau,u$ contained in $\SL_2(\Z/p^n\Z)$, and 
study the fiber of the mod $p^m$ map 
$\SL_2(\Z/p^n\Z)\cap\Conj(\alpha)\longrightarrow 
\SL_2(\Z/p^m\Z)\cap\Conj(\alpha)$, 
where $1\leq m\leq n$ and $\alpha=\sigma,\tau,u$. 
For integers $1\leq m\leq n$, let 
$f_{n,m}:\SL_2(\Z/p^n\Z)\longrightarrow \SL_2(\Z/p^m\Z)$ 
be the $\bmod{\ p^m}$ map. 
If $m<n\leq 2m$ and $\alpha=\sigma,\tau,u$, we see that 
$\alpha^{-1}(f_{n,m}^{-1}(\alpha)\cap\Conj(\alpha))$ 
is a subgroup of 
$(1+p^m\M_2(\Z/p^n\Z))^{\det=1}\cong \M_2(\Z/p^{n-m}\Z)^{\trace=0}$ 
and isomorphic to $(\Z/p^{n-m}\Z)^2$.

In Section 6, we control the number of elements conjugate to 
$\sigma,\tau,u$ contained in $H$, by 
combining the result of Section 5 with the property 
$H\nsupseteq (1+p^{n(p)}\M_2(\Z/p^{n(p)+1}\Z))^{\det=1}$.

In Section 7, we prove $g(X_H)\geq 2$ by using the results of 
Section 4, 5, 6.

I would like to thank my supervisor Professor Takeshi Saito for 
helpful advice and warm encouragement. 
This work was partly supported by 21st Century COE Program 
in The University of Tokyo, A Base for New Developments of 
Mathematics into Science and Technology.


\section{Deduction of Theorem \ref{bound1}}

In order to deduce Theorem \ref{bound1} from Theorem \ref{bound2}, 
we need some facts in 
group theory.

\begin{lemma}
\label{1+p0}

Let $p$ be a prime and $H$ be a closed subgroup of $\GL_2(\Z_p)$. 
Then $H$ contains $\SL_2(\Z_p)$ if and only if 
$H\bmod{p^2}$ contains $\SL_2(\Z/p^2\Z)$. 
In particular, if a subgroup $H'\subseteq \SL_2(\Z/p^n\Z)$ for $n\geq 3$ 
maps surjectively $\bmod{\ p^2}$ onto $\SL_2(\Z/p^2\Z)$, then 
$H'=\SL_2(\Z/p^n\Z)$. 

Assume $p\geq 5$. Then 
$H$ contains $\SL_2(\Z_p)$ if and only if 
$H\bmod{p}$ contains $\SL_2(\Z/p\Z)$. 
In particular, if a subgroup $H'\subseteq \SL_2(\Z/p^n\Z)$ for $n\geq 2$ 
maps surjectively $\bmod{\ p}$ onto $\SL_2(\Z/p\Z)$, then 
$H'=\SL_2(\Z/p^n\Z)$. 

\end{lemma}

\begin{proof}

See [17, IV-23]. 

\end{proof}

\begin{lemma}
\label{mod-p^n+1}

Let $n\geq 1$ be an integer and let $p$  be a prime. 
If $p=2$, assume $n\geq 2$. 
Let $H$ be a closed subgroup of $\GL_2(\Z_p)$.
Then $H$ contains $1+p^n \M_2(\Z_p)$ 
($\resp (1+p^n \M_2(\Z_p))^{\det=1}$)
if and only if 
$H\bmod{p^{n+1}}$ contains $1+p^n \M_2(\Z/p^{n+1}\Z)$ 
($\resp (1+p^n \M_2(\Z/p^{n+1}\Z))^{\det=1}$). 

\end{lemma}

\begin{proof}
The ``only if'' part is trivial. 

Assume $H\bmod{p^{n+1}}$ contains $1+p^n \M_2(\Z/p^{n+1}\Z)$. We will show 
that $H$ contains $1+p^n \M_2(\Z_p)$.
For each $s\in 1+p^n \M_2(\Z_p)$ and for any $a\geq 1$, we will find 
an element 
$x_a\in H$ with $x_a\equiv s\bmod{p^{n+a}}$. Then 
$s=(x_a\bmod{p^{n+a}})_{a\geq 1}\in H$, as desired.
We show the existence of $x_a$'s by induction on $a$.
For $a=1$, it follows from the assumption.
Assume it is true for $a\, (\geq 1)$ i.e. 
$H\bmod{p^{n+a}}\supseteq 1+p^n \M_2(\Z/p^{n+a}\Z)$.
We first show that for $s\in 1+p^n \M_2(\Z_p)$
with $s\equiv 1\bmod{p^{n+a}}$, there exists $x\in H$ with 
$x\equiv s\bmod{p^{n+a+1}}$. We write $s=1+p^{n+a}u$ with some 
$u\in \M_2(\Z_p)$. By the hypothesis of induction, we find 
$y:=1+p^{n+a-1}u+p^{n+a}v\in H$ for some $v\in \M_2(\Z_p)$. 
Put $x:=y^p$. Then $x\in H$ and 
$x=1+p^{n+a}u+p^{n+a+1}\cdots \equiv s\bmod{p^{n+a+1}}$. 
For a general $s\in 1+p^n \M_2(\Z_p)$, the induction hypothesis shows that 
$(s\bmod{p^{n+a}})^{-1}\in H\bmod{p^{n+a}}$. 
Hence we have 
$st\equiv 1\bmod{p^{n+a}}$ for some $t\in H$.
By the above argument, we find $x\in H$ with $x\equiv st\bmod{p^{n+a+1}}$.
Therefore $s\equiv t^{-1}x\bmod{p^{n+a+1}}$ and $t^{-1}x\in H$.

For the ``$\det=1$'' version, it follows from the same argument as in 
the previous lemma ([17, IV-23]). 

\end{proof}

\begin{lemma}
\label{index2}

Let $n\geq 1$ be an integer. 
Let $H$ be a subgroup of $\GL_2(\Z_p)$ containing $1+p^n\M_2(\Z_p)$, and 
$H'$ be a closed subgroup of $\GL_2(\Z_p)$ which is a 
subgroup of $H$ of index $2$. If $p\geq 3$, then 
$H'\supseteq 1+p^n\M_2(\Z_p)$; if $p=2$ and $n\geq 2$, 
$H'\supseteq 1+p^{n+1}\M_2(\Z_p)$. 

\end{lemma}

\begin{proof}

For any $X\in \M_2(\Z_p)$, we have $H\ni 1+p^n X$ and 
$H'\ni (1+p^n X)^2=1+2p^n X+p^{2n}X^2$. This implies 
$H'\bmod{p^{n+1}}\supseteq 1+p^n\M_2(\Z/p^{n+1}\Z)$ if $p\geq 3$ ($\resp
 H'\bmod{p^{n+2}}\supseteq 1+p^{n+1}\M_2(\Z/p^{n+2}\Z)$ 
if $p=2$ and $n\geq 2$). 
Since $H'$ is closed, applying the previous lemma, we get the result. 

\end{proof}

\begin{corollary}
\label{H,-1}

Let $H\subseteq \GL_2(\Z_p)$ be a closed subgroup. 
If $\langle H,\,-1\rangle\supseteq 1+p^n \M_2(\Z_p)$
for some integer $n\geq 1$, then 
\begin{equation*}
H\supseteq \begin{cases}
  1+p^n \M_2(\Z_p) & \text{if $p\geq 3$},\\
  1+p^{n+1} \M_2(\Z_p) & \text{if $p=2$ and $n\geq 2$}.
  \end{cases}
\end{equation*}

\end{corollary}

\begin{lemma}
\label{det1}

Let $n\geq 1$ be an integer and 
assume $p\geq 3$. Let $H$ be a subgroup of $\GL_2(\Z/p^{n+1}\Z)$.
If $\det(H)=(\Z/p^{n+1}\Z)^{\times}$, 
then 
$$\det|_{H\cap (1+p^n \M_2(\Z/p^{n+1}\Z))}:H\cap (1+p^n \M_2(\Z/p^{n+1}\Z))
\longrightarrow 1+p^n\Z/p^{n+1}\Z$$ is surjective.

\end{lemma}

\begin{proof}

We first show the surjectivity when $n=1$. 
Let $x\in H$ be an element such that $\det x$ is a generator of 
$(\Z/p^2\Z)^{\times}$. Put $\overline{x}:=x\bmod{p}$ and 
$a:=\ord \overline{x}$. We first show $(p,a)=1$.
Since $\det \overline{x}$ generates $\F_p^{\times}$, $\overline{x}$
is not a scalar. We consider the commutative subalgebra
$\F_p [\overline{x}]$ of $\M_2(\F_p)$. If we put 
$t:=\tr \overline{x}$ and $d:=\det \overline{x}$, then 
$\overline{x}^2-t\overline{x}+dE=0$ and 
$\F_p [\overline{x}]\cong \F_p[X]/(X^2-tX+d)$.
$X^2-tX+d$ does not have multiple roots because $d$ generates 
$\F_p^{\times}$. Hence we have 
$\F_p [\overline{x}]\cong \F_p \times\F_p$ or $\F_{p^2}$. 
Therefore, $\overline{x}\in \F_p [\overline{x}]^{\times}\cong 
\F_p^{\times} \times\F_p^{\times}$ or $\F_{p^2}^{\times}$ and 
$(p,a)=1$.

Since $\overline{x}^a=1$, $x^a\in 1+p \M_2(\Z/p^2\Z)$.
Therefore, $\ord x^a=1$ or $p$, and 
$\ord (\det x)^a=\ord (\det x^a)=1$ or $p$.
Thus, we have $\ord (\det x)^a=\ord (\det x^a)=p$ because $p|\ord (\det x)$
and $(p,a)=1$. Hence, $\det x^a$ is a generator of $1+p\Z/p^2\Z$.

Next we show the surjectivity for general $n\geq 1$. 
Since $\det |_{H\bmod{p^2}}=(\det |_H)\bmod{p^2}$ is surjective, there 
exists an element $y\in (H\bmod{p^2})\cap (1+p\M_2(\Z/p^2\Z))$ such that 
$\det y$ is a generator of $1+p\Z/p^2\Z$. We choose a lift 
$\tilde{y}\in H$ of $y$. Then $\tilde{y}\in 1+p\M_2(\Z/p^{n+1}\Z)$. 
Taking the $p^{n-1}$-st power, we have 
$\tilde{y}^{p^{n-1}}\in H\cap (1+p^n\M_2(\Z/p^{n+1}\Z))$. 
Since $\det y=\det \tilde{y}\bmod{p^2}$ generates $1+p\Z/p^2\Z$, 
we can write 
$\det \tilde{y}=1+pz$ for some $z\in (\Z/p^{n+1}\Z)^{\times}$. 
Therefore $\det \tilde{y}^{p^{n-1}}=1+p^n z$ generates $1+p^n\Z/p^{n+1}\Z$. 

\end{proof}

\begin{lemma}
\label{det2}

Take two integers $n>r\geq 1$. 
If $p=2$, assume $r\geq 2$.  
Let $H$ be a subgroup of $\GL_2(\Z/p^{n+1}\Z)$.
If $\det(H)$ contains $1+p^r\Z/p^{n+1}\Z$ and if $H$ contains 
$(1+p^{n-r} \M_2(\Z/p^{n+1}\Z))^{\det =1}$, 
then 
$$\det|_{H\cap (1+p^n \M_2(\Z/p^{n+1}\Z))}:H\cap (1+p^n \M_2(\Z/p^{n+1}\Z))
\longrightarrow 1+p^n\Z/p^{n+1}\Z$$ is surjective. 
In particular, $H$ contains 
$1+p^n \M_2(\Z/p^{n+1}\Z)$. 

\end{lemma}

\begin{proof}

Let $x\in H$ be an element such that 
$\det x$ is a generator of $1+p^r\Z/p^{n+1}\Z$. 
Put $y:=x^{(p-1)^2(p+1)}$. 
We see that $y$ is also a generator 
of $1+p^r\Z/p^{n+1}\Z$ 
and $y$ belongs to a $p$-Sylow subgroup of 
$\GL_2(\Z/p^{n+1}\Z)$, which is conjugate to 
$\left\{
\left(\begin{matrix} 1+p*&*\\ p*&1+p* \end{matrix}\right)\right\}$. 
We can write 
$y=1+pz+\epsilon$, where $\epsilon$ is conjugate to 
$\left(\begin{matrix} 0&*\\ 0&0 \end{matrix}\right)$. 
Put $w:=y^{p^{n-r}}$. 
We see that $w$ belongs to $1+p^{n-r}\M_2(\Z/p^{n+1}\Z)$ and 
$\det w$ is a generator of $1+p^n\Z/p^{n+1}\Z$. 
Multiplying $w$ by elements of $(1+p^{n-r} \M_2(\Z/p^{n+1}\Z))^{\det =1}$, 
we get an element $h\in H\cap (1+p^n \M_2(\Z/p^{n+1}\Z))$ such that 
$\det h$ generates $1+p^n\Z/p^{n+1}\Z$. 

\end{proof}

\begin{corollary}
\label{det3}

Let $H\subseteq \GL_2(\Z_p)$ be a closed subgroup 
and $n$, $r\geq 0$ be integers. 
Assume $r\geq 2$ if $p=2$. 
If 
$H\supseteq (1+p^n\M_2(\Z_p))^{\det=1}$ and if 
$\det(H)\supseteq 1+p^r\Z_p$, then 
$H\supseteq 1+p^{n+r}\M_2(\Z_p)$. 

\end{corollary}

\begin{proof}

By Lemma \ref{mod-p^n+1}, 
it suffices to show the result $\bmod{\ p^{n+r+1}}$. 
If $n=0$, it is trivial. 
Assume $n\geq 1$. 
If $r=0$, it follows from Lemma \ref{det1}. 
If $r\geq 1$, it follows from Lemma \ref{det2}. 

\end{proof}

As a consequence of Theorem \ref{bound2}, we get the following. 

\begin{theorem}
\label{bound3}

Let $K$ be a number field, $p$ be a prime. 
Let $n(p)$, 
$\Sigma$ be as in Theorem \ref{bound2}.
Suppose that the image of the $p$-adic cyclotomic character 
$\chi_p:\G_K\longrightarrow \Z_p^{\times}$ 
contains $1+p^r\Z_p$ with $r\geq 0$ an integer. 
Assume $r\geq 2$ if $p=2$. 
Put 
\begin{equation*}
n=r+n(p).
\end{equation*}
Then for an elliptic curve $E$ over $K$, the condition $j(E)\not\in \Sigma$
implies 
$\rho_{E,p}(\G_{K})\supseteq 
1+p^{n} \M_2(\Z_p)$.

\end{theorem}

\begin{proof}

It follows from Theorem \ref{bound2} and Corollary \ref{det3}. 

\end{proof}

We show that there exists a lower bound of the images of 
$\rho_{E,p}$ if we take $E$'s having only finitely many $j$-invariants. 

\begin{lemma}
\label{j-inv}

Let $K$ be a number field. 
Fix an element $j\in K$. Assume that an elliptic curve $E$ 
over $K$ with $j$-invariant 
$j$ has no complex multiplication. 
Take a prime $p$. 
Then there exists a positive integer $n$ 
depending on $p$ and $j$ 
such that for any elliptic curve $E$ over $K$ with $j$-invariant $j$,
we have $\rho_{E,p}(\G_K)\supseteq 1+p^n \M_2(\Z_p)$.

\end{lemma}

\begin{proof}

Take an $E/K$ with $j$-invariant $j$. By Theorem \ref{thSe}, we have 
$\rho_{E,p}(\G_K)\supseteq 1+p^{n_0} \M_2(\Z_p)$ for some $n_0$. 
Let $E'/K$ have $j$-invariant $j$. 
Since $E$ has no complex multiplication, we have 
$j\ne 0,1728$. 
Hence there exists a quadratic extension $L$ of $K$ 
satisfying $E\otimes_K L\cong E'\otimes_K L$ 
(see [22, p.308]). 
Therefore 
$\rho_{E,p}(\G_L)$ is conjugate to $\rho_{E',p}(\G_L)$. 
Since $[\rho_{E,p}(\G_K):\rho_{E,p}(\G_L)]$ is 
$1$ or $2$, applying Lemma \ref{index2}, we get the result. 

\end{proof}

Now we can deduce Theorem \ref{bound1} from Theorem \ref{bound2} 
and some facts in group theory.

\begin{proposition}
\label{deduction}

Theorem \ref{bound2} and Lemma \ref{j-inv} imply Theorem \ref{bound1}. 

\end{proposition}

\begin{proof}

Theorem \ref{bound2} implies Theorem \ref{bound3}. 
Since $\Sigma$ is a finite set, elliptic curves $E$ over $K$ with 
$j(E)\in \Sigma$ have bounded Galois images by Lemma \ref{j-inv}. 

\end{proof}


\section{Modular curves}

We regard an elliptic curve with a specific Galois image as a 
rational point on a certain modular curve, and reduce Theorem \ref{bound2} 
to a genus estimate. 

Now we give a brief review of modular curves. 
For more details, see [1,5]. 
Let $N$ be a positive integer and 
$k$ be a field of characteristic $0$. 
For an elliptic curve $E$ over $k$ and an integer $N\geq 1$, let  
$E[N]=\Ker ([N]:E\longrightarrow E)$ be the kernel of multiplication 
by $N$ on $E$, and let 
$\overline{\rho}_{E,N}:\G_k\longrightarrow
\Aut(E[N](\overline{k}))\cong \GL_2(\Z/N\Z)$ 
be the $\bmod{\; N}$ representation 
determined by the action of $G_k$ on $E[N](\overline{k})$. 
A level $N$-structure on $E$ is 
an isomorphism 
$$\gamma:(\Z/N\Z)^2\longrightarrow E[N].$$ 
Let $Y(N)\longrightarrow \Spec (\Q (\zeta_N))$ 
be the moduli of 
elliptic curves with level $N$-structure. 
We have a right action of 
$G=\SL_2(\Z/N\Z)$ on $Y(N)$ over $\Q (\zeta_N)$ : 
$$[E,\gamma]\mapsto [E,\gamma\circ h]$$ 
where 
$E$ is an elliptic curve over $k$, $\gamma$ a level $N$-structure on $E$ 
and $h\in G$.

For a subgroup $H\subseteq G$, 
put $Y_H$ to be the quotient $Y(N)/H$. The quotient 
$Y_H\longrightarrow \Spec (\Q (\zeta_N))$ is 
an affine smooth curve. 
Let $E$ be an elliptic curve over $k$. 
Choose a basis 
$\langle \epsilon_1,\, \epsilon_2\rangle$ of $E[N](\overline{k})$. 
Then the pair $(E,\langle \epsilon_1,\, \epsilon_2\rangle)$ defines an 
element $P$ of $Y(N)(\overline{k})$. 
Let $Q\in Y_H(\overline{k})$ be the image of $P$ via the map 
$Y(N)(\overline{k})\longrightarrow Y_H(\overline{k})$
induced by the natural map 
$Y(N)\longrightarrow Y_H$. 
If $\overline{\rho}_{E,N}(\G_k)\subseteq H$ with respect to 
$\langle \epsilon_1,\, \epsilon_2\rangle$, then 
$Q$ lies in $Y_H(k)$. 
Note that if two subgroups $H,\,H'\subseteq G$ are 
$\GL_2(\Z/N\Z)$-conjugate, 
then $Y_H$ and $Y_{H'}$ are conjugate.

\begin{lemma}
\label{Y_H}

Let $k$ be a field of characteristic $0$. 
If $Y_H(k)$ is finite, then there exists a finite set 
$\Sigma\subseteq k$ 
satisfying the following condition:

For an elliptic curve $E$ over $k$, 
if a conjugate of 
$\overline{\rho}_{E,N}(\G_k)$ is contained in $H$, then 
$j(E)\in \Sigma$.

\end{lemma}

Let $X_H$ be the smooth compactification of $Y_H$. 
The following is the famous theorem known as Mordell's conjecture proved by 
Faltings. It shows that a curve $X$ over a number field 
has only finitely many rational points 
if its genus $g(X)$ is greater than or equal to $2$. 

\begin{theorem}
\label{Mordell}

([3, Theorem 7])
Let $K$ be a number field and $X/K$ be a proper smooth curve. 
If $g(X)\geq 2$, then $X(K)$ is finite. 

\end{theorem}


Now we compute the genus of $X_H$ explicitly. 
As in Section 1, put 
$$\sigma :=\left(\begin{matrix} 0&1\\ -1&0 \end{matrix}\right),\;
\tau :=\left(\begin{matrix} 1&1\\ -1&0 \end{matrix}\right),\;
u:=\left(\begin{matrix} 1&1\\ 0&1 \end{matrix}\right)
\in \SL_2(\Z).$$
We have $\ord \sigma=4$ and $\ord \tau=6$. 
Note that 
$\sigma^{-1}=\left(\begin{matrix} 0&-1\\ 1&0 \end{matrix}\right)$ 
and $\tau^{-1}=\left(\begin{matrix} 0&-1\\ 1&1 \end{matrix}\right)$. 
For $\alpha\in \SL_2(\Z)$, 
we also use the same letter to denote the reduction of $\alpha$. 
Put $$\fix_{\alpha}:=\{gH\in G/H|\alpha gH=gH\}$$
and 
$$g_H:=1+
  \frac{1}{12}\sharp G/H-\frac{1}{4}\sharp \fix_{\sigma}
  -\frac{1}{3}\sharp \fix_{\tau}
  -\frac{1}{2}\sharp I\backslash G/H.$$
Here $I$ is the cyclic subgroup of $G$ generated by 
$u=\left(\begin{matrix} 1&1\\ 0&1 \end{matrix}\right)$.

\begin{proposition}
\label{gX_H}

([21, Proposition 1.40])
Let $H$ be a subgroup of $G=\SL_2(\Z/N\Z)$. Assume that $H$ contains 
$-1$. Then the genus $g(X_H)$ of the modular curve $X_H$ is given by 
$$g(X_H)=g_H.$$

\end{proposition}

Let $p$ be a prime and consider subgroups of $\GL_2(\Z/p\Z)$. 
A Borel subgroup is a subgroup which is conjugate to 
$\left\{\left(\begin{matrix} *&*\\ 0&* \end{matrix}\right)\right\}$; 
the normalizer of a split Cartan subgroup is conjugate to 
$\left\{\left(\begin{matrix} *&0\\ 0&* \end{matrix}\right),
  \left(\begin{matrix} 0&*\\ *&0 \end{matrix}\right)\right\}$. 
When $p\geq 3$, 
the normalizer of a non-split Cartan subgroup is 
conjugate to
$\left\{\left(\begin{matrix} x&y\\ \lambda y&x \end{matrix}\right),
  \left(\begin{matrix} x&y\\ -\lambda y&-x \end{matrix}\right)
|(x,y)\in \F_p\times \F_p\setminus \{(0,0)\}\right\}$, 
where $\lambda\in \F_p^{\times}\setminus (\F_p^{\times})^2$ is a fixed 
element. 
Assume $p\geq 5$. 
The quotient group $\PGL_2(\Z/p\Z)$ of $\GL_2(\Z/p\Z)$ 
has a subgroup which is 
isomorphic to $S_4$; it has a subgroup which is isomorphic to 
$A_5$ if and only if $p\equiv 0,\pm 1\mod{5}$ 
([18, p.281]). 
Take a subgroup $H$ (of $\GL_2(\Z/p\Z)$) whose order is prime to $p$. 
We call $H$ an exceptional subgroup 
if it is the inverse image of a subgroup which is isomorphic to 
$A_4$, $S_4$ or $A_5$ by the natural surjection 
$\GL_2(\Z/p\Z)\longrightarrow \PGL_2(\Z/p\Z)$.

\begin{proposition}
\label{maxsgp}

([18, p.284])
Let $p\geq 3$ be a prime and 
$H$ be a subgroup of $\GL_2(\Z/p\Z)$. 
If $p$ divides the order of $H$, then 
$H$ contains $\SL_2(\Z/p\Z)$ or 
$H$ is contained in a Borel subgroup of $\GL_2(\Z/p\Z)$. 
If $p$ does not divide the order of $H$, then $H$ is contained in 
the normalizer of a (split or non-split)
Cartan subgroup of $\GL_2(\Z/p\Z)$ or an exceptional 
subgroup of $\GL_2(\Z/p\Z)$. 

\end{proposition}

Put 
\begin{align*}
B&:=\left\{\left(\begin{matrix} *&*\\ 0&* \end{matrix}\right)\right\}
\cap \SL_2(\Z/p\Z),\\
C&:=\left\{\left(\begin{matrix} *&0\\ 0&* \end{matrix}\right),
  \left(\begin{matrix} 0&*\\ *&0 \end{matrix}\right)\right\}
\cap \SL_2(\Z/p\Z),\\
D&:=\left\{\left(\begin{matrix} x&y\\ \lambda y&x \end{matrix}\right),
  \left(\begin{matrix} x&y\\ -\lambda y&-x \end{matrix}\right)\right\}
\cap \SL_2(\Z/p\Z),\\ 
E&:=\text{(an exceptional subgroup)}\cap \SL_2(\Z/p\Z).
\end{align*}
From now on, we use the letter $E$ to denote this subgroup, 
not meaning an elliptic curve. 
The genera of the modular curves corresponding to $B$, $C$, $D$ 
are known as follows. 

\begin{proposition}
\label{g_BCD}

([21, Proposition 1.40, 1.43; 9, p.117])
Let $N=p\geq 5$ be a prime. We have 
\begin{align*}
g_B&=\frac{1}{12}(p-6-3\left(\frac{-1}{p}\right)-4\left(\frac{-3}{p}\right));
\\
g_C&=\frac{1}{24}(p^2-8p+11-4\left(\frac{-3}{p}\right));
\\
g_D&=\frac{1}{24}(p^2-10p+23+6\left(\frac{-1}{p}\right)
+4\left(\frac{-3}{p}\right)).
\end{align*}
We have 
$g_B\geq 2$ if and only if $p\geq 23$; 
we have 
$g_C\geq 2$ if and only if $p\geq 11$; 
we have 
$g_D\geq 2$ if and only if $p\geq 13$. 

\end{proposition}

\begin{remark}

We can also calculate these genera by using Lemma \ref{bcde-conj}.

\end{remark}

Put 
\begin{equation*}
\delta_H:=[G:H]^{-1}(
  \sharp G/H-3\sharp \fix_{\sigma}
  -4\sharp \fix_{\tau}
  -6\sharp I\backslash G/H)
\end{equation*}
so that 
$$g_H=1+\frac{1}{12}[G:H]\delta_H.$$ 
As $g_H$
is an integer, 
we have $g_H\geq 2$ if and only if $\delta_H>0$.

\begin{definition}
\label{slim}

Let $n\geq 1$ be an integer and let 
$H\subseteq \SL_2(\Z/p^n\Z)$ be a subgroup. 
We call $H$ a slim subgroup if 
$$H\nsupseteq (1+p^{n-1}\M_2(\Z/p^n\Z))^{\det=1}.$$

\end{definition}

In this definition, notice that if $n=1$, then a slim subgroup is 
just a proper subgroup.

In order to prove Theorem \ref{bound2}, it suffices to estimate $\delta_H$ 
for any slim subgroup $H$. 

\begin{proposition}
\label{delta}

If $\delta_H>0$ for any slim subgroup 
$H\subseteq \SL_2(\Z/p^{n'(p)+1}\Z)$, 
then Theorem \ref{bound2} holds. 
Here we put 
\begin{equation*}
n'(p):= \begin{cases}
  0  & \text{if $p\geq 23$},\\
  1  & \text{if $p=19,17,13,11$},\\
  2  & \text{if $p=7$},\\
  3  & \text{if $p=5$},\\
  5  & \text{if $p=3$},\\
  10  & \text{if $p=2$}.
  \end{cases}
\end{equation*}

\end{proposition}

\begin{proof}

Put 
\begin{equation*}
\xi:= \begin{cases}
  0  & \text{if $p\geq 3$},\\
  1  & \text{if $p=2$}.
  \end{cases}
\end{equation*}
Let $E$ be an elliptic curve over $K$ satisfying 
$\rho_{E,p}(\G_K)\nsupseteq(1+p^{n'(p)+\xi}\M_2(\Z_p))^{\det=1}$. 
We show that $j(E)\in K$ takes only finitely many values. 
Replacing $K$ by $K(\zeta_{p^{n'(p)+1}})$, we may assume 
$\overline{\rho}_{E,p^{n'(p)+1}}(\G_K)\subseteq 
\SL_2(\Z/p^{n'(p)+1}\Z)$. 
We may also assume that 
$\overline{\rho}_{E,p^{n'(p)+1}}(\G_K)$ 
is contained in 
a slim subgroup 
$H\subseteq \SL_2(\Z/p^{n'(p)+1}\Z)$ 
satisfying 
$H\ni -1$. 
To see this, we consider two cases ($n'(p)=0$ or $n'(p)\geq 1$). 
When $n'(p)=0$ (equivalently $p\geq 23$), we have 
$\overline{\rho}_{E,p}(\G_K)\subsetneq \SL_2(\Z/p\Z)$ 
by Lemma \ref{1+p0}, thus 
$\overline{\rho}_{E,p}(\G_K)\subseteq B, C, D, E$ 
by Proposition \ref{maxsgp}. 
But $B, C, D, E$ contains $-1$. 
When $n'(p)\geq 1$, Corollary \ref{H,-1} shows 
$\langle \rho_{E,p}(\G_K),-1\rangle
\nsupseteq(1+p^{n'(p)}\M_2(\Z_p))^{\det=1}$. 
Thus $H=\langle \rho_{E,p}(\G_K),-1\rangle\bmod{p^{n'(p)+1}}$ 
is a slim subgroup by Lemma \ref{mod-p^n+1}. 

By the hypothesis and Proposition \ref{gX_H}, we have 
$g(X_H)=g_H\geq 2$. 
By Theorem \ref{Mordell}, we see that $X_H(K)$ is finite, hence 
$Y_H(K)$ is also finite. 
Since there are only finitely many subgroups $H$ as above, 
Theorem \ref{bound2} follows from Lemma \ref{Y_H}. 

\end{proof}

We prove $\delta_H>0$ for any subgroup $H$ as in Proposition \ref{delta}. 
More explicitly, we prove the following theorem in Section 7.

\begin{theorem}
\label{main}

1. For a subgroup 
$H\subseteq \SL_2(\Z/p\Z)$, we have $\delta_H>0$ if one of the following 
conditions is satisfied. 
\begin{itemize}
  \item $H\subseteq B$ and $p\geq 23$
  \item $H\subseteq C$ and $p\geq 11$
  \item $H\subseteq D$ and $p\geq 13$
  \item $H\subseteq E$ and $p\geq 17$
\end{itemize}

2. For a slim subgroup 
$H\subseteq \SL_2(\Z/p^2\Z)$, 
we have $\delta_H>0$ if one of the following 
conditions is satisfied. 
\begin{itemize}
  \item $H\bmod{p}\subseteq B$ and $p\geq 11$
  \item $H\bmod{p}\subseteq D$ and $p\geq 11$
  \item $H\bmod{p}\subseteq E$ and $p\geq 11$
\end{itemize}

3. For a slim subgroup 
$H\subseteq \SL_2(\Z/p^3\Z)$, 
we have $\delta_H>0$ if one of the following 
conditions is satisfied. 
\begin{itemize}
  \item $H\bmod{p}\subseteq B$ and $p\geq 7$
  \item $H\bmod{p}\subseteq C$ and $p\geq 5$
  \item $H\bmod{p}\subseteq D$ and $p\geq 7$
  \item $H\bmod{p}\subseteq E$ and $p\geq 7$
\end{itemize}

4. For a slim subgroup 
$H\subseteq \SL_2(\Z/5^4\Z)$, 
we have $\delta_H>0$ if one of the following 
conditions is satisfied. 
\begin{itemize}
  \item $H\bmod{5}\subseteq B$
  \item $H\bmod{5}\subseteq D$
  \item $H\bmod{5}\subseteq E$
\end{itemize}

5. For a slim subgroup 
$H\subseteq \SL_2(\Z/3^6\Z)$, 
we have $\delta_H>0$ if one of the following 
conditions is satisfied. 
\begin{itemize}
  \item $H\bmod{3}\subseteq B$
  \item $H\bmod{3}\subseteq D$
  \item $H\bmod{3}\subseteq E$
  \item $H\bmod{3}=\SL_2(\Z/3\Z)$
\end{itemize}

6. For a slim subgroup 
$H\subseteq \SL_2(\Z/2^{10}\Z)$, 
we have $\delta_H>0$ if one of the following 
conditions is satisfied. 
\begin{itemize}
  \item $H\bmod{2}\subseteq (\text{subgroup of order $3$})$
  \item $H\bmod{2}=\SL_2(\Z/2\Z)$
\end{itemize}

7. For a slim subgroup 
$H\subseteq \SL_2(\Z/2^{11}\Z)$, 
we have $\delta_H>0$ if the following 
condition is satisfied. 
\begin{itemize}
  \item $H\bmod{2}\subseteq B$
\end{itemize}

\end{theorem}


\section{Calculation of conjugate elements in $\SL_2(\Z/p\Z)$}

We rewrite the value $\delta_H$ in terms of the number of 
conjugate elements.

\begin{lemma}
\label{fix-conj}

Let $G$ be a finite group and $H$ be a subgroup of $G$. Take an element $a$ 
of $G$. Then 
$$\sharp \fix_a /[G:H]=\sharp (H\cap\Conj(a))/\sharp \Conj(a).$$
Here we put $\Conj(a):=\{g^{-1}ag\in G|g\in G\}$. 

\end{lemma}

\begin{proof}
Put $\widetilde{\fix}_a:=\{g\in G|gH=agH\}$. Then $H$ acts on 
$\widetilde{\fix}_a$ 
by $g\mapsto gh$ ($g\in \widetilde{\fix}_a,\; h\in H$) and we have 
$\widetilde{\fix}_a/H=\fix_a$. 
Hence, we get $\sharp\widetilde{\fix}_a/\sharp H=\sharp\fix_a$. 
The natural surjection 
$\widetilde{\fix}_a \longrightarrow H\cap\Conj(a) : g\longmapsto g^{-1}ag$ 
induces a bijection 
$\z (a)\backslash \widetilde{\fix}_a \cong H\cap\Conj(a)$, 
where $\z (a)$ is the centralizer of $a$ in $G$, acting on 
$\widetilde{\fix}_a$ by 
$g\mapsto zg$ ($g\in \widetilde{\fix}_a,\; z\in \z (a)$). 
Hence $\sharp \widetilde{\fix}_a/\sharp \z (a)=\sharp(H\cap\Conj(a))$. 
If we take $H$ to be $G$, we have 
$\sharp G/\sharp \z (a)=\sharp\Conj(a)$. 
Combining these three equalities, we get the desired formula. 

\end{proof}

For an integer $N\geq 1$, take a subgroup $H\subseteq G=\SL_2(\Z/N\Z)$. 
By Lemma \ref{fix-conj}, we have 
\begin{equation}
\label{deltaH}
\delta_H
=1
  -3\frac{\sharp H\cap\Conj(\sigma)}{\sharp \Conj(\sigma)}
  -4\frac{\sharp H\cap\Conj(\tau)}{\sharp \Conj(\tau)}
  -6\frac{\sharp I\backslash G/H}{[G:H]}.
\end{equation}
Recall $\Conj(\alpha)
=\{g^{-1}\alpha g \in \SL_2(\Z/N\Z)|g\in \SL_2(\Z/N\Z)\}$ 
for $\alpha\in \SL_2(\Z/N\Z)$.

Assume $n\geq 1$ and $N=p^n$. For $0\leq s\leq n$,  
we have 
$$\fix_{u^{p^s}}=\{gH\in G/H|\sharp I(gH)\leq p^s\},$$ 
where 
$I(gH)$ means the $I$-orbit of $gH$ in $G/H$. We have 
$$G/H=\fix_{u^{p^n}}=\fix_u\amalg
\coprod_{s=0}^{n-1}(\fix_{u^{p^{s+1}}}\setminus \fix_{u^{p^s}})$$
and 
$$\sharp I\backslash G/H=\sharp \fix_u+\sum_{s=0}^{n-1}
\frac{1}{p^{s+1}}(\sharp \fix_{u^{p^{s+1}}}-\sharp \fix_{u^{p^s}})
=\sum_{s=0}^{n-1}\frac{p-1}{p^{s+1}}\sharp \fix_{u^{p^s}}
+\frac{1}{p^n}[G:H].$$
By Lemma \ref{fix-conj}, we get 
\begin{equation}
\label{cusp1}
\frac{\sharp I\backslash G/H}{[G:H]}=\sum_{s=0}^{n-1}
\frac{p-1}{p^{s+1}}\frac{\sharp H\cap\Conj(u^{p^s})}{\sharp \Conj(u^{p^s})}
+\frac{1}{p^n}. 
\end{equation}
In particular, we have 
\begin{equation}
\label{cusp2}
\frac{\sharp I\backslash G/H}{[G:H]}\leq\sum_{s=0}^{t-1}
\frac{p-1}{p^{s+1}}\frac{\sharp H\cap\Conj(u^{p^s})}{\sharp \Conj(u^{p^s})}
+\frac{1}{p^t}
\end{equation}
for $1\leq t\leq n$.

When $p=2$ or $3$, we study the conjugacy classes of the 
maximal subgroups of $\SL_2(\Z/p^2\Z)$.

\begin{lemma}
\label{sl2-gen}

We have 
$$\SL_2(\Z)=\langle \sigma,\tau\rangle 
=\langle \sigma,u\rangle 
=\langle \tau,u\rangle 
=\langle u,\,^t u\rangle. $$

\end{lemma}

\begin{proof}

The equality $\SL_2(\Z)=\langle \sigma,u\rangle$ is well-known 
([21, p.16]).
As $\sigma=u\tau$ and $(^t u)^{-1}u\cdot (^t u)^{-1}=\sigma$, 
we have the other 
equalities.  

\end{proof}

\begin{lemma}
\label{sl2-surj}

([21, Lemma 1.38])
For any positive integer $N$, the map 
$$\mod{N}:\SL_2(\Z)\longrightarrow\\SL_2(\Z/N\Z)$$ 
is surjective. 

\end{lemma}

\begin{lemma}
\label{sl2-u}

Let $N\geq 2$ be an integer and 
$A\subsetneq \SL_2(\Z/N^2\Z)$ be a proper subgroup. Assume $A$ maps 
surjectively $\bmod{\,N}$ onto $\SL_2(\Z/N\Z)$. Then 
$A$ has no element which is $\GL_2(\Z/N^2\Z)$-conjugate to $u$. 

\end{lemma}

\begin{proof}

Suppose $A$ has an element which is $\GL_2(\Z/N^2\Z)$-conjugate to $u$. 
Replacing $A$ by its $\GL_2(\Z/N^2\Z)$-conjugate, 
we may assume 
$A\ni u=\left(\begin{matrix} 1&1\\ 0&1 \end{matrix}\right)$. 
By the hypothesis, we see that 
$A$ contains an element $X$ satisfying 
$X\mod{N}=\left(\begin{matrix} 0&1\\ -1&0 \end{matrix}\right)$. 
Since $X$ has determinant $1$, we can write 
$X=\left(\begin{matrix} Na&1+Nb\\ -1+Nb&Nd \end{matrix}\right)$. 
We have 
$$A\ni 
\left(\begin{matrix} 1&Na\\ 0&1 \end{matrix}\right)
\left(\begin{matrix} Na&1+Nb\\ -1+Nb&Nd \end{matrix}\right)
\left(\begin{matrix} 1&Nd\\ 0&1 \end{matrix}\right)
=\left(\begin{matrix} 0&1+Nb\\ -1+Nb&0 \end{matrix}\right).$$ 
Hence 
$$A\ni 
\left(\begin{matrix} 0&1+Nb\\ -1+Nb&0 \end{matrix}\right)^{-1}
\left(\begin{matrix} 1&1\\ 0&1 \end{matrix}\right)
\left(\begin{matrix} 0&1+Nb\\ -1+Nb&0 \end{matrix}\right)
=\left(\begin{matrix} 1&0\\ -1+2Nb&1 \end{matrix}\right).$$ 
Take an integer $\widetilde{b}$ such that $\widetilde{b}\mod{N^2}=b$, 
and we get 
$A\ni
\left(\begin{matrix} 1&0\\ -1+2Nb&1 \end{matrix}\right)
^{-1-2N\widetilde{b}}
=\left(\begin{matrix} 1&0\\ 1&1 \end{matrix}\right)$. 
Then $A\supseteq \langle u,\,^t u\rangle =\SL_2(\Z/N^2\Z)$, 
a contradiction. 

\end{proof}

\begin{remark}

If $N=p\geq 5$ is a prime, 
there is no subgroup $A$ as in the previous lemma 
by Lemma 2.1.

\end{remark}

\begin{lemma}
\label{table}

The conjugacy classes of $\SL_2(\Z/4\Z)$ are the following:

\begin{align*}
\Conj(1)=&\{1\},\\
\Conj(-1)=&\{-1\},\\
\Conj(\sigma)=&\left\{\sigma,
\left(\begin{matrix} 1&2 \\ -1&-1 \end{matrix}\right),
\left(\begin{matrix} 2&1 \\ -1&2 \end{matrix}\right),
\left(\begin{matrix} -1&2 \\ -1&1 \end{matrix}\right),
\left(\begin{matrix} 1&1 \\ 2&-1 \end{matrix}\right),
\left(\begin{matrix} -1&1 \\ 2&1 \end{matrix}\right)
\right\},\\
\Conj(-\sigma)=&-\Conj(\sigma),\\
\Conj(\tau)=&\left\{\tau,
\left(\begin{matrix} 1&-1 \\ 1&0 \end{matrix}\right),
\left(\begin{matrix} 0&1 \\ -1&1 \end{matrix}\right),
\left(\begin{matrix} 0&-1 \\ 1&1 \end{matrix}\right),
\left(\begin{matrix} -1&1 \\ 1&2 \end{matrix}\right),
\left(\begin{matrix} -1&-1 \\ -1&2 \end{matrix}\right),\right.\\
&\left.\left(\begin{matrix} 2&1 \\ 1&-1 \end{matrix}\right),
\left(\begin{matrix} 2&-1 \\ -1&-1 \end{matrix}\right)
\right\},\\
\Conj(-\tau)=&-\Conj(\tau),\\
\Conj(u)=&\left\{u,
\left(\begin{matrix} 1&0 \\ -1&1 \end{matrix}\right),
\left(\begin{matrix} 2&1 \\ -1&0 \end{matrix}\right),
\left(\begin{matrix} -1&0 \\ -1&-1 \end{matrix}\right),
\left(\begin{matrix} 0&1 \\ -1&2 \end{matrix}\right),
\left(\begin{matrix} -1&1 \\ 0&-1 \end{matrix}\right)
\right\},\\
\Conj(-u)=&-\Conj(u),\\
\Conj(u^2)=&\left\{u^2,
\left(\begin{matrix} 1&0 \\ 2&1 \end{matrix}\right),
\left(\begin{matrix} -1&2 \\ 2&-1 \end{matrix}\right)
\right\},\\
\Conj(-u^2)=&-\Conj(u^2).
\end{align*}

\end{lemma}

Now we determine maximal subgroups of $\SL_2(\Z/4\Z)$ whose image mod $2$
is $\SL_2(\Z/2\Z)$. 

\begin{lemma}
\label{A1-A2}

Let $A\subsetneq \SL_2(\Z/4\Z)$ be a proper subgroup. Assume $A$ maps 
surjectively $\bmod{\,2}$ onto $\SL_2(\Z/2\Z)$. Then $A$ is conjugate to 
$$A_1:=\langle \sigma,
\left(\begin{matrix} 1&1\\ 2&-1 \end{matrix}\right)\rangle,$$
which is a maximal subgroup, and is not a 
normal subgroup. We have 
$\langle A_1, u^2\rangle
=\SL_2(\Z/4\Z)$ and 
$A_1\backslash \SL_2(\Z/4\Z)=
\{1,u,u^2,u^{-1}\}$. 

\end{lemma}

\begin{proof}

By the hypothesis, we see that $A$ contains a lift of 
$u\in \SL_2(\Z/2\Z)$. 
We have 
$$\{X\in \SL_2(\Z/4\Z)|X\bmod{2}=u\}
=\{
u,u^{-1},-u,-u^{-1},
\alpha,\beta,\alpha^{-1},\beta^{-1}
\},$$
where 
$\alpha=\left(\begin{matrix} 1&1\\ 2&-1 \end{matrix}\right)$
and
$\beta=\left(\begin{matrix} -1&1\\ 2&1 \end{matrix}\right)$. 
The first four elements are 
$\GL_2(\Z/4\Z)$-conjugate to $u$. 
The next two elements are $\SL_2(\Z/4\Z)$-conjugate to $\sigma$ and 
the last two to $\sigma^{-1}$. 
Since $A$ has no element which is $\GL_2(\Z/4\Z)$-conjugate to $u$, we have 
$A\cap \Conj(\sigma)\ne\emptyset$. 
Replacing $A$ by its conjugate, we may assume $A\ni \sigma$. 
Repeating the same argument, we see that 
$A$ contains at least one of 
$\alpha, \beta, \alpha^{-1}, \beta^{-1}$. 
We have 
$A_1=\langle\sigma, \alpha\rangle
=\langle\sigma, \alpha^{-1}\rangle$. 
Put 
$A_2:=\langle\sigma, \beta\rangle
=\langle\sigma, \beta^{-1}\rangle$. 
By calculation, we have 
\begin{align*}
A_1=&\left\{
\left(\begin{matrix} 1&0\\ 0&1 \end{matrix}\right),
\left(\begin{matrix} -1&0\\ 0&-1 \end{matrix}\right),
\left(\begin{matrix} 0&1\\ -1&0 \end{matrix}\right),
\left(\begin{matrix} 0&-1\\ 1&0 \end{matrix}\right),
\left(\begin{matrix} 2&1\\ 1&1 \end{matrix}\right),
\left(\begin{matrix} 1&-1\\ -1&2 \end{matrix}\right),\right.\\
&\left.\left(\begin{matrix} 2&-1\\ -1&-1 \end{matrix}\right),
\left(\begin{matrix} -1&1\\ 1&2 \end{matrix}\right),
\left(\begin{matrix} -1&2\\ -1&1 \end{matrix}\right),
\left(\begin{matrix} 1&2\\ 1&-1 \end{matrix}\right),
\left(\begin{matrix} 1&1\\ 2&-1 \end{matrix}\right),
\left(\begin{matrix} -1&-1\\ 2&1 \end{matrix}\right)\right\}
\end{align*}
and 
\begin{align*}
A_2=&\left\{
\left(\begin{matrix} 1&0\\ 0&1 \end{matrix}\right),
\left(\begin{matrix} -1&0\\ 0&-1 \end{matrix}\right),
\left(\begin{matrix} 0&1\\ -1&0 \end{matrix}\right),
\left(\begin{matrix} 0&-1\\ 1&0 \end{matrix}\right),
\left(\begin{matrix} 2&-1\\ -1&1 \end{matrix}\right),
\left(\begin{matrix} 1&1\\ 1&2 \end{matrix}\right),\right.\\
&\left.\left(\begin{matrix} 2&1\\ 1&-1 \end{matrix}\right),
\left(\begin{matrix} -1&-1\\ -1&2 \end{matrix}\right),
\left(\begin{matrix} -1&2\\ 1&1 \end{matrix}\right),
\left(\begin{matrix} 1&2\\ -1&-1 \end{matrix}\right),
\left(\begin{matrix} 1&-1\\ 2&-1 \end{matrix}\right),
\left(\begin{matrix} -1&1\\ 2&1 \end{matrix}\right)\right\}.
\end{align*}
The equalities 
$u^{-1}\sigma u
=\left(\begin{matrix} 1&2\\ -1&-1 \end{matrix}\right)
\in A_2$ 
and 
$u^{-1}\alpha u
=\left(\begin{matrix} -1&1\\ 2&1 \end{matrix}\right)
\in A_2$
show $A_2=u^{-1}A_1 u$. 
We also have 
$\langle A_1,u^2\rangle
=\SL_2(\Z/4\Z)$ by calculation. 

\end{proof}

\begin{remark}
\label{A_1-conj}

In the above lemma, all the subgroups of $\SL_2(\Z/4\Z)$ 
conjugate to $A_1$ are 
$A_1$, $u^{-1}A_1 u$, $u^{-2}A_1 u^2$ and $uA_1 u^{-1}$. 

\end{remark}

Now we calculate the number of elements conjugate to  $\sigma,\tau,u$ 
in the maximal subgroups $B,C,D,E$ introduced in Section 3.

\begin{lemma}
\label{bcde-conj}

In $\SL_2(\Z/p\Z)$, the number of the 
elements conjugate to $\sigma$, $\tau$, $u$ in 
$B$, $C$, $D$, $E$ are as follows. 

1.
\begin{align*}
  B\cap\Conj(\sigma)
  &= \begin{cases}
    \emptyset & \text{if $p\equiv -1 \mod{4}$},\\
    \left\{\left(\begin{matrix} \sqrt{-1}&*\\ 0&-\sqrt{-1} \end{matrix}\right),
    \left(\begin{matrix} -\sqrt{-1}&*\\ 0&\sqrt{-1} \end{matrix}\right)\right\}
    & \text{if $p\equiv 1 \mod{4}$},\\
    \left\{\left(\begin{matrix} 1&1\\ 0&1 \end{matrix}\right)\right\}
    & \text{if $p=2$}.
  \end{cases}
\\
    B\cap\Conj(\tau)
  &= \begin{cases}
    \emptyset & \text{if $p\equiv -1 \mod{3}$},\\
    \left\{\left(\begin{matrix} -\zeta_3&*\\ 0&-\zeta_3^{-1} 
    \end{matrix}\right),
    \left(\begin{matrix} -\zeta_3^{-1}&*\\ 0&-\zeta_3 
    \end{matrix}\right)\right\}
    & \text{if $p\equiv 1 \mod{3}$},\\
    \left\{\left(\begin{matrix} -1&1\\ 0&-1 \end{matrix}\right)
    \right\}
    & \text{if $p=3$}.
  \end{cases}
\\
    B\cap\Conj(u)
  &=\left\{\left(\begin{matrix} 1&s^2\\ 0&1 \end{matrix}\right)|
  (p,s)=1\right\}.
\end{align*}
Hence we have 
\begin{align*}
  \sharp B\cap\Conj(\sigma)
  &= \begin{cases}
    0 & \text{if $p\equiv -1 \mod{4}$},\\
    2p
    & \text{if $p\equiv 1 \mod{4}$},\\
    1
    & \text{if $p=2$}.
  \end{cases}
\\
   \sharp B\cap\Conj(\tau)
  &= \begin{cases}
    0 & \text{if $p\equiv -1 \mod{3}$},\\
    2p
    & \text{if $p\equiv 1 \mod{3}$},\\
    1
    & \text{if $p=3$}.
  \end{cases}
\\
   \sharp B\cap\Conj(u)
  &=\begin{cases}
     \frac{1}{2}(p-1) & \text{if $p\geq 3$},\\
     1 & \text{if $p=2$}.
    \end{cases}
\end{align*}

2.
\begin{align*}
  C\cap
  \Conj(\sigma)
 &= \begin{cases}
    \left\{\left(\begin{matrix} 0&x\\ -x^{-1}&0 \end{matrix}\right)|
    (p,x)=1\right\}\\
    \text{if $p\equiv -1 \mod{4}$},\\
    \left\{\left(\begin{matrix} \sqrt{-1}&0\\ 0&-\sqrt{-1} \end{matrix}\right),
    \left(\begin{matrix} -\sqrt{-1}&0\\ 0&\sqrt{-1} \end{matrix}\right),
    \left(\begin{matrix} 0&x\\ -x^{-1}&0 \end{matrix}\right)|
    (p,x)=1
    \right\}\\
    \text{if $p\equiv 1 \mod{4}$},\\
    \left\{\left(\begin{matrix} 0&1\\ 1&0 \end{matrix}\right)\right\}\\
    \text{if $p=2$}.
  \end{cases}
\\
  C\cap
  \Conj(\tau)
  &= \begin{cases}
    \emptyset & \text{if $p\not\equiv 1 \mod{3}$},\\
    \left\{\left(\begin{matrix} -\zeta_3&0\\ 0&-\zeta_3^{-1} 
    \end{matrix}\right),
    \left(\begin{matrix} -\zeta_3^{-1}&0\\ 0&-\zeta_3 
    \end{matrix}\right)\right\}
    & \text{if $p\equiv 1 \mod{3}$}.
  \end{cases}
\\
  C\cap
  \Conj(u)
  &= \begin{cases}
    \emptyset & \text{if $p\geq 3$},\\
    \left\{\left(\begin{matrix} 0&1\\ 1&0 
    \end{matrix}\right)\right\}& \text{if $p=2$}.
  \end{cases}
\end{align*}
Hence we have 
\begin{align*}
   \sharp C\cap\Conj(\sigma)
  &= \begin{cases}
    p-1
    & \text{if $p\equiv -1 \mod{4}$},\\
    p+1
    & \text{if $p\equiv 1 \mod{4}$},\\
    1
    & \text{if $p=2$}.
  \end{cases}
\\
  \sharp C\cap
  \Conj(\tau)
  &= \begin{cases}
    0 & \text{if $p\not\equiv 1 \mod{3}$},\\
    2
    & \text{if $p\equiv 1 \mod{3}$}.
  \end{cases}
\\
  \sharp C\cap
  \Conj(u)
  &= \begin{cases}
    0 & \text{if $p\geq 3$},\\
    1 & \text{if $p=2$}.
  \end{cases}
\end{align*}

3.
\begin{align*}
  &D\cap
  \Conj(\tau)\\
  &= \begin{cases}
    \left\{
    \left(\begin{matrix} 1/2&1/2\cdot\sqrt{-3/\lambda}
    \\ 1/2\cdot\sqrt{-3\lambda}&1/2 \end{matrix}\right),
    \left(\begin{matrix} 1/2&-1/2\cdot\sqrt{-3/\lambda}
    \\ -1/2\cdot\sqrt{-3\lambda}&1/2 \end{matrix}\right)
    \right\}\\
    \text{if $p\geq 5$ and $p\equiv -1 \mod{3}$},\\
    \emptyset \\
    \text{if $p=3$ or $p\equiv 1 \mod{3}$}.
  \end{cases}
\end{align*}
\begin{align*}
  \sharp D\cap
  \Conj(\sigma)
  &= \begin{cases}
    p+3 & \text{if $p\equiv -1 \mod{4}$},\\
    p+1 & \text{if $p\equiv 1 \mod{4}$},
  \end{cases}
\\
 \sharp D\cap
  \Conj(\tau) 
    &= \begin{cases}
    2
    &\text{if $p\geq 5$ and $p\equiv -1 \mod{3}$},\\
    0
    &\text{if $p=3$ or $p\equiv 1 \mod{3}$}.
  \end{cases}
\\
 \sharp D\cap
  \Conj(u)
  &= 0 \text{ \;\;if $p\geq 3$}.
\end{align*}

4.
\begin{align*}
\sharp E\cap\Conj(\sigma)&\leq
  \begin{cases}
  30 & \text{if $p\equiv \pm 1\mod{5}$},\\
  18 & \text{if $p\geq 5$ and $p\not\equiv \pm 1\mod{5}$}.
  \end{cases}
\\
\sharp E\cap\Conj(\tau)&\leq
  \begin{cases}
  20 & \text{if $p\equiv \pm 1\mod{5}$},\\
  8 & \text{if $p\geq 5$ and $p\not\equiv \pm 1\mod{5}$}.
  \end{cases}
\\
\sharp E\cap\Conj(u)&=0 \text{ \;\;if $p\geq 5$}.
\end{align*}

\end{lemma}

\begin{proof}

1 is easy. 

2. Note that if $p\geq 3$ then 
for any $x\in (\Z/p\Z)^{\times}$ there exist $a,c\in \Z/p\Z$ 
such that $a^2+c^2=x^{-1}$. Thus we have 
$\left(\begin{matrix} a&-xc\\ c&xa \end{matrix}\right)\in \SL_2(\Z/p\Z)$ 
and 
$\left(\begin{matrix} a&-xc\\ c&xa \end{matrix}\right)^{-1}
\left(\begin{matrix} 0&1\\ -1&0 \end{matrix}\right)
\left(\begin{matrix} a&-xc\\ c&xa \end{matrix}\right)=
\left(\begin{matrix} 0&x\\ -x^{-1}&0 \end{matrix}\right)$. 

3. Assume $p\geq 3$. 
We show 
$\sharp \left\{
\left(\begin{matrix} x&y\\ -\lambda y&-x \end{matrix}\right)\right\}
\cap \Conj(\sigma)=p+1$. 
Put 
$D_1=\left\{\left(\begin{matrix} x&y\\ \lambda y&x \end{matrix}\right)
\right\}\cap \SL_2(\Z/p\Z)$ and 
$D_2=\left\{\left(\begin{matrix} x&y\\ -\lambda y&-x \end{matrix}\right)
\right\}\cap \SL_2(\Z/p\Z)$. 
Since $p\geq 3$, we have $\Conj(\sigma)=\SL_2(\Z/p\Z)^{\tr=0}$. 
Thus 
$\left\{
\left(\begin{matrix} x&y\\ -\lambda y&-x \end{matrix}\right)\right\}
\cap \Conj(\sigma)
=D_2=D\setminus D_1$. 
We have 
$D=\Ker (\det:
\left\{\left(\begin{matrix} x&y\\ \lambda y&x \end{matrix}\right),
\left(\begin{matrix} x&y\\ -\lambda y&-x \end{matrix}\right)
\right\}
\longrightarrow (\Z/p\Z)^{\times})$ 
and 
$D_1=\Ker (\det:
\left\{\left(\begin{matrix} x&y\\ \lambda y&x \end{matrix}\right)
\right\}
\longrightarrow (\Z/p\Z)^{\times})$. 
Since $\det:
\left\{\left(\begin{matrix} x&y\\ \lambda y&x \end{matrix}\right)
\right\}
\longrightarrow (\Z/p\Z)^{\times}$ 
is surjective, 
we have 
$\sharp D_1=(p^2-1)/(p-1)=p+1$ 
and 
$\sharp D_1=2(p^2-1)/(p-1)=2(p+1)$. 
Therefore 
$\sharp \left\{
\left(\begin{matrix} x&y\\ -\lambda y&-x \end{matrix}\right)\right\}
\cap \Conj(\sigma)=2(p+1)-(p+1)=p+1$.

4. In $A_5$, there are $15$ elements of order $2$ and 
$20$ elements of order $3$. 
In $S_4$, there are $9$ elements of order $2$ and 
$8$ elements of order $3$. 

\end{proof}

Next we calculate the number of elements conjugate to  $\sigma,\tau,u$ 
in $A_1\subseteq \SL_2(\Z/4\Z)$.

\begin{lemma}
\label{A1-conj}

In $\SL_2(\Z/4\Z)$, we have 
\begin{align*}
A_1\cap \Conj(\sigma)&=\left\{
\sigma,
\left(\begin{matrix} -1&2\\ -1&1 \end{matrix}\right),
\left(\begin{matrix} 1&1\\ 2&-1 \end{matrix}\right)
\right\},\\
A_1\cap \Conj(\tau)&=\left\{
\left(\begin{matrix} 2&-1\\ -1&-1 \end{matrix}\right),
\left(\begin{matrix} -1&1\\ 1&2 \end{matrix}\right)
\right\},\\
A_1\cap \Conj(u)&=\emptyset ,\\
A_1\cap \Conj(u^2)&=\emptyset ,
\end{align*}
where $A_1=\langle \sigma,
\left(\begin{matrix} 1&1\\ 2&-1 \end{matrix}\right)\rangle$. 
Hence we have 
\begin{align*}
\sharp A_1\cap \Conj(\sigma)&=3,\\
\sharp A_1\cap \Conj(\tau)&=2,\\
\sharp A_1\cap \Conj(u)&=0,\\
\sharp A_1\cap \Conj(u^2)&=0.
\end{align*}

\end{lemma}

\begin{proof}

It follows from Lemma \ref{table}. 

\end{proof}


\section{Calculation of conjugate elements in $\SL_2(\Z/p^n\Z)$}

We calculate the number of elements conjugate to  $\sigma,\tau,u$ in 
$\SL_2(\Z/p^n\Z)$.

\begin{lemma}
\label{sl2-conj}

Let $n\geq 1$ be an integer. 
In $\SL_2(\Z/p^n\Z)$  
we have 

\begin{align*}
\sharp\Conj(\sigma)
&= \begin{cases}
    (p-1)p^{2n-1}& \text{if $p\equiv -1 \mod{4}$},\\
    (p+1)p^{2n-1}& \text{if $p\equiv 1 \mod{4}$},\\
    3            & \text{if $p=2$ and $n=1$},\\
    3\cdot 2^{2n-3}& \text{if $p=2$ and $n\geq 2$}, 
  \end{cases}
\\
\sharp\Conj(\tau)
&= \begin{cases}
    (p-1)p^{2n-1}& \text{if $p\equiv -1 \mod{3}$},\\
    (p+1)p^{2n-1}& \text{if $p\equiv 1 \mod{3}$},\\
    4\cdot 3^{2n-2}& \text{if $p=3$},
  \end{cases}
\\
\sharp\Conj(u)
&= \begin{cases}
    \frac{1}{2}(p^2-1)p^{2n-2} & \text{if $p\geq 3$},\\
    3               & \text{if $p=2$ and $n=1$},\\
    6               & \text{if $p=2$ and $n=2$},\\
    3\cdot 2^{2n-4} & \text{if $p=2$ and $n\geq 3$}. 
  \end{cases}
\end{align*}

\end{lemma}

\begin{proof}

As we have seen in the proof of Lemma \ref{fix-conj}, we have an equality 
$\sharp\Conj(\alpha)=\sharp \SL_2(\Z/p^n\Z)/\sharp \z (\alpha)$ for any 
$\alpha\in \SL_2(\Z/p^n\Z)$. As is well known, 
$\sharp \SL_2(\Z/p^n\Z)=(p+1)(p-1) p^{3n-2}$. 
Now we calculate $\sharp \z(\alpha)$ for $\alpha=\sigma,\ \tau,\ u$. 
For $\alpha=\sigma$, we have 
\begin{equation*}
\z(\sigma)=(\Z/p^n\Z[\sigma])^{\det=1}\cong 
  \begin{cases}
    \W_n(\F_{p^2})^{\Norm=1} & \text{if $p\equiv -1 \mod{4}$},\\
    (\Z/p^n\Z\times \Z/p^n\Z)^{\Norm=1} & \text{if $p\equiv 1 \mod{4}$},\\
    (\Z_2[\sqrt{-1}]/2^n\Z_2[\sqrt{-1}])^{\Norm=1} & \text{if $p=2$},
  \end{cases}
\end{equation*}
hence 
\begin{equation*}
\sharp\z(\sigma)=
  \begin{cases}
     (p+1)p^{n-1}& \text{if $p\equiv -1 \mod{4}$},\\
     (p-1) p^{n-1}& \text{if $p\equiv 1 \mod{4}$},\\
     2 & \text{if $p=2$ and $n=1$},\\
     2^{n+1} & \text{if $p=2$ and $n\geq 2$}. 
  \end{cases}
\end{equation*}

For $\alpha=\tau$, we have 
\begin{equation*}
\z(\tau)=(\Z/p^n\Z[\tau])^{\det=1}\cong 
  \begin{cases}
    \W_n(\F_{p^2})^{\Norm=1} & \text{if $p\equiv -1 \mod{3}$},\\
    (\Z/p^n\Z\times \Z/p^n\Z)^{\Norm=1} & \text{if $p\equiv 1 \mod{3}$},\\
    (\Z_3[\zeta_3]/3^n\Z_3[\zeta_3])^{\Norm=1} & \text{if $p=3$},
  \end{cases}
\end{equation*}
where $\zeta_3$ is a primitive third root of unity. Hence 
\begin{equation*}
\sharp\z(\tau)=
  \begin{cases}
     (p+1)p^{n-1}& \text{if $p\equiv -1 \mod{3}$},\\
     (p-1) p^{n-1}& \text{if $p\equiv 1 \mod{3}$},\\
     2\cdot 3^n& \text{if $p=3$}. 
  \end{cases}
\end{equation*}

For $\alpha=u$, we have 
$\z(u)=\left\{\left(\begin{matrix} a&b\\ 0&a \end{matrix}\right)
|a^2=1\right\}$, hence 
\begin{equation*}
\sharp\z(u)
= \begin{cases}
    2p^n            & \text{if $p\geq 3$},\\
    2               & \text{if $p=2$ and $n=1$},\\
    8               & \text{if $p=2$ and $n=2$},\\
    2^{n+2}         & \text{if $p=2$ and $n\geq 3$}, 
  \end{cases}
\end{equation*}
as required. 

\end{proof}

Next we calculate the number of elements conjugate to  $u^{p^r}$ 
in $\SL_2(\Z/p^{r+n}\Z)$.

\begin{lemma}
\label{sl2-upr}

Assume $r\geq 0$ and $n\geq 1$. In $\SL_2(\Z/p^{r+n}\Z)$, we have 
\begin{equation*}
\sharp\Conj(u^{p^r})
= \begin{cases}
    \frac{1}{2}(p^2-1)p^{2n-2} & \text{if $p\geq 3$},\\
    3                   & \text{if $p=2$ and $n=1$},\\
    6                   & \text{if $p=2$ and $n=2$},\\
    3\cdot 2^{2n-4} & \text{if $p=2$ and $n\geq 3$}.
\end{cases}
\end{equation*}
\end{lemma}

\begin{proof}

We have 
$\z(u^{p^r})=\left\{\left(\begin{matrix} a&b\\ c&d \end{matrix}\right)
 |p^r c=0,\ p^r(a-d)=0,\ ad-bc=1\right\}$, 
which maps surjectively 
$\bmod{\ p^n}$ onto 
$\left\{\left(\begin{matrix} a&b\\ 0&a \end{matrix}\right)
\in \SL_2(\Z/p^n\Z)\right\}$, 
and each fiber of the $\bmod{\ p^n}$ map 
consists of $p^{3r}$ elements. 

\end{proof}

We study the fiber of the mod $p^m$ map 
$$\SL_2(\Z/p^n\Z)\cap\Conj(\alpha)\longrightarrow
\SL_2(\Z/p^m\Z)\cap\Conj(\alpha),$$
where $\alpha=\sigma,\tau,u^{p^r}$.

Take two integers $m,n$ with $1\leq m\leq n$. 
As in Section 1, let 
$$f_{n,m} : \SL_2(\Z/p^n\Z)\longrightarrow \SL_2(\Z/p^m\Z)$$ 
be the mod $p^m$ map. 
For $\alpha=\sigma,\tau,u^{p^r}$ ($r\geq 0$), 
put 
$$V_{\alpha}^{r+n,r+m}=
\alpha^{-1}(f_{r+n,r+m}^{-1}(\alpha)\cap \Conj(\alpha))
\subseteq \SL_2(\Z/p^{r+n}\Z).$$ 
When $\alpha=\sigma,\tau$, we always take $r=0$. 
We sometimes omit the superscripts 
$r,n,m$ and 
simply write $V_{\alpha}$.

\begin{lemma}
\label{V-str1}

Let $r\geq 0$ and $1\leq m<n$ be integers. 
Assume $n\leq 2m$. 
If $p\geq 3$, then each $V_{\alpha}^{r+n,r+m}$ 
($\alpha=\sigma,\tau,\, u^{p^r}$) 
is a subgroup of 
$(1+p^{r+m}\M_2(\Z/p^{r+n}\Z))^{\det=1}\cong\M_2(\Z/p^{n-m}\Z)^{\tr=0}$ 
and isomorphic to $(\Z/p^{n-m}\Z)^2$. 
If $p=2$, then $V_{\sigma}^{n,m}$ 
for $m\geq 2$ ($\resp V_{\tau}^{n,m}$ for any $m$, 
$\resp V_{u^{p^r}}^{r+n,r+m}$ for $m\geq 3$) is a subgroup of 
$(1+p^{r+m}\M_2(\Z/p^{r+n}\Z))^{\det=1}\cong\M_2(\Z/p^{n-m}\Z)^{\tr=0}$ 
and isomorphic to $(\Z/p^{n-m}\Z)^2$. 
Explicitly: 

If $p\geq 3$ or ($p=2$ and $m\geq 2$), we have 
$$V_{\sigma}^{n,m}=\left\{1+p^m
  \left(\begin{matrix} a&b\\ b&-a \end{matrix}\right)
  \in \SL_2(\Z/p^n\Z)\right\}\cong \left\{
  \left(\begin{matrix} a&b\\ b&-a \end{matrix}\right)
  \in \M_2(\Z/p^{n-m}\Z)\right\};$$
For any $p$ and $m$, we have 
$$V_{\tau}^{n,m}=\left\{1+p^m
  \left(\begin{matrix} a&b\\ b-a&-a \end{matrix}\right)
  \in \SL_2(\Z/p^n\Z)\right\}\cong \left\{
  \left(\begin{matrix} a&b\\ b-a&-a \end{matrix}\right)
  \in \M_2(\Z/p^{n-m}\Z)\right\};$$
If $p\geq 3$ or ($p=2$ and $m\geq 3$), we have 
$$V_{u^{p^r}}^{r+n,r+m}=\left\{1+p^{r+m}
  \left(\begin{matrix} a&b\\ 0&-a \end{matrix}\right)
  \in \GL_2(\Z/p^{r+n}\Z)\right\}\cong \left\{
  \left(\begin{matrix} a&b\\ 0&-a \end{matrix}\right)
  \in \M_2(\Z/p^{n-m}\Z)\right\}.$$

In particular, the inverse image of one element by the following maps 
consists of $p^2$ elements: 
\begin{itemize}
\item
$\mod{p^m}:\SL_2(\Z/p^{m+1}\Z)\cap\Conj(\sigma)
\longrightarrow \SL_2(\Z/p^m\Z)\cap\Conj(\sigma)$ 
\ if $p\geq 3$ or ($p=2$ and $m\geq 2$), 
\item
$\mod{p^m}:\SL_2(\Z/p^{m+1}\Z)\cap\Conj(\tau)
\longrightarrow \SL_2(\Z/p^m\Z)\cap\Conj(\tau)$ 
\ for any $p$ and any $m\geq 1$, 
\item
$\mod{p^{r+m}}:\SL_2(\Z/p^{r+m+1}\Z)\cap\Conj(u^{p^r})
\longrightarrow \SL_2(\Z/p^{r+m}\Z)\cap\Conj(u^{p^r})$ 
\ if $p\geq 3$ or ($p=2$ and $m\geq 3$). 
\end{itemize}

\end{lemma}

\begin{proof}

Assume $\alpha =\sigma$. We will show
 $$f_{n,m}^{-1}(\sigma)\cap \Conj(\sigma)=
 \left\{\left(\begin{matrix} p^m a&1+p^m b\\ -1+p^m b&-p^m a 
 \end{matrix}\right)
 \in \SL_2(\Z/p^n\Z)\right\}.$$ 
Looking at the determinant and the trace, we get one inclusion 
``$\subseteq$''. In particular, we have 
$\sharp f_{n,m}^{-1}(\sigma ')\cap \Conj(\sigma)\leq p^{2(n-m)}$ for 
any $\sigma '\in \SL_2(\Z/p^m\Z)\cap\Conj(\sigma)$. 
Since $\bmod{\ p^m}:\SL_2(\Z/p^n\Z)\cap\Conj(\sigma)
\longrightarrow \SL_2(\Z/p^m\Z)\cap\Conj(\sigma)$ is surjective, 
Lemma \ref{sl2-conj} shows that the other inclusion ``$\supseteq$'' holds 
if $p\geq 3$ or ($p=2$ and $m\geq 2$).

We can show 
 $$f_{n,m}^{-1}(\tau)\cap \Conj(\tau)=
 \left\{\left(\begin{matrix} 1+p^m a&1+p^m b\\ -1+p^m (b-a)&-p^m a 
 \end{matrix}\right)
 \in \SL_2(\Z/p^n\Z)\right\}$$ 
for any $p$, $m$ and 
 $$f_{r+n,r+m}^{-1}(u^{p^r})\cap \Conj(u^{p^r})=
 \left\{\left(\begin{matrix} 1+p^{r+m} a&p^r+p^{r+m} b\\ 0 &1-p^{r+m} a 
\end{matrix}\right)
 \in \SL_2(\Z/p^{r+n}\Z)\right\}$$ 
for $p\geq 3$ or ($p=2$ and $n\geq 3$) 
similarly by using Lemma \ref{sl2-conj} and \ref{sl2-upr}.

\end{proof}

\begin{remark}
\label{inv=2}

If $p=2$, the inverse image of one element by 
the following maps consists of $2$ elements:
\begin{itemize}
\item
$\mod{2}:\SL_2(\Z/2^2\Z)\cap \Conj(\sigma)\longrightarrow
\SL_2(\Z/2\Z)\cap \Conj(\sigma)$,
\item
$\mod{2^{r+2}}:\SL_2(\Z/2^{r+3}\Z)\cap \Conj(u^{2^r})\longrightarrow
\SL_2(\Z/2^{r+2}\Z)\cap \Conj(u^{2^r})$,
\item
$\mod{2^{r+1}}:\SL_2(\Z/2^{r+2}\Z)\cap \Conj(u^{2^r})\longrightarrow
\SL_2(\Z/2^{r+1}\Z)\cap \Conj(u^{2^r})$.
\end{itemize}

\end{remark}

From now on, we always assume the hypothesis in Lemma \ref{V-str1} when we 
write $V_{\alpha}^{r+n,r+m}$ (for $\alpha=\sigma,\tau,u^{p^r}$ 
and their conjugates: 
defined below), 
so that $V_{\alpha}^{r+n,r+m}$ is a free 
$\Z/p^{n-m}\Z$-submodule of rank $2$ of 
$(1+p^{r+m}\M_2(\Z/p^{r+n}\Z))^{\det=1}\cong \M_2(\Z/p^{n-m}\Z)^{\tr=0}$.

\begin{lemma}
\label{V-str2}

Let $r\geq 0$ and $1\leq m<n$ be integers. 
Assume $n\leq 2m$. 
For $\alpha=\sigma,\, \tau,\, u^{p^r}$, we have 
$$V_{\alpha}^{r+n,r+m}=\{1+p^m(X\alpha^{-1}-\alpha^{-1}X)|
X\in \M_2(\Z/p^{r+n}\Z)\}.$$

\end{lemma}

\begin{proof}

If $X=\left(\begin{matrix} x&y\\ z&w \end{matrix}\right)$, 
we have 
$$X\sigma^{-1}-\sigma^{-1}X
=\left(\begin{matrix} y+z&-x+w\\ w-x&-z-y \end{matrix}\right),$$
$$X\tau^{-1}-\tau^{-1}X
=\left(\begin{matrix} y+z&-x+y+w\\ w-x-z&-z-y \end{matrix}\right)$$
and
$$Xu^{-p^r}-u^{-p^r}X
=\left(\begin{matrix} p^rz&-p^r (x-w)\\ 0&-p^r z \end{matrix}\right).$$

\end{proof}

\begin{lemma}
\label{V1-V2}

Let $r\geq 0$ and $1\leq m<n$ be integers. 
For $\alpha=\sigma,\tau,u^{p^r}$, 
take an element $\alpha '\in \Conj(\alpha)\subseteq\SL_2(\Z/p^{r+m}\Z)$. 
Suppose two elements $\alpha_1,\alpha_2\in 
\Conj(\alpha)\subseteq\SL_2(\Z/p^{r+n}\Z)$ satisfy 
$\alpha_1\equiv\alpha_2\equiv\alpha '\mod{p^{r+m}}$. 
If $n\leq 2m$, 
then we have 
$${\alpha_1}^{-1}(f_{r+n,r+m}^{-1}(\alpha ')\cap \Conj(\alpha))
={\alpha_2}^{-1}(f_{r+n,r+m}^{-1}(\alpha ')\cap \Conj(\alpha)).$$

\end{lemma}

\begin{proof}

It suffices to show the equality in the case $\alpha_2=\alpha$. 
By the hypothesis, we have 
$\alpha_1\in (f_{r+n,r+m}^{-1}(\alpha)\cap \Conj(\alpha))
=\alpha V_{\alpha}$. 
We see that 
${\alpha_1}^{-1}(f_{r+n,r+m}^{-1}(\alpha)\cap \Conj(\alpha))
={\alpha_1}^{-1}\alpha V_{\alpha}
=\{{\alpha_1}^{-1}\alpha+(p^m{\alpha_1}^{-1}\alpha)
(X\alpha^{-1}-\alpha^{-1}X)\}
=\{{\alpha_1}^{-1}\alpha+p^m(X\alpha^{-1}-\alpha^{-1}X)\}$. 
Since 
$({\alpha_1}^{-1}\alpha)^{-1}=\alpha^{-1}\alpha_1\in V_{\alpha}$, 
we have 
${\alpha_1}^{-1}\alpha\in V_{\alpha}
=\{1+p^m(X\alpha^{-1}-\alpha^{-1}X)\}$. 
Therefore 
$\alpha_1^{-1}(f_{r+n,r+m}^{-1}(\alpha)\cap \Conj(\alpha))
=\{1+p^m(X\alpha^{-1}-\alpha^{-1}X)\}
=\alpha^{-1}(f_{r+n,r+m}^{-1}(\alpha)\cap \Conj(\alpha))$, 
as required. 

\end{proof}

In the above lemma, we define 
$$V_{\alpha '}^{r+n,r+m}:=
{\alpha_1}^{-1}(f_{r+n,r+m}^{-1}(\alpha ')\cap \Conj(\alpha)).$$
Note that we have 
$V_{\alpha '}^{r+n,r+m}=\{1+p^m(X{\alpha'}^{-1}-{\alpha'}^{-1}X)\}$. 
For an element $g\in \SL_2(\Z/p^{r+n}\Z)$, 
we have 
$V_{g^{-1}\alpha' g}^{r+n,r+m}=g^{-1}(V_{\alpha'}^{r+n,r+m})g$. 
We see that $V_{\alpha '}^{r+n,r+m}$ depends only on 
$\alpha'\bmod{p^{r+n-m}}$. 
Thus we can define $V_{\alpha''}^{r+n,r+m}$ for 
$\alpha''\in \Conj(\alpha)\subseteq \SL_2(\Z/p^{r+n-m}\Z)$.

For $\alpha=\sigma,\tau,u^{p^r}$ and their conjugates, 
we identify $V_{\alpha}^{r+n,r+m}$ 
with a free submodule of rank $2$ of $\M_2(\Z/p^{n-m}\Z)$ using the 
isomorphisms in Lemma \ref{V-str1}.

\begin{lemma}
\label{ortho}

Let $1\leq m<n$ be integers. 
Take an element $\alpha\in \SL_2(\Z/p^m\Z)$ which is conjugate to 
$\sigma$ or $\tau$. Assume $n\leq 2m$. 
Then $V_{\alpha}^{n,m}$ is the orthogonal complement of 
$\Z/p^{n-m}\Z[\alpha]$ in $\M_2(\Z/p^{n-m}\Z)$ 
with respect to the pairing 
$$\M_2(\Z/p^{n-m}\Z)\times\M_2(\Z/p^{n-m}\Z)
\longrightarrow \Z/p^{n-m}\Z\; :\; (A,B)\mapsto \tr (AB).$$

\end{lemma}

\begin{proof}

We have 
$V_{\alpha}^{n,m}=\{1+p^m(X{\alpha}^{-1}-{\alpha}^{-1}X)
\in 1+p^m\M_2(\Z/p^n\Z)\}
\cong \{X{\alpha}^{-1}-{\alpha}^{-1}X\in \M_2(\Z/p^{n-m}\Z)\}$
and it is orthogonal to $\Z/p^{n-m}\Z[\alpha]$. 
Since 
$\sharp V_{\alpha}^{n,m}=\sharp \Z/p^{n-m}\Z[\alpha]=p^{2(n-m)}$, 
they are the orthogonal complements of each other. 

\end{proof}

Note that 
$$\Z/p^{n-m}\Z[\sigma]=\left\{
\left(\begin{matrix} x&y\\ -y&x \end{matrix}\right)\right\}$$
and 
$$\Z/p^{n-m}\Z[\tau]=\left\{
\left(\begin{matrix} x&y\\ -y&x-y \end{matrix}\right)\right\}.$$

\begin{lemma}
\label{ortho-u}

Let $r\geq 0$ and $1\leq m<n$ be integers. 
Take an element 
$v=1+p^r\epsilon\in \Conj(u^{p^r})\subseteq \SL_2(\Z/p^{r+m}\Z)$. 
Assume $n\leq 2m$. 
Then $V_v^{r+n,r+m}$ is the orthogonal complement of 
$\Z/p^{n-m}\Z[\epsilon]$ in $\M_2(\Z/p^{n-m}\Z)$ 
with respect to the pairing 
$$\M_2(\Z/p^{n-m}\Z)\times\M_2(\Z/p^{n-m}\Z)
\longrightarrow \Z/p^{n-m}\Z\; :\; (A,B)\mapsto \tr (AB).$$

\end{lemma}

\begin{proof}

We have 
$V_{v}^{r+n,r+m}
=\{1+p^m(Xv^{-1}-v^{-1}X)\}
=\{1+p^m(X(1-p^r\epsilon)-(1-p^r\epsilon)X)\}
=\{1+p^{r+m}(X(-\epsilon)-(-\epsilon)X)\in 1+p^{r+m}\M_2(\Z/p^{r+n}\Z)\}
\cong \{X(-\epsilon)-(-\epsilon)X\in \M_2(\Z/p^{n-m}\Z)\}$ 
and it is orthogonal to $\Z/p^{n-m}\Z[\epsilon]$. 
Since 
$\sharp V_v^{r+n,r+m}=\sharp \Z/p^{n-m}\Z[\epsilon]=p^{2(n-m)}$, 
they are the orthogonal complements of each other. 

\end{proof}

Note that 
$$\Z/p^{n-m}\Z[
\left(\begin{matrix} 0&1\\ 0&0 \end{matrix}\right)]
=\left\{
\left(\begin{matrix} x&y\\ 0&x \end{matrix}\right)\right\}.$$

Next we study the condition for the equality $V_{\alpha}=V_{\alpha'}$, 
where $\alpha'\in \Conj(\alpha)$.

\begin{lemma}
\label{recov-sigma}

Let $n\geq 1$ be an integer. 
In $\SL_2(\Z/p^n\Z)$, we have 

  $$\left\{\left(\begin{matrix} x&y\\ -y&x \end{matrix}\right)\right\}
  \cap\Conj(\sigma)$$
\begin{equation*}
  = \begin{cases}
    \{\sigma,\,\sigma^{-1}\}\\ \text{ if $p\geq 3$},\\
    \{\sigma\}\\ \text{ if $p=2$ and $n=1$},\\
    \left\{\sigma,\,\left(\begin{matrix} 2&1\\ -1&2 \end{matrix}\right)
    \right\}\\     \text{ if $p=2$ and $n=2$},\\
    \left\{\sigma,\,
      \left(\begin{matrix} 0&1+2^{n-1}\\ -1+2^{n-1}&0 \end{matrix}\right),\,
      \left(\begin{matrix} 2^{n-1}&1\\ -1&2^{n-1} \end{matrix}\right),\,
      \left(\begin{matrix} 2^{n-1}&1+2^{n-1}\\ -1+2^{n-1}&2^{n-1} 
      \end{matrix}\right)\right\} \\
      \text{ if $p=2$ and $n\geq3$}.
  \end{cases}
\end{equation*}

\end{lemma}

\begin{corollary}
\label{sigma-p}

Let $1\leq m<n$ be integers. 
Take two elements 
$\sigma',\sigma''\in \Conj(\sigma)\subseteq \SL_2(\Z/p^m\Z)$. 
Assume $p\geq 3$ and $n\leq 2m$. 
Then $V_{\sigma'}^{n,m}=V_{\sigma''}^{n,m}$ 
holds if and only if 
$\sigma'\equiv {\sigma''}^{\pm 1}\bmod{p^{n-m}}$. 

\end{corollary}

\begin{proof}

Replacing $\sigma''$ by its conjugate, we may assume 
$\sigma''=\sigma$. 
By Lemma \ref{ortho}, we see that 
$\sigma'\bmod{p^{n-m}}$ is orthogonal to 
$V_{\sigma'}=V_{\sigma}$, and thus we have 
$\sigma'\bmod{p^{n-m}}\in 
\left\{\left(\begin{matrix}x&y\\ -y&x \end{matrix}\right)\right\}$. 
Applying the previous lemma, we get the result. 

\end{proof}

\begin{corollary}
\label{sigma-2}

Let $2\leq m<n$ be integers. 
Fix an element 
$\sigma'\in \Conj(\sigma)\subseteq \SL_2(\Z/2^m\Z)$. 
Assume $n\leq 2m$. 
Then the number of the elements 
$\sigma''\in \Conj(\sigma)\subseteq \SL_2(\Z/2^{n-m}\Z)$
satisfying 
$V_{\sigma''}^{n,m}=V_{\sigma'}^{n,m}$ 
is 
\begin{equation*}
\begin{cases}
1 &\text{if $n-m=1$},\\
2 &\text{if $n-m=2$},\\
4 &\text{if $n-m\geq 3$}.
\end{cases}
\end{equation*}

\end{corollary}

\begin{lemma}
\label{recov-tau}

Let $n\geq 1$ be an integer. 
In $\SL_2(\Z/p^n\Z)$, we have 

  $$\left\{\left(\begin{matrix} x&y\\ -y&x-y \end{matrix}\right)\right\}
  \cap\Conj(\tau)$$
\begin{equation*}
  = \begin{cases}
    \{\tau,\,\tau^{-1}\}& \text{if $p\ne 3$},\\
    \{\tau\}& \text{if $p=3$ and $n=1$},\\
    \left\{\tau,\,
      \left(\begin{matrix} 1+3^{n-1}&1-3^{n-1}\\ -1+3^{n-1}&-3^{n-1} 
      \end{matrix}\right),\,
      \left(\begin{matrix} 1-3^{n-1}&1+3^{n-1}\\ -1-3^{n-1}&3^{n-1} 
      \end{matrix}\right)\right\}
      & \text{if $p=3$ and $n\geq 2$}.
  \end{cases}
\end{equation*}

\end{lemma}

\begin{corollary}
\label{tau-p}

Let $1\leq m<n$ be integers. 
Take two elements 
$\tau',\tau''\in \Conj(\tau)\subseteq \SL_2(\Z/p^m\Z)$. 
Assume $p\ne 3$ and $n\leq 2m$. 
Suppose $m\geq 2$ if $p=2$. 
Then $V_{\tau'}^{n,m}=V_{\tau''}^{n,m}$ 
holds if and only if 
$\tau'\equiv {\tau''}^{\pm 1}\bmod{p^{n-m}}$. 

\end{corollary}

\begin{corollary}
\label{tau-3}

Let $1\leq m<n$ be integers. 
Fix an element 
$\tau'\in \Conj(\tau)\subseteq \SL_2(\Z/3^m\Z)$. 
Assume $n\leq 2m$. 
Then the number of the elements 
$\tau''\in \Conj(\tau)\subseteq \SL_2(\Z/3^{n-m}\Z)$ 
satisfying 
$V_{\tau''}^{n,m}=V_{\tau'}^{n,m}$ 
is 
\begin{equation*}
\begin{cases}
1 &\text{if $n-m=1$},\\
3 &\text{if $n-m\geq 2$}.
\end{cases}
\end{equation*}

\end{corollary}

\begin{lemma}
\label{recov-u}

Let $r\geq 0$ and $n\geq 1$ be integers. 
In $\SL_2(\Z/p^{r+n}\Z)$, we have 

$$\left\{1+p^r\left(\begin{matrix} x&y\\ 0&x \end{matrix}\right)\right\}
  \cap\Conj(u^{p^r})$$
\begin{equation*}
  = \begin{cases}
    \{u^{p^r s^2}|(p,s)=1\}& \text{if $p\geq 3$},\\
    \{u^{2^r}\}& \text{if $p=2$ and $n=1$},\\
    \left\{u^{2^r},\,\left(\begin{matrix} -1&2^r\\ 0&-1 \end{matrix}\right)
    \right\}    & \text{if $p=2$ and $n=2$},\\
    \left\{u^{2^r s^2},\,
      \left(\begin{matrix} 
      1+2^{r+n-1}&2^r s^2\\ 0&1+2^{r+n-1} \end{matrix}\right)|
      (2,s)=1
      \right\} & \text{if $p=2$ and $n\geq 3$}.
  \end{cases}
\end{equation*}

We have 
$\sharp \{u^{p^r s^2}|(p,s)=1\}
=\frac{1}{2}(p-1)p^{n-1}$ if $p\geq 3$, and 

$\sharp \left\{u^{2^r s^2},\,
      \left(\begin{matrix} 
      1+2^{r+n-1}&2^r s^2\\ 0&1+2^{r+n-1} \end{matrix}\right)|
      (2,s)=1
      \right\}
=2^{n-2}$ if $n\geq 3$. 

\end{lemma}

\begin{corollary}
\label{u-p}

Let $r\geq 0$ and $1\leq m<n$ be integers. 
Take two elements 
$v,v'\in \Conj(u^{p^r})\subseteq \SL_2(\Z/p^{r+m}\Z)$. 
Assume $p\geq 3$ and $n\leq 2m$. 
Then $V_v^{r+n,r+m}=V_{v'}^{r+n,r+m}$ 
holds if and only if 
$v\equiv {v'}^{s^2}\bmod{p^{r+n-m}}$ where $(p,s)=1$. 

\end{corollary}

\begin{proof}

Replacing $v'$ by its conjugate, we may assume $v'=u^{p^r}$. 
Write $v=1+p^r\epsilon$. 
By Lemma \ref{ortho-u}, we see that 
$\epsilon\bmod{p^{n-m}}$ is orthogonal to 
$V_v=V_{u^{p^r}}$, and thus we have 
$\epsilon\bmod{p^{n-m}}\in 
\left\{\left(
\begin{matrix}x&y \\ 0&x \end{matrix}\right)\right\}$. 
Hence we have 
$v\bmod{p^{r+n-m}}\in 
\left\{1+p^r\left(
\begin{matrix}x&y \\ 0&x \end{matrix}\right)\right\}$. 
Applying the previous lemma, we get the result. 

\end{proof}

\begin{corollary}
\label{u-2}

Let $r\geq 0$ and $3\leq m<n$ be integers. 
Fix an element 
$v\in \Conj(u^{2^r})\subseteq \SL_2(\Z/2^{r+m}\Z)$. 
Assume $n\leq 2m$. 
Then the number of the elements 
$v'\in \Conj(u^{2^r})\subseteq \SL_2(\Z/2^{r+n-m}\Z)$
satisfying 
$V_{v'}^{r+n,r+m}=V_v^{r+n,r+m}$ 
is 
\begin{equation*}
\begin{cases}
1 &\text{if $n-m=1$},\\
2 &\text{if $n-m=2$},\\
2^{n-m-2} &\text{if $n-m\geq 3$}.
\end{cases}
\end{equation*}

\end{corollary}


\section{Control of inverse images}

We control the number of elements conjugate to $\sigma,\tau,u^{p^r}$ 
contained in a slim subgroup $H$.

Let $n\geq 1$ be an integer and 
let $H$ be a subgroup of $\SL_2(\Z/p^n\Z)$. 
For an integer $1\leq s\leq n$, put 
$$H_s:=H\cap (1+p^s\M_2(\Z/p^{n}\Z))
=\Ker (\bmod{p^s}:H\longrightarrow \SL_2(\Z/p^s\Z)).$$ 
We identify 
$H/H_s$ with $H\bmod{p^s}$. 
For two integers $s,t$ with $1\leq s\leq t\leq n$ and for 
$\alpha=\sigma,\tau,u^{p^r}$, 
let 
$$f_{t,s}^{H,\alpha}:(H/H_t)\cap \Conj(\alpha)\longrightarrow
 (H/H_s)\cap \Conj(\alpha)$$ 
be the mod $p^s$ map. 
Here we assume $s>r$ when $\alpha=u^{p^r}$.

Recall that a subgroup $H\subseteq \SL_2(\Z/p^n\Z)$ 
is called a slim subgroup if 
$H\nsupseteq (1+p^{n-1}\M_2(\Z/p^n\Z))^{\det=1}$, 
equivalently $\sharp H_{n-1}\leq p^2$. 
We prepare for controlling $\sharp H\cap\Conj(\alpha)$ 
for a slim subgroup $H$.

\begin{lemma}
\label{sharpH}

Let $n\geq 2$ be an integer and 
let $H\subseteq \SL_2(\Z/p^n\Z)$ be a slim subgroup. 
Take two integers $s,t$ with $1\leq t<s\leq n$. 
Assume $t\geq 2$ if $p=2$. 
Then we have 
$\sharp H_t/H_s\leq p^{2(s-t)}$. 

\end{lemma}

\begin{proof}

Put $G=\SL_2(\Z/p^n\Z)$. 
Take an integer $i$ with $2\leq i\leq n-1$. 
Assume $i\geq 3$ if $p=2$. 
Then the $p$-th power map 
$\eta :G_{i-1}/G_i\longrightarrow G_i/G_{i+1}$ 
is surjective. 
As $\sharp G_{i-1}/G_i=\sharp G_i/G_{i+1}=p^3$, 
we see that $\eta$ is an isomorphism. 
Let $\eta':H_{i-1}/H_i\longrightarrow H_i/H_{i+1}$ 
be the $p$-th power map. 
The commutative diagram 
\begin{equation*}
  \begin{CD}
  H_{i-1}/H_i   @>\text{$\subseteq$}>>  G_{i-1}/G_i  \\
   @V\text{$\eta'$}VV                  @V\text{$\eta$}VV  \\
  H_i/H_{i+1}   @>\text{$\subseteq$}>>  G_i/G_{i+1}
  \end{CD}
\end{equation*}
shows that $\eta'$ is injective. 
By the hypothesis 
$H\nsupseteq (1+p^{n-1}\M_2(\Z/p^n\Z))^{\det=1}$, 
we have 
$H_{n-1}/H_n\subsetneq G_{n-1}/G_n$. 
Hence $\sharp H_{n-1}/H_n\leq p^2$. 
Therefore 
$\sharp H_1/H_2\leq \sharp H_2/H_3\leq \cdots \leq 
\sharp H_{n-1}/H_n\leq p^2$ 
if $p\geq 3$, while 
$\sharp H_2/H_3\leq \sharp H_3/H_4\leq \cdots \leq 
\sharp H_{n-1}/H_n\leq p^2$ 
if $p=2$. 
Consequently, we get 
$\sharp H_t/H_s=\prod_{i=t+1}^s\sharp H_{i-1}/H_i\leq p^{2(s-t)}$. 

\end{proof}

From now to the end of this section, 
we use the letter $\alpha$ to denote any of 
$\sigma,\tau,u^{p^r}$, where $r\geq 0$ be an integer. 
As usual, assume $r=0$ if $\alpha=\sigma,\tau$.

\begin{corollary}
\label{H-str}

Let $n\geq 2$ be an integer and 
let $H\subseteq \SL_2(\Z/p^{r+n}\Z)$ be a slim subgroup. 
Take two integers $s,t$ with $1\leq t<s\leq n$. 
Take an element 
$\alpha'\in \Conj(\alpha)\subseteq \SL_2(\Z/p^{r+t}\Z)$. 
Assume $s\leq 2t$. 
When $p=2$, further assume $r+t\geq 2$. 
If $H/H_{r+s}\supseteq V_{\alpha'}^{r+s,r+t}$, then 
$H_{r+t}/H_{r+s}=V_{\alpha'}^{r+s,r+t}$. 

\end{corollary}

\begin{proof}

If $H/H_{r+s}\supseteq V_{\alpha'}^{r+s,r+t}$, then 
$H_t/H_{r+s}\supseteq V_{\alpha'}^{r+s,r+t}$. 
By Lemma \ref{V-str1}, we have $\sharp V_{\alpha'}^{r+s,r+t}=p^{2(s-t)}$. 
We have seen 
$\sharp H_{r+t}/H_{r+s}\leq p^{2(s-t)}$ 
in the previous lemma. 
Thus we get 
$H_{r+t}/H_{r+s}=V_{\alpha'}^{r+s,r+t}$. 

\end{proof}

\begin{remark}
\label{rem-mod-p^n+1}

By using the diagram in the proof of Lemma \ref{sharpH}, 
we can show Lemma \ref{mod-p^n+1}. 

\end{remark}

\begin{lemma}
\label{iden}

Let $n\geq 2$ be an integer and 
let $H\subseteq \SL_2(\Z/p^{r+n}\Z)$ be a slim subgroup. 
Take three integers $s,t,i$ satisfying 
$1\leq t<s\leq n$ and $i\geq 1$. 
Take an element 
$\alpha'\in \Conj(\alpha)\subseteq \SL_2(\Z/p^{r+t}\Z)$. 
Assume $s\leq 2t$ and $s+i\leq n$. 
When $p=2$ and $\alpha=\tau$, further assume $s\leq 2t-1$. 
If 
$H_{r+t}/H_{r+s}=V_{\alpha'}^{r+s,r+t}$, 
then 
$H_{r+t+i}/H_{r+s+i}=V_{\alpha'}^{r+s+i,r+t+i}$. 

\end{lemma}

\begin{proof}

Since 
$H_{r+t}\bmod{p^{r+s}}
=\{1+p^t(X{\alpha'}^{-1}-{\alpha'}^{-1}X)|
X\in \M_2(\Z/p^{r+s}\Z)\}$, 
we have 
$H\ni 1+p^t(X{\alpha'}^{-1}-{\alpha'}^{-1}X)+p^{r+s}Y$ 
for some $Y\in \M_2(\Z/p^{r+n}\Z)$. 

Case $\alpha=\sigma$ or $\tau$. In this case $r=0$. 
Suppose $p\geq 3$ or ($p=2$ and $s\leq 2t-1$). 
Then we have 
$H\ni (1+p^t(X{\alpha'}^{-1}-{\alpha'}^{-1}X)+p^sY)^{p^i}
=1+p^{t+i}(X{\alpha'}^{-1}-{\alpha'}^{-1}X)+p^{s+i}Y'$ 
for some $Y'\in \M_2(\Z/p^n\Z)$. 
Thus 
$H_{t+i}/H_{s+i}\supseteq 
\{1+p^{t+i}(X{\alpha'}^{-1}-{\alpha'}^{-1}X)\}=V_{\alpha'}^{s+i,t+i}$. 
By Corollary \ref{H-str}, 
we have the equality 
$H_{t+i}/H_{s+i}=V_{\alpha'}^{s+i,t+i}$. 
Next suppose $\alpha=\sigma$, $p=2$ and $s=2t$. 
Notice that we are assuming $t\geq 2$ in this case. 
Replacing $\alpha'$ by its conjugate, we may assume 
$\alpha'=\sigma$. 
Then we have 
$H\ni 1+2^t\left(\begin{matrix} a&b \\ b&-a \end{matrix}\right)
+2^{2t}Y$, 
where 
$X\sigma^{-1}-\sigma^{-1}X
=\left(\begin{matrix} a&b \\ b&-a \end{matrix}\right)$. 
Taking the second power, we have 
$H\ni 1+2^{t+1}\left(
\begin{matrix} a+2^{t-1}(a^2+b^2)&b \\ b&-a+2^{t-1}(a^2+b^2) 
\end{matrix}\right)
+2^{2t+1}Y'$ 
for some $Y'\in \M_2(\Z/2^n\Z)$. 
Put 
$R(a,b):=\left(
\begin{matrix} a+2^{t-1}(a^2+b^2)&b \\ b&-a+2^{t-1}(a^2+b^2) 
\end{matrix}\right)$. 
The elements $R(a,b)$ generate 
$V_{\sigma}^{2t+1,t+1}\cong
\left\{
\left(\begin{matrix} x&y \\ y&-x \end{matrix}\right)
\in \M_2(\Z/2^t\Z)\right\}$ 
because 
$(1+2^{t-1})R(1,0)=\left(\begin{matrix} 1&0 \\ 0&-1 \end{matrix}\right)$ 
and 
$R(0,1)+2^{t-1}R(1,0)
=\left(\begin{matrix} 0&1 \\ 1&0 \end{matrix}\right)$. 
Hence 
$H_{t+1}/H_{2t+1}\supseteq V_{\sigma}^{2t+1,t+1}$, 
and this implies 
$H_{t+1}/H_{2t+1}=V_{\sigma}^{2t+1,t+1}$ 
by Corollary \ref{H-str}. 
Therefore we get 
$H_{t+i}/H_{2t+i}=V_{\sigma}^{2t+i,t+i}$ 
as we have just seen.

Case $\alpha=u^{p^r}$. 
We may assume $\alpha'=u^{p^r}$. 
Thus we have 
$H\ni 1+p^{r+t}\left(\begin{matrix} a&b \\ 0&-a \end{matrix}\right)
+p^{r+s}Y$, 
where 
$Xu^{-p^r}-u^{-p^r}X
=p^r\left(\begin{matrix} a&b \\ 0&-a \end{matrix}\right)$. 
Suppose $p\geq 3$ or ($p=2$ and $r+s\leq 2(r+t)-1$). 
Then we have 
$H\ni (1+p^{r+t}\left(\begin{matrix} a&b \\ 0&-a \end{matrix}\right)
+p^{r+s}Y)^{p^i}
=1+p^{r+t+i}\left(\begin{matrix} a&b \\ 0&-a \end{matrix}\right)
+p^{r+s+i}Y'$
for some $Y'\in \M_2(\Z/p^{r+n}\Z)$. 
Hence we get 
$H_{r+t+i}/H_{r+s+i}=V_{u^{p^r}}^{r+s+i,r+t+i}$. 
Next suppose 
$p=2$ and $r+s=2(r+t)$, so that $r=0$ and $s=2t$. 
Notice that we are assuming $t\geq 3$ in this case. 
We have 
$H\ni (1+2^t\left(\begin{matrix} a&b \\ 0&-a \end{matrix}\right)
+2^{2t}Y)^2
=1+2^{t+1}\left(
\begin{matrix} a+2^{t-1}a^2&b \\ 0&-a+2^{t-1}a^2 \end{matrix}
\right)+2^{2t+1}Y'$ 
for some $Y'\in \M_2(\Z/2^n\Z)$. 
Put 
$S(a,b):=\left(
\begin{matrix} a+2^{t-1}a^2&b \\ 0&-a+2^{t-1}a^2 \end{matrix}
\right)$. 
The elements $S(a,b)$ generate 
$V_u^{2t+1,t+1}\cong\left\{
\left(\begin{matrix} x&y \\ 0&-x \end{matrix}\right)
\in \M_2(\Z/2^t\Z)\right\}$ 
because 
$(1+2^{t-1})S(1,0)=\left(\begin{matrix} 1&0 \\ 0&-1 \end{matrix}\right)$ 
and 
$S(0,1)=\left(\begin{matrix} 0&1 \\ 0&0 \end{matrix}\right)$. 
Therefore 
$H_{t+1}/H_{2t+1}=V_u^{2t+1,t+1}$ 
and 
$H_{t+i}/H_{2t+i}=V_u^{2t+i,t+i}$. 

\end{proof}

\begin{lemma}
\label{leq-p}

Let $n\geq 2$ be an integer and 
let $H\subseteq \SL_2(\Z/p^{r+n}\Z)$ be a subgroup. 
Take two integers $t,i$ with $t\geq 1$, $i\geq 1$ and 
$t+i\leq n$. 
Take an element $\alpha'\in (H/H_{r+t})\cap\Conj(\alpha)$. 
Assume $i\leq t$. 
If $H_{r+t}/H_{r+t+i}\ne V_{\alpha'}^{r+t+i,r+t}$, then
we have 
$\sharp f_{r+t+i,r+t+1}((f_{r+t+i,r+t}^{H,\alpha})^{-1}(\alpha'))\leq p$. 

\end{lemma}

\begin{proof}

We show the lemma only when $r=0$. 
Put 
$X:=f_{t+i,t+1}((f_{t+i,t}^{H,\alpha})^{-1}(\alpha'))$. 
Suppose 
$(f_{t+i,t}^{H,\alpha})^{-1}(\alpha')\ne\emptyset$. 
Take an element 
$\widetilde{\alpha}\in (f_{t+i,t}^{H,\alpha})^{-1}(\alpha')$. 
The natural surjection 
mod $p^{t+1}:(f_{t+i,t}^{H,\alpha})^{-1}(\alpha')
\longrightarrow X$ 
induces the following commutative diagram : 
\begin{equation*}
  \begin{CD}
V_{\alpha'}^{t+i,t}\cap (H/H_{t+i})@>\text{$\bmod{p^{t+1}}$}>> 
\widetilde{\alpha}^{-1}X \\
            @V\text{$\subsetneq$}VV    @V\text{$\subseteq$}VV \\
V_{\alpha'}^{t+i,t}                @>\text{$\bmod{p^{t+1}}$}>> 
V_{\alpha'}^{t+1,t} \\
            @V\text{$\cong$}VV            @V\text{$\cong$}VV \\
(\Z/p^i\Z)^{\oplus 2}             @>\text{$\bmod{p}$}>> 
(\Z/p\Z)^{\oplus 2}.
  \end{CD}
\end{equation*}
All the horizontal maps in the above diagram are surjective. 
Since 
$V_{\alpha'}^{t+i,t}\cap (H/H_{t+i})
=V_{\alpha'}^{t+i,t}\cap (H_t/H_{t+i})$ 
is a proper subgroup of 
$V_{\alpha'}^{t+i,t}\cong (\Z/p^i\Z)^{\oplus 2}$, 
we have 
$\sharp V_{\alpha'}^{t+i,t}\cap (H/H_{t+i})\leq p^{2i-1}$. 
Since no proper subgroup of $(\Z/p^i\Z)^{\oplus 2}$ maps surjectively 
$\bmod{p}$ onto $(\Z/p\Z)^{\oplus 2}$, 
we get  
$\sharp \widetilde{\alpha}^{-1}X\leq p$. 

\end{proof}

\begin{corollary}
\label{leq-p^n-1}

Let $n\geq 2$ be an integer and 
let $H\subseteq \SL_2(\Z/p^{r+n}\Z)$ be a slim subgroup. 
Take two integers $i\geq 1$, $\delta\geq 0$ with $i+\delta\leq n$. 
Take an element 
$\alpha'\in (H/H_{r+i+\delta})\cap\Conj(\alpha)$. 
Assume $2i+\delta\leq n$. 
When $p=2$ and $\alpha=\tau$, further assume $\delta\geq 1$. 
If $H_{r+n-i}/H_{r+n}\ne V_{\alpha'}^{r+n,r+n-i}$, 
then we have 
$\sharp (f_{r+n,r+i+\delta}^{H,\alpha})^{-1}(\alpha')\leq p^{n-1-\delta}$. 

\end{corollary}

\begin{proof}

We show the corollary only when $r=0$. 
For an integer $t$ with $i\leq t\leq n$, put 
$X_t:=f_{n,t}((f_{n,i}^{H,\alpha})^{-1}(\alpha'))
\subseteq (f_{t,i}^{H,\alpha})^{-1}(\alpha')$. 
For $i\leq t\leq s\leq n$, let 
$f_{X_s,X_t}:X_s\longrightarrow X_t$ 
be the $\bmod{p^t}$ reduction map. 
Note that $f_{X_s,X_t}$ is surjective. 
For each element $\alpha''\in X_{n-i+1}$, we have 
$\sharp f_{X_n,X_{n-i+1}}^{-1}(\alpha'')\leq p^{2(i-1)}$ 
by Lemma \ref{V-str1}. 
Thus we have 
$\sharp X_n\leq p^{2(i-1)}\sharp X_{n-i+1}$. 
For $i+\delta\leq t\leq n-i$ and $\alpha''\in X_t$, 
we show $\sharp f_{X_{t+1},X_t}^{-1}(\alpha'')\leq p$. 
We have 
$f_{X_{t+1},X_t}^{-1}(\alpha'')
=f_{t+i,t+1}(f_{X_{t+i},X_t}^{-1}(\alpha''))$. 
By the assumption $H_{n-i}/H_n\ne V_{\alpha'}^{n,n-i}$ 
and by Lemma \ref{iden}, we have 
$H_t/H_{t+i}\ne V_{\alpha''}^{t+i,t}$. 
Since  
$f_{X_{t+i},X_t}^{-1}(\alpha'')\subseteq 
(f_{t+i,t}^{H,\alpha})^{-1}(\alpha'')$, 
we have the inequality 
$\sharp f_{t+i,t+1}(f_{X_{t+i},X_t}^{-1}(\alpha''))\leq p$ 
by Lemma \ref{leq-p}. 
Hence we have 
$\sharp f_{X_{t+1},X_t}^{-1}(\alpha'')\leq p$ 
for $i+\delta\leq t\leq n-i$ and $\alpha''\in X_t$. 
Thus we get 
$\sharp X_{t+1}\leq p\sharp X_t$ for $i+\delta\leq t\leq n-i$. 
Consequently we have 
$\sharp (f_{n,i}^{H,\alpha})^{-1}(\alpha')
=\sharp X_n\leq p^{2(i-1)}\cdot p^{n-2i+1-\delta}\sharp X_i
\leq p^{n-1-\delta}$. 

\end{proof}

Let $n\geq 2$ be an integer and let
$H\subseteq \SL_2(\Z/p^{r+n}\Z)$ be a slim subgroup. 
Put 
\begin{equation*}
l:=
  \begin{cases}
    \frac{n}{2} & \text{if $n$ is even}, \\
    \frac{n-1}{2} & \text{if $n$ is odd}.
  \end{cases}
\end{equation*}
Define a decreasing sequence 
$$Y_0\supseteq Y_1\supseteq \cdots \supseteq Y_i\supseteq Y_{i+1}
\supseteq\cdots\supseteq Y_l$$ 
by 
\begin{equation*}
Y_i:=
  \begin{cases}
    H\cap\Conj(\alpha) & \text{if $i=0$}, \\
    \{\alpha'\in H\cap\Conj(\alpha)|H_{r+n-i}=V_{\alpha'}^{r+n,r+n-i}\}
    & \text{if $1\leq i\leq l$}.
  \end{cases}
\end{equation*}

When $p\geq 3$, we use $Y_i$ for $0\leq i\leq l$; 
when $p=2$, we use $Y_i$ only for 
\begin{equation*}
  \begin{cases}
    i=0,1 & \text{if $\alpha=\sigma$ and $3\leq n\leq 5$}, \\
    0\leq i\leq l & \text{if $\alpha=\sigma$ and $n\geq 6$}, \\
    i=0,\ 3\leq i\leq l & \text{if $\alpha=u^{2^r}$ and $n\geq 6$},
  \end{cases}
\end{equation*}
so that the hypothesis in Lemma \ref{V-str1} is satisfied.

\begin{lemma}
\label{setminus}

For three integers $i,i',i''$ satisfying $0\leq i'<i\leq i''\leq l$, 
we have 
$$(Y_{i'}\setminus Y_i)\bmod{p^{r+i''}}=
(Y_{i'}\bmod{p^{r+i''}})\setminus (Y_i\bmod{p^{r+i''}}).$$

\end{lemma}

\begin{proof}

The inclusion ``$\supseteq$'' is trivial. 
Now we show the other inclusion ``$\subseteq$''. 
Take an element 
$\alpha'\bmod{p^{r+i''}}\in (Y_{i'}\setminus Y_i)\bmod{p^{r+i''}}$ 
where $\alpha'\in Y_{i'}\setminus Y_i$. 
Assume $\alpha'\bmod{p^{r+i''}}\in Y_i\bmod{p^{r+i''}}$. 
Then we find an element $\alpha''\in Y_i$ such that 
$\alpha'\equiv \alpha''\bmod{p^{r+i''}}$. 
Since $i\leq i''$, we have 
$V_{\alpha'}^{r+n,r+n-i}=V_{\alpha''}^{r+n,r+n-i}=H_i$. 
But this contradicts the assumption
$\alpha'\not\in Y_i$. 
Hence we get the inclusion ``$\subseteq$''. 

\end{proof}

\begin{proposition}
\label{estimate}

Let $n\geq 2$ be an integer and let
$H\subseteq \SL_2(\Z/p^{r+n}\Z)$ be a slim subgroup. 
Take an integer $s$ satisfying $1\leq s\leq l$. 
Then we have 
\begin{align*}
&\sharp H\cap\Conj(\alpha)\\
\leq 
&(p^{2(n-l)}-p^{n-1})\sharp(Y_l\bmod{p^{r+l}})
+(p^2-1)p^{n-1}\sum_{i=s}^{l-1}\sharp(Y_i\bmod{p^{r+i}})\\
&+p^{n-1}\sharp(H/H_{r+s})\cap\Conj(\alpha).
\end{align*}

\end{proposition}

\begin{proof}

For simplicity, we show the inequality only when $r=0$. 
Decompose 
$$\sharp Y_0
=\sharp Y_l+\sum_{i=s+1}^l\sharp(Y_{i-1}\setminus Y_i)
+\sharp(Y_0\setminus Y_s).$$
We have 
$$\sharp Y_l\leq p^{2(n-l)}\sharp(Y_l\bmod{p^l})$$
by Lemma \ref{V-str1}. 
We see 
$$\sharp(Y_{i-1}\setminus Y_i)\leq 
p^{n-1}\sharp((Y_{i-1}\setminus Y_i)\bmod{p^i})
=p^{n-1}(\sharp(Y_{i-1}\bmod{p^i})-\sharp(Y_i\bmod{p^i}))$$
and 
$$\sharp(Y_0\setminus Y_s)\leq 
p^{n-1}\sharp((Y_0\setminus Y_s)\bmod{p^s})
=p^{n-1}(\sharp(Y_0\bmod{p^s})-\sharp(Y_s\bmod{p^s}))$$
by the definition of $Y_i$, Corollary \ref{leq-p^n-1} 
and Lemma \ref{setminus}. 
Hence 
\begin{align*}
\sharp Y_0 \leq 
&p^{2(n-l)}\sharp(Y_l\bmod{p^l})
+p^{n-1}\sum_{i=s+1}^l(\sharp(Y_{i-1}\bmod{p^i})-\sharp(Y_i\bmod{p^i}))\\
&+p^{n-1}(\sharp(Y_0\bmod{p^s})-\sharp(Y_s\bmod{p^s}))\\
=&(p^{2(n-l)}-p^{n-1})\sharp(Y_l\bmod{p^l})
+p^{n-1}\sum_{i=s}^{l-1}
(\sharp(Y_i\bmod{p^{i+1}})-\sharp(Y_i\bmod{p^i}))\\
&+p^{n-1}\sharp(Y_0\bmod{p^s}).
\end{align*}
By Lemma \ref{V-str1}, we see 
$$\sharp(Y_i\bmod{p^{i+1}})-\sharp(Y_i\bmod{p^i})
\leq (p^2-1)\sharp(Y_i\bmod{p^i}).$$
Since 
$$Y_0\bmod{p^s}=(H/H_s)\cap\Conj(\alpha),$$ 
we get the desired inequality. 

\end{proof}

We control the number of elements conjugate to $\sigma$ 
in a slim subgroup $H\subseteq \SL_2(\Z/p^{n}\Z)$ 
when $p\geq 3$.

For $p\geq 3$, 
define a sequence $\{a(\sigma,p)_n\}_{n\geq 2}$ as follows:
\begin{equation*}
a(\sigma,p)_n:=
2p^{2(n-l)}+2(l-1)(p^2-1)p^{n-1},
\end{equation*}
where $n=2l$ or $2l+1$.

\begin{corollary}
\label{asp}

Let $n\geq 2$ be an integer and let
$H\subseteq \SL_2(\Z/p^{n}\Z)$ be a slim subgroup. 
Assume $p\geq 3$. 
Then we have 
$$\sharp H\cap\Conj(\sigma)\leq a(\sigma,p)_n
+p^{n-1}(\sharp(H/H_1)\cap\Conj(\sigma)-2).$$

\end{corollary}

\begin{proof}

By Corollary \ref{sigma-p}, we have 
$$\sharp(Y_i\bmod{p^i})\leq 2$$ 
for $1\leq i\leq l$. 
Applying Proposition \ref{estimate} (put $s=1$), we get 
\begin{align*}
&\sharp H\cap\Conj(\sigma) \\
\leq 
&(p^{2(n-l)}-p^{n-1})\cdot 2+(p^2-1)p^{n-1}\cdot 2(l-1)
+p^{n-1}\sharp(H/H_1)\cap\Conj(\sigma)\\
=&a(\sigma,p)_n+p^{n-1}(\sharp(H/H_1)\cap\Conj(\sigma)-2).
\end{align*}

\end{proof}

We control the number of elements conjugate to $\tau$ 
in a slim subgroup $H\subseteq \SL_2(\Z/p^{n}\Z)$ 
when $p\geq 3$. 

For $p\geq 5$, define a sequence $\{a(\tau,p)_n\}_{n\geq 2}$ 
by 
$$a(\tau,p)_n:=a(\sigma,p)_n.$$

\begin{corollary}
\label{atp}

Let $n\geq 2$ be an integer and let
$H\subseteq \SL_2(\Z/p^{n}\Z)$ be a slim subgroup. 
Assume $p\geq 5$. 
Then we have 
$$\sharp H\cap\Conj(\tau)\leq a(\tau,p)_n
+p^{n-1}(\sharp(H/H_1)\cap\Conj(\tau)-2).$$

\end{corollary}

\begin{proof}

The same as the proof of Corollary \ref{asp}.

\end{proof}

Define a sequence $\{a(\tau,3)_n\}_{n\geq 2}$ as follows:
\begin{equation*}
a(\tau,3)_n:=
\begin{cases}
3^2 &\text{if $n=2$},\\
(4n-11)\cdot 3^n &\text{if $n=2l\geq 4$},\\
(4n-9)\cdot 3^n &\text{if $n=2l+1$}.
\end{cases}
\end{equation*}

\begin{corollary}
\label{at3}

Let $n\geq 2$ be an integer and let
$H\subseteq \SL_2(\Z/3^n\Z)$ be a slim subgroup. 
Then we have 
$$\sharp H\cap\Conj(\tau)\leq a(\tau,3)_n
+3^{n-1}(\sharp(H/H_1)\cap\Conj(\tau)-1).$$

\end{corollary}

\begin{proof}

By Corollary \ref{tau-3}, we have 
\begin{equation*}
\sharp(Y_i\bmod{p^i})\leq 
\begin{cases}
  1 &\text{if $i=1$},\\
  3 &\text{if $i\geq 2$}.
\end{cases}
\end{equation*}
Applying Proposition \ref{estimate} (put $s=1$), we get 
\begin{align*}
&\sharp H\cap\Conj(\tau)\\
\leq &
  \begin{cases}
    (3^{2(n-l)}-3^{n-1})\cdot 1+3^{n-1}\sharp(H/H_1)\cap\Conj(\tau)\\
    \text{\ if $l=1$},\\ 
    (3^{2(n-l)}-3^{n-1})\cdot 3+(3^2-1)\cdot 3^{n-1}\cdot 3(l-2)
    +(3^2-1)\cdot 3^{n-1}\cdot 1\\
    +3^{n-1}\sharp(H/H_1)\cap\Conj(\tau)\\
    \text{\ if $l\geq 2$}
  \end{cases}
\\
=&a(\tau,3)_n+3^{n-1}(\sharp(H/H_1)\cap\Conj(\tau)-1).
\end{align*}

\end{proof}

We control the number of elements conjugate to $u^{p^r}$ 
in a slim subgroup $H\subseteq \SL_2(\Z/p^{r+n}\Z)$.

For $p\geq 3$, 
define a sequence $\{a(u,p)_n\}_{n\geq 2}$ as follows:
\begin{equation*}
a(u,p)_n:=
\begin{cases}
\frac{1}{2}(p-1)(2\cdot p^{3l-1}-p^n) &\text{if $n=2l$},\\
\frac{1}{2}(p-1)(p^{3l+1}+p^{3l}-p^n) &\text{if $n=2l+1$}.
\end{cases}
\end{equation*}

\begin{corollary}
\label{aup}

Let $r\geq 0$, $n\geq 2$ be integers and 
let $H\subseteq \SL_2(\Z/p^{r+n}\Z)$ be a slim subgroup. 
Assume $p\geq 3$. 
Then we have 
$$\sharp H\cap\Conj(u^{p^r})\leq a(u,p)_n
+p^{n-1}(\sharp(H/H_{r+1})\cap\Conj(u^{p^r})-\frac{1}{2}(p-1)).$$

\end{corollary}

\begin{proof}

By Corollary \ref{u-p}, we have 
$$\sharp(Y_i\bmod{p^{r+i}})\leq \frac{1}{2}(p-1)p^{i-1}$$ 
for $1\leq i\leq l$. 
Applying Proposition \ref{estimate} (put $s=1$), we get 
\begin{align*}
&\sharp H\cap\Conj(u^{p^r})\\
\leq 
&(p^{2(n-l)}-p^{n-1})\cdot\frac{1}{2}(p-1)p^{l-1}
+(p^2-1)p^{n-1}\sum_{i=1}^{l-1}\frac{1}{2}(p-1)p^{i-1}\\
&+p^{n-1}\sharp(H/H_{r+1})\cap\Conj(u^{p^r})\\
=&a(u,p)_n
+p^{n-1}(\sharp(H/H_{r+1})\cap\Conj(u^{p^r})-\frac{1}{2}(p-1)).
\end{align*}

\end{proof}

Define a sequence $\{a(u,2)_n\}_{n\geq 6}$ as follows:
\begin{equation*}
a(u,2)_n:=
\begin{cases}
2^{3l-1}-2^{n+1} &\text{if $n=2l$},\\
3\cdot 2^{3l-1}-2^{n+1} &\text{if $n=2l+1$}.
\end{cases}
\end{equation*}

\begin{corollary}
\label{au2}

Let $r\geq 0$, $n\geq 6$ be integers and let
$H\subseteq \SL_2(\Z/2^{r+n}\Z)$ be a slim subgroup. 
Then we have 
$$\sharp H\cap\Conj(u^{2^r})\leq a(u,2)_n
+2^{n-1}(\sharp(H/H_{r+3})\cap\Conj(u^{2^r})-2).$$

\end{corollary}

\begin{proof}

By Corollary \ref{u-2}, we have 
$$\sharp(Y_i\bmod{2^{r+i}})\leq 2^{i-2}$$ 
for $3\leq i\leq l$. 
Applying Proposition \ref{estimate} (put $s=3$), we get 
\begin{align*}
&\sharp H\cap\Conj(u^{2^r})\\
\leq &(2^{2(n-l)}-2^{n-1})\cdot 2^{l-2}
+(2^2-1)\cdot 2^{n-1}\sum_{i=3}^{l-1}2^{i-2}
+2^{n-1}\sharp(H/H_{r+3})\cap\Conj(u^{2^r})\\
=&a(u,2)_n+2^{n-1}(\sharp(H/H_{r+3})\cap\Conj(u^{2^r})-2).
\end{align*}

\end{proof}

We control the number of elements conjugate to $\sigma$ 
in a slim subgroup $H\subseteq \SL_2(\Z/2^{n}\Z)$.

Define a sequence $\{a(\sigma,2)_n\}_{n\geq 3}$ as follows:
\begin{equation*}
a(\sigma,2)_n:=
\begin{cases}
  2^3 &\text{if $n=3$},\\
  2^5 &\text{if $n=4$},\\
  3(l-2)\cdot 2^{n+1} &\text{if $n=2l\geq 6$},\\
  (3l-4)\cdot 2^{n+1} &\text{if $n=2l+1\geq 5$}.
\end{cases}
\end{equation*}

\begin{proposition}
\label{as2}

Let $n\geq 3$ be an integer and let
$H\subseteq \SL_2(\Z/2^n\Z)$ be a slim subgroup. 
Then we have 
$$\sharp H\cap\Conj(\sigma)
\leq a(\sigma,2)_n+2^{n-2}(\sharp(H/H_2)\cap\Conj(\sigma)-2).$$

\end{proposition}

\begin{proof}

Decompose 
$$\sharp Y_0=\sharp Y_1+\sharp(Y_0\setminus Y_1).$$
We have 
$$\sharp Y_1\leq 2^{2(n-2)}\sharp(Y_1\bmod{2^2})$$ 
by Lemma \ref{V-str1}. 
We see 
$$\sharp(Y_0\setminus Y_1)\leq 
2^{n-2}\sharp((Y_0\setminus Y_1)\bmod{2^2})
=2^{n-2}(\sharp(Y_0\bmod{2^2})-\sharp(Y_1\bmod{2^2}))$$
by the definition of $Y_1$, Corollary \ref{leq-p^n-1} 
and Lemma \ref{setminus}. 
Hence 
$$\sharp Y_0\leq 
(2^{2(n-2)}-2^{n-2})\sharp(Y_1\bmod{2^2})+2^{n-2}\sharp(Y_0\bmod{2^2}).$$
As 
$$\sharp(Y_1\bmod{2^2})\leq 2$$ 
by Corollary \ref{sigma-2}
and Remark \ref{inv=2}, 
we get 
$$\sharp Y_0\leq 
(2^{2(n-2)}-2^{n-2})\cdot 2+2^{n-2}\sharp(Y_0\bmod{2^2})
=2^{2n-3}+2^{n-2}(\sharp(H/H_2)\cap\Conj(\sigma)-2).$$
This shows the desired inequality for $n=3,4,5$.

Assume $n\geq 6$, so that we have $l\geq 3$. 
Decompose 
$$\sharp Y_0=\sharp Y_l+\sum_{i=1}^l\sharp(Y_{i-1}\setminus Y_i).$$
We have 
$$\sharp Y_l\leq 2^{2(n-l)}\sharp(Y_l\bmod{2^l})$$
by Lemma \ref{V-str1}. 
We see 
\begin{align*}
&\sharp(Y_{i-1}\setminus Y_i)\\
\leq &
  \begin{cases}
    2^{n-1}\sharp((Y_{i-1}\setminus Y_i)\bmod{2^i})
    =2^{n-1}(\sharp(Y_{i-1}\bmod{2^i})-\sharp(Y_i\bmod{2^i}))
    & \text{if $i\geq 2$},\\
    2^{n-2}\sharp((Y_0\setminus Y_1)\bmod{2^2})
    =2^{n-2}(\sharp(Y_0\bmod{2^2})-\sharp(Y_1\bmod{2^2}))
    & \text{if $i=1$}
  \end{cases}
\end{align*}
by the definition of $Y_i$, Corollary \ref{leq-p^n-1} and 
Lemma \ref{setminus}. 
Hence 
\begin{align*}
\sharp Y_0\leq 
&2^{2(n-l)}\sharp(Y_l\bmod{2^l})
+2^{n-1}\sum_{i=2}^l 
(\sharp(Y_{i-1}\bmod{2^i})-\sharp(Y_i\bmod{2^i}))\\
&+2^{n-2}(\sharp(Y_0\bmod{2^2})-\sharp(Y_1\bmod{2^2}))\\
=&(2^{2(n-l)}-2^{n-1})\sharp(Y_l\bmod{2^l})
+2^{n-1}\sum_{i=2}^{l-1}
(\sharp(Y_i\bmod{2^{i+1}})-\sharp(Y_i\bmod{2^{i}}))\\
&+2^{n-2}\sharp(Y_1\bmod{2^2})
+2^{n-2}\sharp(Y_0\bmod{2^2}).
\end{align*}
By Lemma \ref{V-str1}, we see 
$$\sharp(Y_i\bmod{2^{i+1}})-\sharp(Y_i\bmod{2^i})
\leq (2^2-1)\sharp(Y_i\bmod{2^i})
=3\sharp(Y_i\bmod{2^i})$$
for $2\leq i\leq l-1$. 
Thus we have 
\begin{align*}
\sharp Y_0\leq 
&(2^{2(n-l)}-2^{n-1})\sharp(Y_l\bmod{2^l})
+3\cdot 2^{n-1}(\sum_{i=2}^{l-1}\sharp(Y_i\bmod{2^i}))\\
&+2^{n-2}\sharp(Y_1\bmod{2^2})
+2^{n-2}\sharp(Y_0\bmod{2^2}).
\end{align*}
As 
\begin{equation*}
\sharp(Y_i\bmod{2^i})\leq 
\begin{cases}
  2 &\text{if $i=2$},\\
  4 &\text{if $i\geq 3$}
\end{cases}
\end{equation*}
and 
$$\sharp(Y_1\bmod{2^2})\leq 2$$ 
by Corollary \ref{sigma-2} and Remark \ref{inv=2}, 
we get 
\begin{align*}
\sharp Y_0 \leq 
&(2^{2(n-l)}-2^{n-1})\cdot 4+3\cdot 2^{n-1}\cdot 4(l-3)
+3\cdot 2^{n-1}\cdot 2+2^{n-2}\cdot 2\\
&+2^{n-2}\sharp(Y_0\bmod{2^2})\\
=&a(\sigma,2)_n+2^{n-2}(\sharp(H/H_2)\cap\Conj(\sigma)-2).
\end{align*}

\end{proof}

We also use a slightly different way to control 
$\sharp H\cap\Conj(\tau)$ and 
$\sharp H\cap\Conj(u^{p^r})$ when $p=2$. 

Let $n\geq 2$ be an integer and let
$H\subseteq \SL_2(\Z/2^{r+n}\Z)$ be a slim subgroup. 
Put 
\begin{equation*}
l':=
  \begin{cases}
    \frac{n}{2} & \text{if $n$ is even}, \\
    \frac{n+1}{2} & \text{if $n$ is odd}.
  \end{cases}
\end{equation*}
Define a decreasing sequence 
$$Z_0\supseteq Z_1\supseteq \cdots \supseteq Z_i\supseteq Z_{i+1}
\supseteq\cdots\supseteq Z_{l'-1}$$ 
by 
\begin{equation*}
Z_i:=
  \begin{cases}
    H\cap\Conj(\alpha) & \text{if $i=0$}, \\
    \{\alpha'\in H\cap\Conj(\alpha)|H_{r+n-i}=V_{\alpha'}^{r+n,r+n-i}\}
    & \text{if $1\leq i\leq l'-1$}.
  \end{cases}
\end{equation*}

When $p=2$, we use $Z_i$ only for 
\begin{equation*}
  \begin{cases}
    i=0,\ 2\leq i\leq l'-1 & \text{if $\alpha=\tau$ and $n\geq 5$}, \\
    i=0,1 & \text{if $\alpha=u^{2^r}$ and $4\leq n\leq 6$}, \\
    i=0,1,\ 3\leq i\leq l'-1 & \text{if $\alpha=u^{2^r}$ and $n\geq 7$},
  \end{cases}
\end{equation*}
so that the hypothesis in Lemma \ref{V-str1} is satisfied.

\begin{lemma}
\label{setminus2}

For three integers $i,i',i''$ satisfying $0\leq i'<i\leq i''\leq l$, 
we have 
$$(Z_{i'}\setminus Z_i)\bmod{p^{r+i''}}=
(Z_{i'}\bmod{p^{r+i''}})\setminus (Z_i\bmod{p^{r+i''}}).$$

\end{lemma}

\begin{proof}

The same as the proof of Lemma \ref{setminus}. 

\end{proof}

We control the number of elements conjugate to $\tau$ 
in a slim subgroup $H\subseteq \SL_2(\Z/2^n\Z)$.

Define a sequence $\{a(\tau,2)_n\}_{n\geq 5}$ as follows:
\begin{equation*}
a(\tau,2)_n:=
\begin{cases}
(3l'-5)\cdot 2^{n+1} &\text{if $n=2l'$},\\
(3l'-7)\cdot 2^{n+1} &\text{if $n=2l'-1$}.
\end{cases}
\end{equation*}

\begin{proposition}
\label{at2}

Let $n\geq 5$ be an integer and let
$H\subseteq \SL_2(\Z/2^n\Z)$ be a slim subgroup. 
Then we have 
$$\sharp H\cap\Conj(\tau)\leq a(\tau,2)_n
+2^{n-2}(\sharp(H/H_3)\cap\Conj(\tau)-8).$$

\end{proposition}

\begin{proof}

Decompose 
$$\sharp Z_0=\sharp Z_{l'-1}
+\sum_{i=3}^{l'-1}\sharp(Z_{i-1}\setminus Z_i)
+\sharp(Z_0\setminus Z_2).$$
We have 
$$\sharp Z_{l'-1}\leq 2^{2(n-l')}\sharp(Z_{l'-1}\bmod{2^{l'}})$$
by Lemma \ref{V-str1}. 
By the definition of $Z_i$, Corollary \ref{leq-p^n-1} 
and Lemma \ref{setminus2}, we see 
\begin{align*}
\sharp(Z_{i-1}\setminus Z_i)\leq 
&2^{n-2}\sharp((Z_{i-1}\setminus Z_i)\bmod{2^{i+1}})\\
=&2^{n-2}(\sharp(Z_{i-1}\bmod{2^{i+1}})-\sharp(Z_i\bmod{2^{i+1}}))
\end{align*}
for $3\leq i\leq l'-1$
and 
$$\sharp(Z_0\setminus Z_2)\leq 
2^{n-2}\sharp((Z_0\setminus Z_2)\bmod{2^3})
=2^{n-2}(\sharp(Z_0\bmod{2^3})-\sharp(Z_2\bmod{2^3})).$$
Hence 
\begin{align*}
\sharp Z_0\leq 
&2^{2(n-l')}\sharp(Z_{l'-1}\bmod{2^{l'}})
+2^{n-2}\sum_{i=3}^{l'-1}
(\sharp(Z_{i-1}\bmod{2^{i+1}})-\sharp(Z_i\bmod{2^{i+1}}))\\
&+2^{n-2}(\sharp(Z_0\bmod{2^3})-\sharp(Z_2\bmod{2^3}))\\
=&(2^{2(n-l')}-2^{n-2})\sharp(Z_{l'-1}\bmod{2^{l'}})\\
&+2^{n-2}\sum_{i=2}^{l'-2}
(\sharp(Z_i\bmod{2^{i+2}})-\sharp(Z_i\bmod{2^{i+1}}))
+2^{n-2}\sharp(Z_0\bmod{2^3}).
\end{align*}
By Lemma \ref{V-str1}, we see 
$$\sharp(Z_i\bmod{2^{i+2}})-\sharp(Z_i\bmod{2^{i+1}})
\leq (2^2-1)\sharp(Z_i\bmod{2^{i+1}})
=3\sharp(Z_i\bmod{2^{i+1}})$$ 
for $2\leq i\leq l'-2$. 
Therefore
\begin{align*}
\sharp Z_0\leq 
&(2^{2(n-l')}-2^{n-2})\sharp(Z_{l'-1}\bmod{2^{l'}})
+3\cdot 2^{n-2}\sum_{i=2}^{l'-2}\sharp(Z_i\bmod{2^{i+1}})\\
&+2^{n-2}\sharp(Z_0\bmod{2^3}).
\end{align*}
By Corollary \ref{tau-p} and Lemma \ref{V-str1}, 
we have 
$$\sharp(Z_i\bmod{2^{i+1}})\leq 2^3$$ 
for $2\leq i\leq l'-1$. 
Consequently, we get 
\begin{align*}
\sharp Z_0 &\leq 
(2^{2(n-l')}-2^{n-2})\cdot 2^3
+3\cdot 2^{n-2}(l'-3)\cdot 2^3
+2^{n-2}\sharp(Z_0\bmod{2^3})\\
&=a(\tau,2)_n+2^{n-2}(\sharp(H/H_3)\cap\Conj(\tau)-8).
\end{align*}

\end{proof}

We find another control of 
the number of elements conjugate to $u^{2^r}$ 
in a slim subgroup $H\subseteq \SL_2(\Z/2^{r+n}\Z)$. 
It is better than a previous one for a small $n$.

Define a sequence $\{b(u,2)_n\}_{n\geq 4}$ as follows:
\begin{equation*}
b(u,2)_n:=
\begin{cases}
3\cdot 2^{3l'-2}-2^{n+1} &\text{if $n=2l'$},\\
2^{3l'-2}-2^{n+1} &\text{if $n=2l'-1$}.
\end{cases}
\end{equation*}

\begin{proposition}
\label{bu2}

Let $r\geq 0$, $n\geq 4$ be integers and let
$H\subseteq \SL_2(\Z/2^{r+n}\Z)$ be a slim subgroup. 
Then we have 
$$\sharp H\cap\Conj(u^{2^r})\leq b(u,2)_n
+2^{n-3}(\sharp(H/H_{r+3})\cap\Conj(u^{2^r})-4).$$

\end{proposition}

\begin{proof}

Decompose 
$$\sharp Z_0=\sharp Z_1+\sharp(Z_0\setminus Z_1).$$
We have 
$$\sharp Z_1\leq 2^{2(n-3)}\sharp(Z_1\bmod{2^{r+3}})$$
by Lemma \ref{V-str1}. 
We see 
$$\sharp(Z_0\setminus Z_1)\leq 
2^{n-3}\sharp((Z_0\setminus Z_1)\bmod{2^{r+3}})
=2^{n-3}(\sharp(Z_0\bmod{2^{r+3}})-\sharp(Z_1\bmod{2^{r+3}}))$$
by the definition of $Z_1$, Corollary \ref{leq-p^n-1} 
and Lemma \ref{setminus2}. 
Hence 
$$\sharp Z_0\leq 
(2^{2(n-3)}-2^{n-3})\sharp(Z_1\bmod{2^{r+3}})
+2^{n-3}\sharp(Z_0\bmod{2^{r+3}}).$$
As 
$$\sharp(Z_1\bmod{2^{r+3}})\leq 2^2$$
by Corollary \ref{u-2} and Remark \ref{inv=2}, we get 
\begin{align*}
\sharp Z_0 &\leq 
(2^{2(n-3)}-2^{n-3})\cdot 2^2+2^{n-3}\sharp(Z_0\bmod{2^{r+3}})\\
&=2^{2n-4}+2^{n-3}(\sharp(H/H_{r+3})\cap\Conj(u^{2^r})-4).
\end{align*}
This shows the desired inequality for $n=4,5,6$.

Assume $n\geq 7$, so that 
we have $l'\geq 4$. 
Decompose 
$$\sharp Z_0=\sharp Z_{l'-1}
+\sum_{i=4}^{l'-1}\sharp(Z_{i-1}\setminus Z_i)
+\sharp(Z_1\setminus Z_3)
+\sharp(Z_0\setminus Z_1).$$
We have 
$$\sharp Z_{l'-1}\leq 2^{2(n-l')}\sharp(Z_{l'-1}\bmod{2^{r+l'}})$$
by Lemma \ref{V-str1}. 
By the definition of $Z_i$, Corollary \ref{leq-p^n-1} 
and Lemma \ref{setminus2}, we see 
\begin{align*}
\sharp(Z_{i-1}\setminus Z_i)\leq 
&2^{n-2}\sharp((Z_{i-1}\setminus Z_i)\bmod{2^{r+i+1}})\\
=&2^{n-2}(\sharp(Z_{i-1}\bmod{2^{r+i+1}})-\sharp(Z_i\bmod{2^{r+i+1}}))
\end{align*}
for $4\leq i\leq l'-1$, 
$$\sharp(Z_1\setminus Z_3)\leq 
2^{n-2}\sharp((Z_1\setminus Z_3)\bmod{2^{r+4}})
=2^{n-2}(\sharp(Z_1\bmod{2^{r+4}})-\sharp(Z_3\bmod{2^{r+4}}))$$
and 
$$\sharp(Z_0\setminus Z_1)\leq 
2^{n-3}\sharp((Z_0\setminus Z_1)\bmod{2^{r+3}})
=2^{n-3}(\sharp(Z_0\bmod{2^{r+3}})-\sharp(Z_1\bmod{2^{r+3}})).$$
Hence 
\begin{align*}
\sharp Z_0
\leq 
&2^{2(n-l')}\sharp(Z_{l'-1}\bmod{2^{r+l'}})
+2^{n-2}\sum_{i=4}^{l'-1}
(\sharp(Z_{i-1}\bmod{2^{r+i+1}})-\sharp(Z_i\bmod{2^{r+i+1}}))\\
&+2^{n-2}(\sharp(Z_1\bmod{2^{r+4}})-\sharp(Z_3\bmod{2^{r+4}}))\\
&+2^{n-3}(\sharp(Z_0\bmod{2^{r+3}})-\sharp(Z_1\bmod{2^{r+3}}))\\
=&(2^{2(n-l')}-2^{n-2})\sharp(Z_{l'-1}\bmod{2^{r+l'}})\\
&+2^{n-2}\sum_{i=3}^{l'-2}
(\sharp(Z_i\bmod{2^{r+i+2}})-\sharp(Z_i\bmod{2^{r+i+1}}))\\
&+2^{n-3}(2\sharp(Z_1\bmod{2^{r+4}})-\sharp(Z_1\bmod{2^{r+3}}))
+2^{n-3}\sharp(Z_0\bmod{2^{r+3}}).
\end{align*}
By Lemma \ref{V-str1}, we see 
\begin{align*}
\sharp(Z_i\bmod{2^{r+i+2}})-\sharp(Z_i\bmod{2^{r+i+1}})
\leq &(2^2-1)\sharp(Z_i\bmod{2^{r+i+1}})\\
=&3\sharp(Z_i\bmod{2^{r+i+1}})
\end{align*}
for $3\leq i\leq l'-2$ and 
$$2\sharp(Z_1\bmod{2^{r+4}})-\sharp(Z_1\bmod{2^{r+3}})
\leq (2\cdot 2^2-1)\sharp(Z_1\bmod{2^{r+3}})
=7\sharp(Z_1\bmod{2^{r+3}}).$$
Therefore
\begin{align*}
\sharp Z_0\leq 
&(2^{2(n-l')}-2^{n-2})\sharp(Z_{l'-1}\bmod{2^{r+l'}})
+3\cdot 2^{n-2}\sum_{i=3}^{l'-2}\sharp(Z_i\bmod{2^{r+i+1}})\\
&+7\cdot 2^{n-3}\sharp(Z_1\bmod{2^{r+3}})
+2^{n-3}\sharp(Z_0\bmod{2^{r+3}}).
\end{align*}
By Corollary \ref{u-2} and Lemma \ref{V-str1}, 
we have 
$$\sharp(Z_i\bmod{2^{r+i+1}})\leq 2^i$$ 
for $3\leq i\leq l'-1$. 
We have seen 
$$\sharp(Z_1\bmod{2^{r+3}})\leq 2^2$$
before. 
Consequently, we get 
\begin{align*}
\sharp Z_0 &\leq 
(2^{2(n-l')}-2^{n-2})\cdot 2^{l'-1}
+3\cdot 2^{n-2}\sum_{i=3}^{l'-2}2^i
+7\cdot 2^{n-3}\cdot 2^2
+2^{n-3}\sharp(Z_0\bmod{2^{r+3}})\\
&=b(u,2)_n+2^{n-3}(\sharp(H/H_{r+3})\cap\Conj(u^{2^r})-4).
\end{align*}

\end{proof}


\section{Proof of the main theorem}

Now we prove Theorem \ref{main}. 
Recall the four types of groups 
\begin{align*}
B&=\left\{\left(\begin{matrix} *&*\\ 0&* \end{matrix}\right)\right\}
\cap \SL_2(\Z/p\Z),\\
C&=\left\{\left(\begin{matrix} *&0\\ 0&* \end{matrix}\right),
  \left(\begin{matrix} 0&*\\ *&0 \end{matrix}\right)\right\}
\cap \SL_2(\Z/p\Z),\\
D&=\left\{\left(\begin{matrix} x&y\\ \lambda y&x \end{matrix}\right),
  \left(\begin{matrix} x&y\\ -\lambda y&-x \end{matrix}\right)\right\}
\cap \SL_2(\Z/p\Z),\\
E&=\text{(an exceptional subgroup)}\cap \SL_2(\Z/p\Z).
\end{align*}

For a subgroup 
$H\subseteq G=\SL_2(\Z/p^n\Z)$, 
recall (\ref{deltaH}):
\begin{equation*}
\delta_H
=1
  -3\frac{\sharp H\cap\Conj(\sigma)}{\sharp \Conj(\sigma)}
  -4\frac{\sharp H\cap\Conj(\tau)}{\sharp \Conj(\tau)}
  -6\frac{\sharp I\backslash G/H}{[G:H]},
\end{equation*}
and (\ref{cusp1}), (\ref{cusp2}):
\begin{equation*}
\frac{\sharp I\backslash G/H}{[G:H]}=\sum_{s=0}^{n-1}
\frac{p-1}{p^{s+1}}\frac{\sharp H\cap\Conj(u^{p^s})}{\sharp \Conj(u^{p^s})}
+\frac{1}{p^n}
\leq \sum_{s=0}^{t-1}
\frac{p-1}{p^{s+1}}\frac{\sharp H\cap\Conj(u^{p^s})}{\sharp \Conj(u^{p^s})}
+\frac{1}{p^t},
\end{equation*}
where $1\leq t\leq n$.

\begin{lemma}
\label{17E}

Let $p\geq 5$ and assume a subgroup 
$H\subseteq \SL_2(\Z/p\Z)$ is contained in an exceptional subgroup $E$. 
If $p\geq 17$, then $\delta_H>0$. 

\end{lemma}

\begin{proof}

In $\SL_2(\Z/p\Z)$, we have 
$\sharp\Conj(\sigma)\geq (p-1)p$ 
and 
$\sharp\Conj(\tau)\geq (p-1)p$ 
by Lemma \ref{sl2-conj}. 
Lemma \ref{bcde-conj} shows that 
in $\SL_2(\Z/p\Z)$ we have 
$\sharp E\cap \Conj(\sigma)
\leq 30$, 
$\sharp E\cap \Conj(\tau)
\leq 20$ and 
$\sharp E\cap \Conj(u)=0$. 
Since $\sharp E\cap \Conj(u)=0$, 
we have $\frac{\sharp I\backslash G/H}{[G:H]}=\frac{1}{p}$. 
Therefore 
$\delta_H\geq 
1-3\cdot\frac{30}{(p-1)p}
-4\cdot\frac{20}{(p-1)p}
-6\cdot\frac{1}{p}
=\frac{p^2-7p-164}{(p-1)p}>0$ 
if $p\geq 17$. 

\end{proof}


\begin{proposition}
\label{19B}

Assume $p=19$. For any slim subgroup $H\subseteq \SL_2(\Z/19^2\Z)$
satisfying $H/H_1\subseteq B$, 
we have $\delta_H>0$. 

\end{proposition}

\begin{proof}

In $\SL_2(\Z/19^2\Z)$, we have 
$\sharp\Conj(\tau)=20\cdot 19^3$ 
and 
$\sharp\Conj(u)=\frac{1}{2}(p^2-1)p^2$ 
by Lemma \ref{sl2-conj}. 
Lemma \ref{bcde-conj} shows that 
in $\SL_2(\Z/19\Z)$ we have  
$\sharp B
\cap \Conj(\sigma)=0$, 
$\sharp B
\cap \Conj(\tau)=2p=38$ and 
$\sharp B
\cap \Conj(u)=\frac{1}{2}(p-1)=9$. 
By Corollary \ref{atp}, we have 
$\sharp H\cap \Conj(\tau)\leq
a(\tau,p)_2+p(38-2)=74\cdot 19$.
Hence $\frac{\sharp H\cap\Conj(\tau)}{\sharp\Conj(\tau)}\leq 
\frac{74\cdot 19}{20\cdot 19^3}
=\frac{37}{10\cdot 19^2}$. 
We have 
$\sharp H\cap\Conj(u)\leq p^2\sharp(H/H_1)\cap\Conj(u)$ 
by Lemma \ref{V-str1}. 
Thus 
$\frac{\sharp H\cap\Conj(u)}{\sharp\Conj(u)}\leq 
\frac{p^2\cdot\frac{1}{2}(p-1)}{\frac{1}{2}(p^2-1)p^2}=\frac{1}{p+1}$. 
Therefore 
$\frac{\sharp I\backslash G/H}{[G:H]}\leq 
\frac{p-1}{p}\cdot\frac{1}{p+1}+\frac{1}{p}
=\frac{2}{p+1}=\frac{1}{10}$. 
Consequently, 
$\delta_H \geq 1-3\cdot 0
-4\cdot\frac{37}{10\cdot 19^2}-6\cdot\frac{1}{10}
=\frac{1805-74-1083}{5\cdot 19^2}>0$, as required. 

\end{proof}


\begin{proposition}
\label{17B}

Assume $p=17$. For any slim subgroup $H\subseteq \SL_2(\Z/17^2\Z)$ 
satisfying $H/H_1\subseteq B$, 
we have $\delta_H>0$. 

\end{proposition}

\begin{proof}

In $\SL_2(\Z/17^2\Z)$ we have 
$\sharp\Conj(\sigma)=18\cdot 17^3$ 
and 
$\sharp\Conj(u)=\frac{1}{2}(p^2-1)p^2$ 
by Lemma \ref{sl2-conj}. 
Lemma \ref{bcde-conj} shows that 
in $\SL_2(\Z/17\Z)$ we have 
$\sharp B
\cap \Conj(\sigma)=2p=34$, 
$\sharp B
\cap \Conj(\tau)=0$ and 
$\sharp B
\cap \Conj(u)=\frac{1}{2}(p-1)=8$. 
By Corollary \ref{asp}, we have 
$\sharp H\cap \Conj(\sigma)\leq
a(\sigma,p)_2+p(34-2)=66\cdot 17$. 
Hence $\frac{\sharp H\cap\Conj(\sigma)}{\sharp\Conj(\sigma)}\leq 
\frac{66\cdot 17}{18\cdot 17^3}
=\frac{11}{3\cdot 17^2}$. 
We also have 
$\frac{\sharp H\cap\Conj(u)}{\sharp\Conj(u)}\leq \frac{1}{p+1}$ 
by the same reason as in Proposition \ref{19B}. 
Thus 
$\frac{I\backslash G/H}{[G:H]}\leq 
\frac{p-1}{p}\cdot\frac{1}{p+1}+\frac{1}{p}
=\frac{2}{p+1}=\frac{1}{9}$. 
Consequently, 
$\delta_H \geq 1-3\cdot \frac{11}{3\cdot 17^2}
-4\cdot 0-6\cdot\frac{1}{9}
=\frac{867-33-578}{3\cdot 17^2}>0$. 

\end{proof}


\begin{proposition}
\label{13BE}

Assume $p=13$. Take a slim subgroup $H\subseteq \SL_2(\Z/13^2\Z)$. 
If $H/H_1$ is contained in 
$B$ or 
$E$, then $\delta_H>0$. 

\end{proposition}

\begin{proof}

In $\SL_2(\Z/13^2\Z)$ we have 
$\sharp\Conj(\sigma)=\sharp\Conj(\tau)=14\cdot 13^3$, 
$\sharp\Conj(u)=\frac{1}{2}(p^2-1)p^2$ 
and 
$\sharp\Conj(u^p)=\frac{1}{2}(p^2-1)$
by Lemma \ref{sl2-conj}, \ref{sl2-upr}. 


Suppose $H/H_1\subseteq B$. 
Lemma \ref{bcde-conj} shows that 
in $\SL_2(\Z/13\Z)$ we have 
$\sharp B
\cap \Conj(\sigma)=2p=26$, 
$\sharp B
\cap \Conj(\tau)=2p=26$ and 
$\sharp B
\cap \Conj(u)=\frac{1}{2}(p-1)=6$. 
By Corollary \ref{asp}, we have 
$\sharp H\cap \Conj(\sigma)\leq
a(\sigma,p)_2+p(26-2)=50\cdot 13$. 
Hence $\frac{\sharp H\cap\Conj(\sigma)}{\sharp\Conj(\sigma)}\leq 
\frac{50\cdot 13}{14\cdot 13^3}
=\frac{25}{7\cdot 13^2}$. 
The same calculation shows 
$\frac{\sharp H\cap\Conj(\tau)}{\sharp\Conj(\tau)}\leq 
\frac{25}{7\cdot 13^2}$. 
We have $B\cap
  \Conj(u)
  =\left\{\left(\begin{matrix} 1&s^2\\ 0&1 \end{matrix}\right)
  \in \SL_2(\Z/p\Z)|
  (p,s)=1\right\}$ 
by Lemma \ref{bcde-conj}, 
and 
$V_{u'}=V_u =\left\{1+p
  \left(\begin{matrix} a&b\\ 0&-a \end{matrix}\right)\right\}$ 
for any $u'\in \left\{\left(\begin{matrix} 1&s^2\\ 0&1 \end{matrix}\right)|
  (p,s)=1\right\}$ 
by Lemma \ref{V-str1} and Corollary \ref{u-p}.

If $H$ contains $V_u$, 
then the $(2,1)$ entry of an element of $H\cap \Conj(u^p)$ must be zero 
since $H$ is slim, 
thus 
$H\cap \Conj(u^p)\subseteq 
\left\{\left(\begin{matrix} 1&ps^2\\ 0&1 \end{matrix}\right)|
  (p,s)=1\right\}$. 
Hence 
$\frac{\sharp H\cap\Conj(u^p)}{\sharp\Conj(u^p)}\leq 
\frac{\frac{1}{2}(p-1)}{\frac{1}{2}(p^2-1)}=\frac{1}{p+1}$. 
We also have 
$\frac{\sharp H\cap\Conj(u)}{\sharp\Conj(u)}\leq\frac{1}{p+1}$ 
as in Proposition \ref{19B}. 
Thus 
$\frac{\sharp I\backslash G/H}{[G:H]}\leq 
\frac{p-1}{p}\cdot\frac{1}{p+1}+\frac{p-1}{p^2}\cdot\frac{1}{p+1}
+\frac{1}{p^2}=\frac{1}{p}=\frac{1}{13}$. 
Therefore 
$\delta_H\geq 1-3\cdot\frac{25}{7\cdot 13^2}
-4\cdot\frac{25}{7\cdot 13^2}-6\cdot\frac{1}{13}
=\frac{1183-75-100-546}{7\cdot 13^2}>0$.

If $H$ does not contain $V_u$, 
for each element 
$x \in (H/H_1) \cap \Conj(u)\subseteq B\cap\Conj(u)$ 
(there are at most $\frac{1}{2}(p-1)$ such elements 
by Lemma \ref{bcde-conj}) 
we have 
$\sharp (f^{H,u}_{2,1})^{-1}(x)\leq p=13$ by Lemma \ref{leq-p}. 
Hence $\frac{\sharp H\cap\Conj(u)}{\sharp\Conj(u)}\leq 
\frac{p\cdot\frac{1}{2}(p-1)}{\frac{1}{2}(p^2-1)p^2}=\frac{1}{p(p+1)}$ and 
$\frac{\sharp I\backslash G/H}{[G:H]}\leq 
\frac{p-1}{p}\cdot\frac{1}{p(p+1)}+\frac{1}{p}
=\frac{p^2+2p-1}{(p+1)p^2}
=\frac{97}{7\cdot 13^2}$. 
Therefore 
$\delta_H\geq 1-3\cdot\frac{25}{7\cdot 13^2}
-4\cdot\frac{25}{7\cdot 13^2}
-6\cdot\frac{97}{7\cdot 13^2}
=\frac{1183-75-100-582}{7\cdot 13^2}>0$.


Next suppose $H/H_1\subseteq E$. 
Lemma \ref{bcde-conj} shows that 
in $\SL_2(\Z/13\Z)$ we have 
$\sharp E\cap \Conj(\sigma)
\leq 18$, 
$\sharp E\cap \Conj(\tau)
\leq 8$ and 
$\sharp E\cap \Conj(u)=0$. 
By Corollary \ref{asp}, we have 
$\sharp H\cap \Conj(\sigma)\leq
a(\sigma,p)_2+p(18-2)=42\cdot 13$. 
Hence $\frac{\sharp H\cap\Conj(\sigma)}{\sharp\Conj(\sigma)}\leq 
\frac{42\cdot 13}{14\cdot 13^3}
=\frac{3}{13^2}$. 
By Corollary \ref{atp}, we have 
$\sharp H\cap \Conj(\tau)\leq
a(\tau,p)_2+p(8-2)=32\cdot 13$. 
Hence $\frac{\sharp H\cap\Conj(\tau)}{\sharp\Conj(\tau)}\leq 
\frac{32\cdot 13}{14\cdot 13^3}
=\frac{16}{7\cdot 13^2}$. 
Since $\sharp(H/H_1)\cap \Conj(u)=0$, we have 
$\frac{\sharp I\backslash G/H}{[G:H]}\leq \frac{1}{p}
=\frac{1}{13}$ 
(take $t=1$ in (4.3)). 
Therefore 
$\delta_H\geq 1-3\cdot\frac{3}{13^2}
-4\cdot\frac{16}{7\cdot 13^2}
-6\cdot\frac{1}{13}
=\frac{1183-63-64-546}{7\cdot 13^2}>0$. 

\end{proof}


\begin{proposition}
\label{11BDE}

Assume $p=11$. Take a slim subgroup $H\subseteq \SL_2(\Z/11^2\Z)$. 
If $H/H_1$ is contained in 
$B$, 
$D$ 
or $E$, then $\delta_H>0$. 

\end{proposition}

\begin{proof}

In $\SL_2(\Z/11^2\Z)$ we have 
$\sharp\Conj(\sigma)=\sharp\Conj(\tau)=10\cdot 11^3$, 
$\sharp\Conj(u)=\frac{1}{2}(p^2-1)p^2$ 
and 
$\sharp\Conj(u^p)=\frac{1}{2}(p^2-1)$ 
by Lemma \ref{sl2-conj}, \ref{sl2-upr}. 


Suppose $H/H_1\subseteq B$. 
Lemma \ref{bcde-conj} shows that 
in $\SL_2(\Z/11\Z)$ we have 
$\sharp B
\cap \Conj(\sigma)=0$, 
$\sharp B
\cap \Conj(\tau)=0$ and 
$\sharp B
\cap \Conj(u)=\frac{1}{2}(p-1)=5$.

If $H$ contains $V_u$, 
then the same calculation as in Proposition \ref{13BE} shows 
$\frac{\sharp I\backslash G/H}{[G:H]}\leq 
\frac{1}{p}=\frac{1}{11}$. 
Therefore 
$\delta_H\geq 1-3\cdot 0
-4\cdot 0-6\cdot\frac{1}{11}
=\frac{11-6}{11}>0$.

If $H$ does not contain $V_u$, 
as in Proposition \ref{13BE}, 
$\frac{\sharp I\backslash G/H}{[G:H]}\leq 
\frac{p^2+2p-1}{(p+1)p^2}
=\frac{71}{2\cdot 3\cdot 11^2}$. 
Therefore 
$\delta_H\geq 1-3\cdot 0
-4\cdot 0
-6\cdot\frac{71}{2\cdot 3\cdot 11^2}
=\frac{121-71}{11^2}>0$.


Next suppose $H/H_1\subseteq E$. 
Lemma \ref{bcde-conj} shows that 
in $\SL_2(\Z/11\Z)$ we have 
$\sharp E\cap \Conj(\sigma)
\leq 30$, 
$\sharp E\cap \Conj(\tau)
\leq 20$ and 
$\sharp E\cap \Conj(u)=0$. 
By Corollary \ref{asp}, we have 
$\sharp H\cap \Conj(\sigma)\leq
a(\sigma,p)_2+p(30-2)=50\cdot 11$. 
Hence $\frac{\sharp H\cap\Conj(\sigma)}{\sharp\Conj(\sigma)}\leq 
\frac{50\cdot 11}{10\cdot 11^3}
=\frac{5}{11^2}$. 
By Corollary \ref{atp}, we have 
$\sharp H\cap \Conj(\tau)\leq
a(\tau,p)_2+p(20-2)=40\cdot 11$. 
Hence $\frac{\sharp H\cap\Conj(\tau)}{\sharp\Conj(\tau)}\leq 
\frac{40\cdot 11}{10\cdot 11^3}
=\frac{4}{11^2}$. 
Since $\sharp(H/H_1)\cap \Conj(u)=0$, we have 
$\frac{\sharp I\backslash G/H}{[G:H]}\leq \frac{1}{p}
=\frac{1}{11}$. 
Therefore 
$\delta_H\geq 1-3\cdot\frac{5}{11^2}
-4\cdot\frac{4}{11^2}
-6\cdot\frac{1}{11}
=\frac{121-15-16-66}{11^2}>0$. 


Finally, suppose 
$H/H_1\subseteq D$. 
Lemma \ref{bcde-conj} shows that 
in $\SL_2(\Z/11\Z)$ we have 
$\sharp 
D
\cap \Conj(\sigma)=p+3=14$, 
$\sharp 
D
\cap \Conj(\tau)=2$ and 
$\sharp 
D
\cap \Conj(u)=0$. 
Since 
$\sharp D\cap \Conj(\sigma)<30$, 
$\sharp D\cap \Conj(\tau)<20$ and 
$\sharp D\cap \Conj(u)=0$, 
the calculation for $E$ shows $\delta_H>0$ 
in this case. 

\end{proof}


\begin{proposition}
\label{7BCDE}

Assume $p=7$. Take a slim subgroup $H\subseteq \SL_2(\Z/7^3\Z)$. 
If $H/H_1$ is contained in 
$B$, 
$C$, 
$D$ 
or $E$ , then $\delta_H>0$. 

\end{proposition}

\begin{proof}

In $\SL_2(\Z/7^3\Z)$ we have 
$\sharp\Conj(\sigma)=6\cdot 7^5$, 
$\sharp\Conj(\tau)=8\cdot 7^5$, 
$\sharp\Conj(u)=\frac{1}{2}(p^2-1)p^4$, 
$\sharp\Conj(u^p)=\frac{1}{2}(p^2-1)p^2$ 
and 
$\sharp\Conj(u^{p^2})=\frac{1}{2}(p^2-1)$
by Lemma \ref{sl2-conj}, \ref{sl2-upr}. 


Suppose $H/H_1\subseteq B$. 
Lemma \ref{bcde-conj} shows that 
in $\SL_2(\Z/7\Z)$ we have 
$\sharp B
\cap \Conj(\sigma)=0$, 
$\sharp B
\cap \Conj(\tau)=2p=14$ and 
$\sharp B
\cap \Conj(u)=\frac{1}{2}(p-1)=3$. 
By Corollary \ref{atp}, we have 
$\sharp H\cap \Conj(\tau)
\leq a(\tau,p)_3+p^2(14-2)=110\cdot 7^2$. 
Hence $\frac{\sharp H\cap\Conj(\tau)}{\sharp\Conj(\tau)}\leq 
\frac{110\cdot 7^2}{8\cdot 7^5}
=\frac{55}{2^2\cdot 7^3}$. 

If $H/H_2$ contains 
$V_u^{2,1}=\left\{1+p
  \left(\begin{matrix} a&b\\ 0&-a \end{matrix}\right)\right\}$, 
then the $(2,1)$ entry of an element of $(H/H_2)\cap \Conj(u^p)$ must be zero 
as we have seen in Proposition \ref{13BE},
thus 
$(H/H_2)\cap \Conj(u^p)\subseteq 
\left\{\left(\begin{matrix} 1&ps^2\\ 0&1 \end{matrix}\right)|
  (p,s)=1\right\}$. 
We have 
$\sharp H\cap\Conj(u^p)\leq p^2\sharp(H/H_2)\cap\Conj(u^p)$ 
by Lemma \ref{V-str1}. 
Hence 
$\frac{\sharp H\cap\Conj(u^p)}{\sharp\Conj(u^p)}\leq 
\frac{p^2\cdot\frac{1}{2}(p-1)}{\frac{1}{2}(p^2-1)p^2}=\frac{1}{p+1}$. 
As $H/H_2$ contains  
$\left\{1+p
  \left(\begin{matrix} a&b\\ 0&-a \end{matrix}\right)\right\}$, 
the subgroup $H$ contains 
$\left\{1+p^2
  \left(\begin{matrix} a&b\\ 0&-a \end{matrix}\right)\right\}$. 
Thus the $(2,1)$ entry of an element of $H\cap \Conj(u^{p^2})$ 
must be zero since $H$ is slim, 
and 
$H\cap \Conj(u^{p^2})\subseteq 
\left\{\left(\begin{matrix} 1&p^2s^2\\ 0&1 \end{matrix}\right)|
  (p,s)=1\right\}$. 
Hence 
$\frac{\sharp H\cap\Conj(u^{p^2})}{\sharp\Conj(u^{p^2})}\leq 
\frac{\frac{1}{2}(p-1)}{\frac{1}{2}(p^2-1)}=\frac{1}{p+1}$. 
We also have 
$\frac{\sharp H\cap\Conj(u)}{\sharp\Conj(u)}\leq \frac{1}{p+1}$ 
as in Proposition \ref{19B}. 
Therefore 
$\frac{\sharp I\backslash G/H}{[G:H]}\leq 
\frac{p-1}{p}\cdot \frac{1}{p+1}+\frac{p-1}{p^2}\cdot \frac{1}{p+1}
+\frac{p-1}{p^3}\cdot \frac{1}{p+1}+\frac{1}{p^3}
=\frac{p^2+1}{(p+1)p^2}=\frac{25}{2^2\cdot 7^2}$. 
Consequently, 
$\delta_H\geq 1-3\cdot 0-4\cdot \frac{55}{2^2\cdot 7^3}
-6\cdot \frac{25}{2^2\cdot 7^2}
=\frac{686-110-525}{2\cdot 7^3}>0$. 

If $H/H_2$ does not contain $V_u^{2,1}$, 
for each element 
$x \in (H/H_1) \cap \Conj(u)\subseteq B\cap\Conj(u)$ 
(there are at most 
$\frac{1}{2}(p-1)$ such elements by Lemma \ref{bcde-conj}) we have 
$\sharp (f^{H,u}_{3,1})^{-1}(x)\leq
 p^2\sharp (f^{H,u}_{2,1})^{-1}(x)\leq p^3$ 
by Lemma \ref{V-str1}, \ref{leq-p}. 
Thus 
$\frac{\sharp H\cap\Conj(u)}{\sharp\Conj(u)}\leq 
\frac{p^3\cdot\frac{1}{2}(p-1)}{\frac{1}{2}(p^2-1)p^4}
=\frac{1}{(p+1)p}=\frac{1}{2^3\cdot 7}$. 
By Corollary \ref{aup}, we have 
$\sharp H\cap \Conj(u^p)
\leq a(u,p)_2+p(\frac{1}{2}(p^2-1)-\frac{1}{2}(p-1))
=(p-1)p^2$. 
Thus 
$\frac{\sharp H\cap\Conj(u^p)}{\sharp\Conj(u^p)}\leq 
\frac{(p-1)p^2}
{\frac{1}{2}(p^2-1)p^2}
=\frac{2}{p+1}
=\frac{1}{4}$. 
Therefore 
$\frac{\sharp I\backslash G/H}{[G:H]}\leq 
\frac{p-1}{p}\cdot \frac{1}{(p+1)p}+\frac{p-1}{p^2}\cdot \frac{2}{p+1}
+\frac{1}{p^2}=\frac{4p-2}{(p+1)p^2}=\frac{13}{2^2\cdot 7^2}$. 
Consequently, 
$\delta_H\geq 1-3\cdot 0-4\cdot \frac{55}{2^2\cdot 7^3}
-6\cdot \frac{13}{2^2\cdot 7^2}
=\frac{686-110-273}{2\cdot 7^3}>0$.


Next suppose $H/H_1\subseteq E$. 
Lemma \ref{bcde-conj} shows that 
in $\SL_2(\Z/7\Z)$ we have 
$\sharp E\cap \Conj(\sigma)
\leq 18$, 
$\sharp E\cap \Conj(\tau)
\leq 8$ and 
$\sharp E\cap \Conj(u)=0$. 
By Corollary \ref{asp}, we have 
$\sharp H\cap \Conj(\sigma)
\leq a(\sigma,p)_3+p^2(18-2)=114\cdot 7^2$. 
Hence $\frac{\sharp H\cap\Conj(\sigma)}{\sharp\Conj(\sigma)}\leq 
\frac{114\cdot 7^2}{6\cdot 7^5}
=\frac{19}{7^3}$. 
By Corollary \ref{atp}, we have 
$\sharp H\cap \Conj(\tau)
\leq a(\tau,p)_3+p^2(8-2)=104\cdot 7^2$. 
Hence $\frac{\sharp H\cap\Conj(\tau)}{\sharp\Conj(\tau)}\leq 
\frac{104\cdot 7^2}{8\cdot 7^5}
=\frac{13}{7^3}$. 
The same calculation as above 
(the case that $H/H_1\subseteq 
B$ and 
$H/H_2$ does not contain $V_u^{2,1}$) shows the inequality 
$\frac{\sharp H\cap\Conj(u^p)}{\sharp\Conj(u^p)}\leq 
\frac{2}{p+1}
=\frac{1}{4}$. 
Therefore 
$\frac{\sharp I\backslash G/H}{[G:H]}\leq 
\frac{p-1}{p}\cdot 0+\frac{p-1}{p^2}\cdot \frac{2}{p+1}
+\frac{1}{p^2}=\frac{3p-1}{(p+1)p^2}=\frac{5}{2\cdot 7^2}$. 
Consequently, 
$\delta_H\geq 1-3\cdot \frac{19}{7^3}
-4\cdot \frac{13}{7^3}
-6\cdot \frac{5}{2\cdot 7^2}
=\frac{343-57-52-105}{7^3}>0$.


Suppose 
$H/H_1\subseteq C$. 
Lemma \ref{bcde-conj} shows that in $\SL_2(\Z/7\Z)$ 
we have 
$\sharp
C
\cap \Conj(\sigma)=p-1=6<18$, 
$\sharp
C
\cap \Conj(\tau)=2<8$, 
$\sharp
C
\cap \Conj(u)=0$. 
The calculation for $E$ shows $\delta_H>0$ in this case.


Finally, suppose 
$H/H_1\subseteq D$. 
Lemma \ref{bcde-conj} shows that in $\SL_2(\Z/7\Z)$ 
we have 
$\sharp 
D
\cap \Conj(\sigma)=p+3=10<18$, 
$\sharp 
D
\cap \Conj(\tau)=0<8$ and 
$\sharp 
D
\cap \Conj(u)=0$. 
The calculation for $E$ 
also shows $\delta_H>0$ in this case. 

\end{proof}



\begin{proposition}
\label{5C}

Assume $p=5$. Take a slim subgroup $H\subseteq \SL_2(\Z/5^3\Z)$. 
If $H/H_1$ is contained in 
$C$, 
then $\delta_H>0$. 

\end{proposition}

\begin{proof}

In $\SL_2(\Z/5^3\Z)$ we have 
$\sharp\Conj(\sigma)=6\cdot 5^5$, 
and 
$\sharp\Conj(u^p)=\frac{1}{2}(5^2-1)5^2=12\cdot 5^2$
by Lemma \ref{sl2-conj}, \ref{sl2-upr}. 
Lemma \ref{bcde-conj} shows that 
in $\SL_2(\Z/5\Z)$ we have  
$\sharp C
\cap \Conj(\sigma)=p+1=6$, 
$\sharp C
\cap \Conj(\tau)=0$ and 
$\sharp C
\cap \Conj(u)=0$. 
By Corollary \ref{asp}, we have 
$\sharp H\cap \Conj(\sigma)
\leq a(\sigma,p)_3+p^3(6-2)=54\cdot 5^2$. 
Hence $\frac{\sharp H\cap\Conj(\sigma)}{\sharp\Conj(\sigma)}\leq 
\frac{54\cdot 5^2}{6\cdot 5^5}
=\frac{9}{5^3}$. 
In $\SL_2(\Z/5^2\Z)$ we have 
$\sharp\Conj(u^5)=12$
by Lemma \ref{sl2-upr}. 
By Corollary \ref{aup}, we have 
$\sharp H\cap \Conj(u^p)
\leq a(u,p)_2+p(12-\frac{1}{2}(5-1))
=4\cdot 5^2$. 
Hence 
$\frac{\sharp H\cap\Conj(u^p)}{\sharp\Conj(u^p)}\leq 
\frac{4\cdot 5^2}{12\cdot 5^2}
=\frac{1}{3}$. 
Therefore 
$\frac{\sharp I\backslash G/H}{[G:H]}\leq 
\frac{5-1}{5}\cdot 0+\frac{5-1}{5^2}\cdot \frac{1}{3}
+\frac{1}{5^2}=\frac{7}{3\cdot 5^2}$. 
Consequently, 
$\delta_H\geq 1-3\cdot \frac{9}{5^3}
-4\cdot 0
-6\cdot \frac{7}{3\cdot 5^2}
=\frac{125-27-70}{5^3}>0$. 

\end{proof}

\begin{proposition}
\label{5BDE}

Assume $p=5$. Take a slim subgroup $H\subseteq \SL_2(\Z/5^4\Z)$. 
If $H/H_1$ is contained in 
$B$, $D$ 
or $E$, then $\delta_H>0$. 

\end{proposition}

\begin{proof}

In $\SL_2(\Z/5^4\Z)$ we have 
$\sharp\Conj(\sigma)=6\cdot 5^7$, 
$\sharp\Conj(\tau)=4\cdot 5^7$, 
$\sharp\Conj(u)=12\cdot 5^6$, 
$\sharp\Conj(u^p)=12\cdot 5^4$
and 
$\sharp\Conj(u^{p^2})=12\cdot 5^2$
by Lemma \ref{sl2-conj}, \ref{sl2-upr}. 


Suppose $H/H_1\subseteq B$. 
Lemma \ref{bcde-conj} shows that 
in $\SL_2(\Z/5\Z)$ we have  
$\sharp B
\cap \Conj(\sigma)=2p=10$, 
$\sharp B
\cap \Conj(\tau)=0$ and 
$\sharp B
\cap \Conj(u)=\frac{1}{2}(p-1)=2$. 
By Corollary \ref{asp}, we have 
$\sharp H\cap \Conj(\sigma)
\leq a(\sigma,p)_4+p^3(10-2)=66\cdot 5^3$. 
Hence $\frac{\sharp H\cap\Conj(\sigma)}{\sharp\Conj(\sigma)}\leq 
\frac{66\cdot 5^3}{6\cdot 5^7}
=\frac{11}{5^4}$. 
By Corollary \ref{aup}, we have 
$\sharp H\cap \Conj(u)
\leq a(u,p)_4=18\cdot 5^4$. 
Hence 
$\frac{\sharp H\cap\Conj(u)}{\sharp\Conj(u)}\leq 
\frac{18\cdot 5^4}{12\cdot 5^6}
=\frac{3}{2\cdot 5^2}$. 
In $\SL_2(\Z/5^2\Z)$ we have 
$\sharp\Conj(u^p)=12$
by Lemma \ref{sl2-upr}. 
Again by Corollary \ref{aup}, we have 
$\sharp H\cap \Conj(u^p)
\leq a(u,p)_3+p^2(12-\frac{1}{2}(5-1))
=12\cdot 5^3$. 
Hence 
$\frac{\sharp H\cap\Conj(u^p)}{\sharp\Conj(u^p)}\leq 
\frac{12\cdot 5^3}{12\cdot 5^4}
=\frac{1}{5}$. 
Similarly 
$\sharp H\cap \Conj(u^{p^2})
\leq a(u,p)_2+p(12-\frac{1}{2}(5-1))
=4\cdot 5^2$. 
Hence 
$\frac{\sharp H\cap\Conj(u^{p^2})}{\sharp\Conj(u^{p^2})}\leq 
\frac{4\cdot 5^2}{12\cdot 5^2}
=\frac{1}{3}$. 
Therefore 
$\frac{\sharp I\backslash G/H}{[G:H]}\leq 
\frac{5-1}{5}\cdot \frac{3}{2\cdot 5^2}+\frac{5-1}{5^2}\cdot \frac{1}{5}
+\frac{5-1}{5^3}\cdot \frac{1}{3}
+\frac{1}{5^3}=\frac{37}{3\cdot 5^3}$. 
Consequently, 
$\delta_H\geq 1-3\cdot \frac{11}{5^4}
-4\cdot 0
-6\cdot \frac{37}{3\cdot 5^3}
=\frac{625-33-370}{5^4}>0$.


Next suppose $H/H_1\subseteq E$. 
Lemma \ref{bcde-conj} shows that 
in $\SL_2(\Z/5\Z)$ we have 
$\sharp E\cap \Conj(\sigma)
\leq 18$, 
$\sharp E\cap \Conj(\tau)
\leq 8$ and 
$\sharp E\cap \Conj(u)=0$. 
By Corollary \ref{asp}, we have 
$\sharp H\cap \Conj(\sigma)
\leq a(\sigma,p)_4+p^3(18-2)=74\cdot 5^3$. 
Hence $\frac{\sharp H\cap\Conj(\sigma)}{\sharp\Conj(\sigma)}\leq 
\frac{74\cdot 5^3}{6\cdot 5^7}
=\frac{37}{3\cdot 5^4}$. 
By Corollary \ref{atp}, we have 
$\sharp H\cap \Conj(\tau)
\leq a(\tau,p)_4+p^3(8-2)=64\cdot 5^3$. 
Hence $\frac{\sharp H\cap\Conj(\tau)}{\sharp\Conj(\tau)}\leq 
\frac{64\cdot 5^3}{4\cdot 5^7}
=\frac{16}{5^4}$. 
The calculation in the Borel case shows 
$\frac{\sharp I\backslash G/H}{[G:H]}\leq 
\frac{5-1}{5}\cdot 0+\frac{5-1}{5^2}\cdot \frac{1}{5}
+\frac{5-1}{5^3}\cdot \frac{1}{3}
+\frac{1}{5^3}=\frac{19}{3\cdot 5^3}$. 
Consequently, 
$\delta_H\geq 1-3\cdot \frac{37}{3\cdot 5^4}
-4\cdot \frac{16}{5^4}
-6\cdot \frac{19}{3\cdot 5^3}
=\frac{625-37-64-190}{5^4}>0$.


Finally, suppose 
$H/H_1\subseteq D$. 
Lemma \ref{bcde-conj} shows that in $\SL_2(\Z/5\Z)$ 
we have 
$\sharp 
D
\cap \Conj(\sigma)=p+1=6<18$, 
$\sharp 
D
\cap \Conj(\tau)=2<8$ and 
$\sharp 
D
\cap \Conj(u)=0$. 
The calculation for $E$ 
shows $\delta_H>0$ in this case. 

\end{proof}


\begin{proposition}
\label{3BCDSL}

Assume $p=3$. Take a slim subgroup $H\subseteq \SL_2(\Z/3^6\Z)$. 
If $H/H_1$ is contained in 
$B$, $C$, $D$ 
or $H/H_1=\SL_2(\Z/3\Z)$, 
then $\delta_H>0$. 

\end{proposition}

\begin{proof}

In $\SL_2(\Z/3^6\Z)$ we have 
$\sharp\Conj(\sigma)=2\cdot 3^{11}$, 
$\sharp\Conj(\tau)=4\cdot 3^{10}$ 
and 
$\sharp\Conj(u^{3^i})=4\cdot 3^{10-2i}$ 
for $0\leq i\leq 4$ 
by Lemma \ref{sl2-conj}, \ref{sl2-upr}. 


Suppose $H/H_1\subseteq B$. 
Lemma \ref{bcde-conj} shows that 
in $\SL_2(\Z/3\Z)$ we have 
$\sharp B
\cap \Conj(\sigma)=0$, 
$\sharp B
\cap \Conj(\tau)=1$ and 
$\sharp B
\cap \Conj(u)=1$. 
By Corollary \ref{at3}, we have 
$\sharp H\cap\Conj(\tau)\leq 
a(\tau,3)_6=13\cdot 3^6$. 
Hence $\frac{\sharp H\cap\Conj(\tau)}{\sharp\Conj(\tau)}\leq 
\frac{13\cdot 3^6}{4\cdot 3^{10}}
=\frac{13}{2^2\cdot 3^4}$. 
By Corollary \ref{aup}, we have 
$\sharp H\cap\Conj(u)\leq
a(u,3)_6=17\cdot 3^6$. 
Hence $\frac{\sharp H\cap\Conj(u)}{\sharp\Conj(u)}\leq 
\frac{17\cdot 3^6}{4\cdot 3^{10}}
=\frac{17}{2^2\cdot 3^4}$. 
In $\SL_2(\Z/3^{r+1}\Z)$ we have 
$\sharp\Conj(u^{3^r})=\frac{1}{2}(3^2-1)=4$ 
by Lemma \ref{sl2-upr}. 
Again 
by Corollary \ref{aup}, we have 
$\sharp H\cap\Conj(u^{3^i})\leq
a(u,3)_{6-i}+3^{5-i}(4-\frac{1}{2}(3-1))
=a(u,3)_{6-i}+3^{6-i}$ 
for $1\leq i\leq 4$. 
Hence 
$$\frac{\sharp H\cap\Conj(u^{3^i})}{\sharp\Conj(u^{3^i})}\leq 
\frac{a(u,3)_{6-i}+3^{6-i}}{4\cdot 3^{10-2i}}$$ 
for $1\leq i\leq 4$. 
Therefore 
\begin{align*}
\frac{\sharp I\backslash G/H}{[G:H]}\leq 
&\frac{3-1}{3}\cdot \frac{17}{2^2\cdot 3^4}
+\sum_{i=1}^4 \frac{3-1}{3^{i+1}}\cdot 
\frac{a(u,3)_{6-i}+3^{6-i}}{4\cdot 3^{10-2i}}
+\frac{1}{3^5}\\
=&\frac{17+12+6+4+2+2}{2\cdot 3^5}
=\frac{43}{2\cdot 3^5}.
\end{align*}
Consequently, 
$\delta_H\geq 1-3\cdot 0
-4\cdot \frac{13}{2^2\cdot 3^4}
-6\cdot \frac{43}{2\cdot 3^5}
=\frac{81-13-43}{3^4}>0$.


Next suppose $H/H_1=\SL_2(\Z/3\Z)$. 
Since $H/H_2\subsetneq \SL_2(\Z/3^2\Z)$ by Lemma 2.1, 
we have $H\cap \Conj(u)=\emptyset$ by Lemma \ref{sl2-u}. 
Lemma \ref{sl2-conj} shows that in $\SL_2(\Z/3\Z)$ we have 
$\sharp \Conj(\sigma)=6$ and 
$\sharp \Conj(\tau)=4$. 
By Corollary \ref{asp}, we have 
$\sharp H\cap\Conj(\sigma)\leq
a(\sigma,3)_6+3^5(6-2)=14\cdot 3^6$. 
Hence $\frac{\sharp H\cap\Conj(\sigma)}{\sharp\Conj(\sigma)}\leq 
\frac{14\cdot 3^6}{2\cdot 3^{11}}
=\frac{7}{3^5}$. 
By Corollary \ref{at3}, we have 
$\sharp H\cap\Conj(\tau)\leq
a(\tau,3)_6+3^5(4-1)=14\cdot 3^6$. 
Hence 
$\frac{\sharp H\cap\Conj(\tau)}{\sharp\Conj(\tau)}\leq 
\frac{14\cdot 3^6}{4\cdot 3^{10}}
=\frac{7}{2\cdot 3^4}$. 
The same calculation as in the Borel case shows 
$$\frac{\sharp I\backslash G/H}{[G:H]}\leq 
\frac{3-1}{3}\cdot 0
+\sum_{i=1}^4 \frac{3-1}{3^{i+1}}\cdot 
\frac{a(u,3)_{6-i}+3^{6-i}}{4\cdot 3^{10-2i}}
+\frac{1}{3^5}
=\frac{13}{3^5}.$$ 
Consequently, 
$\delta_H\geq 1-3\cdot \frac{7}{3^5}
-4\cdot \frac{7}{2\cdot 3^4}
-6\cdot \frac{13}{3^5}
=\frac{81-7-14-26}{3^4}>0$. 
%

%


Suppose 
$H/H_1\subseteq C$. 
Lemma \ref{bcde-conj} shows that in $\SL_2(\Z/3\Z)$
we have 
$\sharp
C
\cap \Conj(\sigma)=2<6$, 
$\sharp
C
\cap \Conj(\tau)=0<4$, 
$\sharp
C
\cap \Conj(u)=0$. 
The calculation for the case $H/H_1=\SL_2(\Z/3\Z)$ shows $\delta_H>0$ 
in this case.


Finally, suppose 
$H/H_1\subseteq D$. 
Lemma \ref{bcde-conj} shows that in $\SL_2(\Z/3\Z)$
we have 
$\sharp 
D
\cap \Conj(\sigma)=6$, 
$\sharp 
D
\cap \Conj(\tau)=0<4$ and 
$\sharp 
D
\cap \Conj(u)=0$. 
The calculation for the case $H/H_1=\SL_2(\Z/3\Z)$ also shows $\delta_H>0$ 
in this case. 

\end{proof}


As is well-known, we have $\SL_2(\Z/2\Z)\cong S_3$ as a group. 
All non-trivial subgroups of $\SL_2(\Z/2\Z)$ are 
the following:
\begin{align*}
B&=\left\{
\left(\begin{matrix} 1&0\\ 0&1 \end{matrix}\right),\,
\left(\begin{matrix} 1&1\\ 0&1 \end{matrix}\right)\right\},\,
\left\{
\left(\begin{matrix} 1&0\\ 0&1 \end{matrix}\right),\,
\left(\begin{matrix} 1&0\\ 1&1 \end{matrix}\right)\right\},\,
\left\{
\left(\begin{matrix} 1&0\\ 0&1 \end{matrix}\right),\,
\left(\begin{matrix} 0&1\\ 1&0 \end{matrix}\right)\right\},\\
F&=\left\{
\left(\begin{matrix} 1&0\\ 0&1 \end{matrix}\right),\,
\left(\begin{matrix} 1&1\\ 1&0 \end{matrix}\right),\,
\left(\begin{matrix} 0&1\\ 1&1 \end{matrix}\right)\right\}.
\end{align*}
The first three groups are conjugate.

\begin{proposition}
\label{2B}

Assume $p=2$. Take a slim subgroup $H\subseteq \SL_2(\Z/2^{11}\Z)$. 
If $H/H_1$ is contained in 
$B$, 
then $\delta_H>0$. 

\end{proposition}

\begin{proof}

In $\SL_2(\Z/2^{11}\Z)$ we have 
$\sharp\Conj(\sigma)=3\cdot 2^{19}$ 
and 
$\sharp\Conj(u^{2^i})=3\cdot 2^{18-2i}$ 
for $0\leq i\leq 7$ 
by Lemma \ref{sl2-conj}, \ref{sl2-upr}. 
%
%
By Lemma \ref{table}, we see that 
in $\SL_2(\Z/4\Z)$ we have 
$\sharp f_{2,1}^{-1}(B)
\cap \Conj(\sigma)=2$, 
$\sharp f_{2,1}^{-1}(B)
\cap \Conj(\tau)=0$, 
$\sharp f_{2,1}^{-1}(B)
\cap \Conj(u)=2$ and 
$\sharp f_{2,1}^{-1}(B)
\cap \Conj(u^2)=3$. 
By Proposition \ref{as2}, we have 
$\sharp H\cap\Conj(\sigma)\leq 
a(\sigma,2)_{11}=11\cdot 2^{12}$.
Hence 
$\frac{\sharp H\cap\Conj(\sigma)}{\sharp\Conj(\sigma)}\leq 
\frac{11\cdot 2^{12}}{3\cdot 2^{19}}
=\frac{11}{3\cdot 2^7}$. 
We have 
$\sharp(H/H_3)\cap\Conj(u)\leq \sharp f_{3,1}^{-1}(B)\cap\Conj(u)=4$ 
by Remark \ref{inv=2}. 
By Corollary \ref{au2}, we have 
$\sharp H\cap\Conj(u)\leq 
a(u,2)_{11}+2^{10}(4-2)
=23\cdot 2^{11}$.
Hence 
$\frac{\sharp H\cap\Conj(u)}{\sharp\Conj(u)}\leq 
\frac{23\cdot 2^{11}}{3\cdot 2^{18}}
=\frac{23}{3\cdot 2^7}$. 
In $\SL_2(\Z/2^{r+3}\Z)$ we have 
$\sharp\Conj(u^{2^r})=3\cdot 2^2=12$ 
by Lemma \ref{sl2-upr}. 
Again by Corollary \ref{au2}, we have 
$\sharp H\cap\Conj(u^2)\leq 
a(u,2)_{10}+2^9(12-2)
=19\cdot 2^{10}$.
Hence 
$\frac{\sharp H\cap\Conj(u^2)}{\sharp\Conj(u^2)}\leq 
\frac{19\cdot 2^{10}}{3\cdot 2^{16}}
=\frac{19}{3\cdot 2^6}$. 
By Proposition \ref{bu2}, we have 
$\sharp H\cap\Conj(u^{2^i})\leq 
b(u,2)_{11-i}+2^{8-i}(12-4)
=b(u,2)_{11-i}+2^{11-i}$ for $2\leq i\leq 7$. 
Hence 
$$\frac{\sharp H\cap\Conj(u^{2^i})}{\sharp\Conj(u^{2^i})}\leq 
\frac{b(u,2)_{11-i}+2^{11-i}}{3\cdot 2^{18-2i}}$$
for $2\leq i\leq 7$. 
Therefore 
\begin{align*}
\frac{\sharp I\backslash G/H}{[G:H]} &\leq 
\frac{2-1}{2}\cdot \frac{23}{3\cdot 2^7}
+\frac{2-1}{2^2}\cdot \frac{19}{3\cdot 2^6}
+\sum_{i=2}^7 \frac{2-1}{2^{i+1}}\cdot 
\frac{b(u,2)_{11-i}+2^{11-i}}{3\cdot 2^{18-2i}}
+\frac{1}{2^8}\\
&=\frac{23+19+15+11+7+5+3+2+3}{3\cdot 2^8}
=\frac{11}{3\cdot 2^5}.
\end{align*}
Consequently, 
$\delta_H\geq 1-3\cdot \frac{11}{3\cdot 2^7}
-4\cdot 0
-6\cdot \frac{11}{3\cdot 2^5}
=\frac{128-11-88}{2^7}>0$. 

\end{proof}

\begin{proposition}
\label{2FSL}

Assume $p=2$. Take a slim subgroup $H\subseteq \SL_2(\Z/2^{10}\Z)$. 
If $H/H_1$ is contained in 
$F$ or $H/H_1=\SL_2(\Z/2\Z)$, 
then $\delta_H>0$. 

\end{proposition}

\begin{proof}

In $\SL_2(\Z/2^{10}\Z)$ we have 
$\sharp\Conj(\sigma)=3\cdot 2^{17}$, 
$\sharp\Conj(\tau)=2^{19}$ 
and 
$\sharp\Conj(u^{2^i})=3\cdot 2^{16-2i}$ 
for $1\leq i\leq 6$ 
by Lemma \ref{sl2-conj}, \ref{sl2-upr}. 

Suppose $H/H_1\subseteq F$. 
By Lemma \ref{table}, we see that 
in $\SL_2(\Z/4\Z)$ we have 
$\sharp f_{2,1}^{-1}(F)
\cap \Conj(\sigma)=0$, 
$\sharp F
\cap \Conj(\tau)=2$, 
$\sharp f_{2,1}^{-1}(F)
\cap \Conj(u)=0$ and 
$\sharp f_{2,1}^{-1}(F)
\cap \Conj(u^2)=3$. 
We have 
$\sharp(H/H_3)\cap\Conj(\tau)\leq 
\sharp f_{3,1}^{-1}(F)\cap\Conj(\tau)
=2\cdot 2^4=32$ 
by Lemma \ref{V-str1}. 
By Proposition \ref{at2}, we have 
$\sharp H\cap\Conj(\tau)\leq 
a(\tau,2)_{10}+2^8(32-8)=13\cdot 2^{11}$.
Hence 
$\frac{\sharp H\cap\Conj(\tau)}{\sharp\Conj(\tau)}\leq 
\frac{13\cdot 2^{11}}{2^{19}}
=\frac{13}{2^8}$. 
The same calculation as in the previous lemma shows 
$$\frac{\sharp H\cap\Conj(u^{2^i})}{\sharp\Conj(u^{2^i})}
\leq \frac{b(u,2)_{10-i}+2^{10-i}}{3\cdot 2^{16-2i}}$$
for $1\leq i\leq 6$. 
Therefore 
\begin{align*}
\frac{\sharp I\backslash G/H}{[G:H]} &\leq 
\frac{2-1}{2}\cdot 0
+\sum_{i=1}^6 \frac{2-1}{2^{i+1}}\cdot 
\frac{b(u,2)_{10-i}+2^{10-i}}{3\cdot 2^{16-2i}}
+\frac{1}{2^7}\\
&=\frac{15+11+7+5+3+2+3}{3\cdot 2^7}
=\frac{23}{3\cdot 2^6}.
\end{align*}
Consequently, 
$\delta_H\geq 1-3\cdot 0
-4\cdot \frac{13}{2^8}
-6\cdot \frac{23}{3\cdot 2^6}
=\frac{64-13-46}{2^6}>0$.

Next suppose $H/H_1=\SL_2(\Z/2\Z)$. 
Since $H/H_2\subsetneq \SL_2(\Z/4\Z)$ by Lemma 2.1, 
it is conjugate to $A_1$ 
by Lemma \ref{A1-A2}. 
Lemma \ref{A1-conj} shows that 
in $\SL_2(\Z/4\Z)$ we have 
$\sharp A_1
\cap \Conj(\sigma)=3$, 
$\sharp A_1
\cap \Conj(\tau)=2$, 
$\sharp A_1
\cap \Conj(u)=0$ and 
$\sharp A_1
\cap \Conj(u^2)=0$. 
By Proposition \ref{as2}, we have 
$\sharp H\cap\Conj(\sigma)\leq 
a(\sigma,2)_{10}+2^8(3-2)=73\cdot 2^8$.
Hence 
$\frac{\sharp H\cap\Conj(\sigma)}{\sharp\Conj(\sigma)}\leq 
\frac{73\cdot 2^8}{3\cdot 2^{17}}
=\frac{73}{3\cdot 2^9}$. 
We have 
$\sharp(H/H_3)\cap\Conj(\tau)\leq 
\sharp f_{3,2}^{-1}(A_1)\cap\Conj(\tau)
=2\cdot 2^2=8$ 
by Lemma \ref{V-str1}. 
By Proposition \ref{at2}, we have 
$\sharp H\cap\Conj(\tau)\leq 
a(\tau,2)_{10}=5\cdot 2^{12}$.
Hence 
$\frac{\sharp H\cap\Conj(\tau)}{\sharp\Conj(\tau)}\leq 
\frac{5\cdot 2^{12}}{2^{19}}
=\frac{5}{2^7}$. 
The same calculation as in the case 
$H/H_1\subseteq F$ shows 
$$\frac{\sharp I\backslash G/H}{[G:H]}\leq 
\frac{2-1}{2}\cdot 0
+\frac{2-1}{2^2}\cdot 0
+\sum_{i=2}^6 \frac{2-1}{2^{i+1}}\cdot 
\frac{b(u,2)_{10-i}+2^{10-i}}{3\cdot 2^{16-2i}}
+\frac{1}{2^7}
=\frac{31}{3\cdot 2^7}.$$ 
Consequently, 
$\delta_H\geq 1-3\cdot \frac{73}{3\cdot 2^9}
-4\cdot \frac{5}{2^7}
-6\cdot \frac{31}{3\cdot 2^7}
=\frac{512-73-80-248}{2^9}>0$. 

\end{proof}

Proposition \ref{g_BCD}, Lemma \ref{17E} and 
Proposition \ref{19B} - \ref{2FSL} imply Theorem \ref{main}.

\end{document}